\documentclass[10pt,oneside]{book}
\usepackage{amsmath}\usepackage{amsfonts}
\usepackage[mathscr]{euscript}
\usepackage[all]{xy}
\xyoption{dvips}
\usepackage{XYADDS}
\usepackage{array}

\input ARTHEAD

\DeclareMathAlphabet      {\mathsfi}{OT1}{cmss}{m}{sl}
\DeclareMathAlphabet      {\mathsfb}{OT1}{cmss}{bx}{n}
\def\sfi{\ifmmode\mathsfi\else\sffamily\itshape\fi}
\def\sfb{\ifmmode\mathsfb\else\sffamily\bfseries\fi}
\def\scr{\mathcal}
\def\bbold{\mathbb}
\def\tsize{\displaystyle}
\def\labelstyle{\tsize}
\long\def\mpark#1{\relax}

\makeatletter
\def\skd{\mathsfi} 
\def\catname{\mathbf}
\def\tar{\mathop{\text{target}}}
\def\sr{\mathop{\text{source}}}
\def\eval{\mathop{\text{eval}}}

\def\I{\mathop{\bf I\fake}\nolimits}
\def\one{\mathop{\sf 1\fake}}

\let\prd\prod 

\def\fake{\hspace*{.01pt}} 

\def\@xiiipt{13.2}
 \DeclareMathSizes{\@xiiipt}{\@xiiipt}{8}{6}

\DeclareFontShape{OT1}{cmr}{bx}{n}{
      <5> cmbx5
      <6> cmbx6
      <7> cmbx7
      <8> cmbx8
      <9> cmbx9
      <10> <10.95> cmbx10
      <12><13.2> <14.4> <17.28> <20.74> <24.88>cmbx12
       }{}

\DeclareFontShape{OT1}{cmr}{m}{n}{
      <5> cmr5
      <6> cmr6
      <7> cmr7
      <8> cmr8
      <9> cmr9
      <10> <10.95> cmr10
      <12><13.2> <14.4> <17.28> <20.74> <24.88>cmr12
       }{}

\def\logopn#1{\expandafter\def\csname#1\endcsname{%
\mathop{\mathsf{#1\fake}}\nolimits}}
\def\logord#1{\expandafter\def\csname#1\endcsname{%
\mathord{\mathsf{#1\fake}}}}
\def\logvary#1{\expandafter\def\csname#1\endcsname{%
\mathop{\mathsfi{#1\fake}}\nolimits}}
\def\csspace#1{\expandafter\def\csname#1\endcsname{%
\mathop{\mathsfb{#1\fake}}\nolimits}}
\def\skname#1{\expandafter\def\csname#1\endcsname{%
\mathop{\mathsfi{#1\fake}}\nolimits}}
\def\frakop#1{\expandafter\def\csname#1\endcsname{%
\mathop{\mathfrak{#1\fake}}\nolimits}}

\newcommand{\Catsk}{\mathsfi{Cat}}

\csspace{Cat}
\csspace{BinProd}
\csspace{CCC}
\csspace{FinProd}
\csspace{FinLim}
\def\Eb{\mathop{\mathsfb E\fake}\nolimits}

\frakop{C}
\frakop{I}
\frakop{M}
\frakop{U}
\def\Ff{\mathop{{\frak F\fake}}\nolimits}

\logopn{AmbSp}
\logopn{Arrows}
\logopn{Arr}
\logopn{Ass}
\logopn{BsDiag}
\logopn{CatTh}
\logopn{CatStruc}
\logopn{Cod}
\logopn{Cones}
\logopn{Constants}
\logopn{ConstSp}
\logopn{Ddcn}
\logopn{Dia}
\logopn{Diagrams}
\logopn{Dom}
\logopn{Expr}
\logopn{Fillin}
\logopn{FLTh}
\logopn{FinLimTh}
\logopn{FinLimUnivMod}
\logopn{Form}
\logopn{FPTh}
\logopn{Func}
\logopn{Graph}
\logopn{Id}
\logopn{Incl}
\logopn{Init}
\logopn{InputTypes}
\logopn{Inp}
\logopn{Insert}
\logopn{Length}
\logopn{Left}
\logopn{LimCone}
\logopn{Lim}
\logopn{Lin}
\logopn{LinTh}
\logopn{LinUnivMod}
\logopn{List}
\logopn{Name}
\logopn{Node}
\logopn{Nodes}
\logopn{Oprns}
\logopn{Outp}
\logopn{Par}
\logopn{Pbl}
\logopn{Pbr}
\logopn{PEqAr}
\logopn{PPEqAr}
\logopn{Proj}
\logopn{Right}
\logopn{Rng}
\logopn{Sep}
\logopn{ShpGr}
\logopn{Sk}
\logopn{SynCat}
\logopn{Th}
\logopn{TypeList}
\logopn{TypeSet}
\logopn{Types}
\logopn{Type}
\logopn{UndCat}
\logopn{UndSk}
\logopn{UndGr}
\logopn{UnivMod}
\logopn{VarList}
\logopn{VarSet}
\logopn{Var}
\logopn{Vbl}
\logopn{Vertex}

\logvary{claim}
\logvary{claimcon}
\logvary{f}
\logvary{hyp}
\logvary{hypcon}
\logvary{verif}
\logvary{wksp}

\skname{E}
\skname{F}
\skname{S}
\skname{T}

\logord{arrow}
\logord{base}
\logord{bot}
\logord{comp}
\logord{cone}
\logord{csource}
\logord{ctarget}
\logord{curry}
\logord{ebot}
\logord{econe}
\logord{efid}
\logord{eq}
\logord{equ}
\logord{esoco}
\logord{etaco}
\logord{etop}
\logord{eufid}
\logord{efia}
\logord{ev}
\logord{exp}
\logord{inc}
\logord{fia}
\logord{fid}
\logord{fs}
\logord{fss}
\logord{id}
\logord{lam}
\logord{lass}
\logord{lfac}
\logord{lproj}
\logord{lunit}
\logord{mfac}
\logord{modpon}
\logord{n}
\logord{ob}
\logord{ppair}
\logord{prod}
\logord{proj}
\logord{rass}
\logord{rfac}
\logord{runit}
\logord{rproj}
\logord{s}
\logord{soco}
\logord{soob}
\logord{source}
\logord{target}
\logord{ta}
\logord{ter}
\logord{taco}
\logord{top}
\logord{tsource}
\logord{ttarget}
\logord{twovf}
\logord{unit}
\logord{ufid}
\logord{vertex}
\def\vsfi{\mathord{{\mathsfi v}}}
\def\Fa{\Name[\F]}

\def\arr{\mathord{\mathsf{ar\fake}}}

\def\@listI{\leftmargin\leftmargini \parsep 5pt plus 2.5pt minus 1pt\topsep
10pt plus 4pt minus 6pt\itemsep 3pt plus 2pt minus 1pt}
\let\@listi\@listI
\@listi
\makeatletter
\def\biggg{\bBigg@4}
\def\Biggg{\bBigg@{4.5}}
\def\bigggg{\bBigg@{5.0}}
\def\Biggg{\bBigg@{5}}
\def\Bigggg{\bBigg@{7}}
\def\biggggg{\bBigg@{8}}
\def\Biggggg{\bBigg@{9}}
\makeatother

\def\brspace{\hspace{3em}}

\makeatother
\def\coneto{\hspace{-1pt}\mathrel{%
\raisebox{1ex}%
{\xymatrix@R=.5ex@C=1em@M=0pt@W=0pt{ &  \ar@{-}[dd] \\
\ar@{-}[ur]\ar@{-}[dr] & \\
&  \\
}}}}\makeatletter

\renewcommand\chapter{\if@openright\cleardoublepage\else\clearpage\fi
                    \thispagestyle{plain}%
                    \global\@topnum\z@
                    \@afterindentfalse
                    \secdef\@chapter\@schapter}
\def\@chapter[#1]#2{\ifnum \c@secnumdepth >\m@ne
                       \if@mainmatter
                         \refstepcounter{chapter}%
                         \typeout{\@chapapp\space\thechapter.}%
                         \addcontentsline{toc}{chapter}%
                                   {\protect\numberline{\thechapter}#1}%
                       \else
                         \addcontentsline{toc}{chapter}{#1}%
                       \fi
                    \else
                      \addcontentsline{toc}{chapter}{#1}%
                    \fi
                    \chaptermark{#1}%
                    \addtocontents{lof}{\protect\addvspace{10\p@}}%
                    \addtocontents{lot}{\protect\addvspace{10\p@}}%
                    \if@twocolumn
                      \@topnewpage[\@makechapterhead{#2}]%
                    \else
                      \@makechapterhead{#2}%
                      \@afterheading
                    \fi}
\def\@makechapterhead#1{%
  {\parindent \z@ \raggedright \normalfont
    \ifnum \c@secnumdepth >\m@ne
      \if@mainmatter
        \Large\bfseries \@chapapp\space \thechapter
        \par\nobreak
         \vskip 1ex
      \fi
    \fi
    \interlinepenalty\@M
    \LARGE \bfseries #1\par\nobreak
    \vskip 2ex
  }}
\def\@schapter#1{\if@twocolumn
                   \@topnewpage[\@makeschapterhead{#1}]%
                 \else
                   \@makeschapterhead{#1}%
                   \@afterheading
                 \fi}
\def\@makeschapterhead#1{%
  {\parindent \z@ \raggedright
    \normalfont
    \interlinepenalty\@M
    \LARGE \bfseries  #1\par\nobreak
    \vskip 2ex
  }}

\def\@seccntformat#1{\csname the#1\endcsname\hspace*{.6em}}

\def\section{\@startsection {section}{1}{\z@}{-3ex plus -1ex minus
 -.2ex}{1ex plus .2ex}{\Large\bf}}
\def\subsection{\@startsection{subsection}{2}{\z@}%
{-1.5ex plus -.3ex minus -.1ex}{.2pt}{\normalsize\bf}}
\def\subsubsection{\@startsection{subsubsection}{3}{\z@}%
{-2ex plus -.1ex minus -.1ex}{-1em}{\normalsize\bf}}
\def\specialsubsection{\@startsection{subsection}{3}{\z@}%
{-1.5ex plus -.1ex minus -.1ex}{-1em}{\normalsize\bf}}

\def\@listI{\leftmargin\leftmargini \parsep 5pt plus 2.5pt minus 1pt\topsep
10pt plus 4pt minus 6pt\itemsep 3pt plus 2pt minus 1pt}
\let\@listi\@listI
\@listi

\makeatother

\renewenvironment{squishlist}
{\begin{list}{\arabic{lister})}%
{\usecounter{lister}\parsep0pt\itemsep 0pt plus 2pt
\listparindent 1.5em}}
{\end{list}}
\renewcommand{\sqlist}{\begin{squishlist}}
\renewcommand{\esqlist}{\end{squishlist}\noindent}

\mathchardef\eqs="603D
\mathchardef\lra="622C
\mathchardef\xord="0202 
\let\xo\xord
\def\Frac#1#2{{\def\arraystretch{1.2}
\begin{array}{c}{#1}\\[6pt] \hline
\\[-12pt]{#2}\end{array}}}

\def\Rule#1#2#3#4{\mbox{\rm #3}%
\hspace*{1.5em}{\Frac{\kern-2pt #1}%
{\kern2pt #2}}
\hspace*{1.5em} \parbox{2in}{\raggedright #4}}

\def\artwo#1#2{\ar@<1ex>[r]^{#1}\ar@<-1ex>[r]_{#2}}
\let\rtwo\artwo

\def\epf{$\Box$\par\addvspace{\medskipamount}}

\def\ceq{\mathrel{\colon\hspace{-2pt}\eqs}}
\def\ceqv{\mathrel{\colon\hspace{-2pt}\lra}}
\def\astep{\subsubsection*{Translation as a construction}}
\def\bstep{\subsubsection*{Expression as actual factorization}}

\def\mk{\underline}
\def\ds#1#2#3{#1^{#2}_{#3}}
\def\dsu#1#2#3{\mk{#1}^{#2}_{#3}}

\def\crs#1#2{\prd_{i=1}^{#2}{#1}_i}
\def\typevar{\gamma}
\def\ab#1{\Length[#1]}

\def\narrower{\spreaddiagramcolumns{-.5pc}}
\def\wider{\spreaddiagramcolumns{.5pc}}
\def\muchwider{\spreaddiagramcolumns{2.5pc}}
\def\taller{\spreaddiagramrows{.5pc}}
\def\shorter{\spreaddiagramrows{-.5pc}}
\def\Narrower{\spreaddiagramcolumns{-2pc}}
\def\Wider{\spreaddiagramcolumns{2pc}}

\def\pf{\specialsubsection*{Proof}}
\makeatletter
\def\mld#1$${\null\,\vcenter\bgroup\def\\{\cr&}\openup9pt\m@th
\ialign\bgroup\strut\hfil$\displaystyle{##}$%
&$\displaystyle{{}##}$\hfil\crcr #1
\crcr\egroup\egroup\,$$}
\makeatother

\def\inv{^{-1}}
\expandafter\ifx\csname amssym.def\endcsname\relax \else \fi
\expandafter\edef\csname amssym.def\endcsname{%
\catcode`\noexpand\@=\the\catcode`\@\space}
\catcode`\@=11
\def\undefine#1{\let#1\undefined}
\def\newsymbol#1#2#3#4#5{\let\next@\relax
 \ifnum#2=\@ne\let\next@\msafam@\else
 \ifnum#2=\tw@\let\next@\msbfam@\fi\fi
 \mathchardef#1="#3\next@#4#5}
\def\mathhexbox@#1#2#3{\relax
 \ifmmode\mathpalette{}{\m@th\mathchar"#1#2#3}%
 \else\leavevmode\hbox{$\m@th\mathchar"#1#2#3$}\fi}
\def\hexnumber@#1{\ifcase#1 0\or 1\or 2\or 3\or 4\or 5\or 6\or 7\or 8\or
 9\or A\or B\or C\or D\or E\or F\fi}

\def\ds#1#2#3{#1^{#2}_{#3}}
\def\dsu#1#2#3{\underline{#1}^{#2}_{#3}}

\def\scr{\mathcal}
\def\Asc{\mathord{\scr A}}
\def\Bsc{\mathord{\scr B}}
\def\Csc{\mathord{\scr C}}

\def\Fsc{\mathord{\scr F}}

\def\Psc{\mathord{\scr P}}

\def\Ssc{\mathord{\scr S}}
\def\Tsc{\mathord{\scr T}}
\def\Usc{\mathord{\scr U}}

\CompileMatrices

\DeclareFontFamily{U}{euf}{}
\DeclareFontShape{U}{euf}{m}{n}{
  <-6> eufm5 <6-8> eufm7 <8-> eufm10
   }{} \DeclareFontShape{U}{euf}{b}{n}{
  <-6> eufb5 <6-8> eufb7 <8-> eufb10   }{}
\def\todo#1{\textbf{[To do: \emph{#1}]}\par\vspace{6pt}}
\long\def\ig#1{\relax}
\let\mpar\ig
\newcommand{\oc}%
{\raise.2ex\hbox{$\scriptstyle\mskip\tinymuskip\circ\mskip\tinymuskip$}}
\begin{document}
\title{\textbf{Graph-based Logic\\ and Sketches}\emph{\\Preliminary version: Please refer others to this website instead of giving them copies of this version.}}

\author{Atish Bagchi\thanks{The first author is
grateful to the Mathematical Sciences Research Institute
for providing the opportunity to work on this book there.
Research at MSRI was supported in part by NSF grant
\#DMS-9022140.}\ \ and Charles Wells}

\maketitle

\tableofcontents

\chapter*{Preface}
\addcontentsline{toc}{chapter}{Preface}
\textbf{\emph{This is a preliminary version of this monograph.  In certain places changes and additions to be made for the final version are marked in a paragraph beginning \lq To do\rq.}}
\vspace*{6pt}

This monograph
presents the basic idea of forms (a generalization of
Ehresmann's sketches) and their models. We also provide a
formal logic that gives an intrinsically categorial
definition of the concepts of assertion and proof for any
particular type of form.  We provide detailed examples of
the machinery that enables the construction of forms, as
well as examples of proofs in our formalism of certain
specific assertions.

Our monograph requires familiarity with the basic notions of
mathematical logic as in Chapters 2 through 5 of
\cite{EFT}\mpark{Do the reader really need as much logic as
it says here?}, and with category theory and sketches as in
Chapters 1 through 10 of \cite{ctcs}. We specifically
presuppose that finite-limit sketches and their models are known. Some
notation for these ideas is established in
Section~\ref{oldskdefs}.

This work is a combination and revision of the work in \cite{logstr} and \cite{logstr2}.\mpar{Revised paragraph.}  It
is better because of conversations we have
had with Robin Cockett and Colin McLarty. We are grateful
to Frank Piessens and the referees for careful readings of
earlier versions
that uncovered errors,
and to Max Kelly, Anders Kock and Steve Lack for supplying
references. The names ``string-based logic'' and
``graph-based logic'' were suggested by Peter Freyd. The
diagrams were prepared using K. Rose's {\tt xypic}.

{\flushright Atish Bagchi\\ Charles Wells\\}

\mainmatter

\chapter{Introduction}\mpark{I
have rewritten the introduction because of the merger.  We
should go over it in detail again later.}\mpar{All through the text I have changed "section" to "Section" when it is followed by a section number for consistency.}

\section{Brief outline}

Sketches as a method of specification of mathematical
structures are an alternative to the string-based
specification employed in mathematical logic. They have
been proposed as suitable for the specification of data
types and programs \cite{gray87}, \cite{gray89},
\cite{wellsbarr}, \cite{duvsen}, and some work on
implementation has been carried out \cite{gray90},
\cite{duvrey1}, \cite{duvrey2}.

In \cite{gensk} the second author introduced the notion of
\textbf{form}, a graph-based method of specification of
mathematical structures that generalizes Ehresmann's
sketches. Forms are a proper generalization of sketches: a
form can have a model category that cannot be the model
category of a sketch (Section~\ref{modthform}). Forms were
generalized to $2$-categories by Power and Wells
\shortcite{powerwells} (where the word ``sketch'' was used
instead of ``form''). Other generalizations of sketches are
given in \cite{lairtrames}, \cite{msdp} and \cite{kpt}.

Sketch theory has been criticized as being lacunary when
contrasted with logic because it apparently has nothing
corresponding to proof theory. In this monograph, we
outline\mpar{Changed "provide" to "outline"} a uniform proof theory for all types of sketches
and forms.  We show that, in the case of finite-product
sketches, this results in a system with the same power as
equational logic.

\section{Types of Logic}

\subsection{String-based logic}\label{sblogic} Traditional
treatments of formal logic provide:
\begin{labeledlist}{SBL}
\item \label{form} A syntax for terms and formulas.  The
formulas are typically strings of symbols defined
recursively by a production system (grammar), and the
complete syntax of each term or formula is provided by the
corresponding parsing tree (formation tree).  To deduce the tree from the
string of symbols requires fairly sophisticated pattern
matching by the reader, or else a parsing mechanism in a
computer program. \item\label{infrel} Inference relations
between sets of formulas. This may be given by structural
induction on the construction of the formulas, so that to
verify an inference relation requires understanding the
parse trees of the formulas involved. \item\label{sem}
Rules for assigning meaning to formulas (semantics) that
are sound with respect to the inference relation. The
semantics may also be given by structural induction on the
construction of the formulas.
\end{labeledlist}

First order logic, the logic and semantics of programming
languages, and the languages that have been formulated for
various kinds of categories are all commonly described by SBL 1--3.

The  strings of symbols that constitute the terms and
formulas are ultimately modeled on the sentences
mathematicians speak and write when proving theorems.  Thus
first order logic in particular is a \textit{mathematical
model} of the way in which mathematicians reason.  The terms and
formulas are ordinary mathematical objects.

\subsection{Graph-based logic} Mathematicians with a
category-theoretic point of view frequently state and prove
theorems using graphs and diagrams (described in Chapter~\ref{prelims}). The graphs, diagrams,
cones and other data of a sketch or form are formal objects
that correspond to the graphs and diagrams used by such
mathematicians in much the same  way in which the formulas of
traditional logic correspond to the sentences
mathematicians use in proofs.

The  functorial semantics of sketches and forms corresponds
to item~SBL.\ref{sem} in the list in Section~\ref{sblogic}.
This semantics is sound in the informal sense that it
preserves by definition the structures given in the sketch
or form. The analogy to the semantics of traditional model
theory is close enough that sketches and forms and their
models fit the definition of ``institution''
(\cite{gogburst}), which is an abstract notion of a logical
system having syntactic and semantic components. This is
described in detail for the case of sketches in
\cite{ctcs}, Section 10.3. Note that the soundness of
functorial semantics appears trivial when contrasted with
the inductive proofs of soundness that occur in
string-based logic because the semantics functor is not
defined recursively.

This monograph exhibits a structure in the theory of
sketches and forms that corresponds to items~SBL.\ref{form}
and~SBL.\ref{infrel}. The data making up the structure we
give do not correspond in any simple way to the data
involved in items~SBL.\ref{form} and~SBL.\ref{infrel} of
traditional logic; we discuss the relationship in
Section~\ref{ruled}.

\section{Forms}\label{formssec}\mpar{Revised 6 Oct 07.} Forms are
parametrized by the type of constructions they allow. Let $\Catsk$ be the finite-limit sketch for categories
(described in detail in Section~\ref{realcatsk}). Let
$\E$ denote a finite-limit sketch whose models are a type of category with structure that is essentially algebraic over categories (this means that $\E$
contains  $\Catsk$ as a subsketch with the property that all the
objects of $\E$ are limits of finite diagrams in $\Catsk$). Let $\Eb$ denote the finite-limit theory
generated by $\E$. $\E$ is called a \textbf{constructor space sketch} and $\Eb$ is a \textbf{constructor space}.

$\Eb$-categories are then the models of
$\Eb$ in the category of sets. (This is described in detail
in Chapter~\ref{skgeneral}.) The kinds of categories that
can be described in this way include categories with
specified finite products, categories with specified limits
or colimits over any particular set of diagrams, cartesian
closed categories, regular categories, toposes, and many
others (always with specified structure rather than
structure determined only up to isomorphism --- see
Chapter~\ref{futwor}). Finite limit sketches for several
specific instances of $\Eb$ are given in
Chapter~\ref{conspsk}.

An $\Eb$-form $\F$  is a graph-based structure that
allows the specification of any kind of
construction that can be made in any
$\Eb$-category.  A model of $\F$ in an
$\Eb$-category $\Csc$ is, informally, an instance
of that construction in $\Csc$ (see the remarks
in~\ref{modrem1}). Forms are defined precisely in
Section~\ref{skedef}. Ordinary sketches can also be
realized as forms.

As an example, let $\CCC$ be a finite-limit theory for
Cartesian closed categories (one is outlined in
Section~\ref{cccsec}). It is possible to require that
a certain object in a $\CCC$-form $\F$ be the ``formal
function space'' $A^B$ of two other objects $A$ and
$B$ of the form. This means that in any model ${\frak
M}$ of the form in a Cartesian closed category $\Csc$, the value of ${\frak
M}(A^B)$ is in fact the function space ${\frak
M}(A)^{{\frak M}(B)}$ in $\Csc$. The object $A^B$ is
not itself a function space; it is an object of a form
(that is, a generalized sketch), not of a Cartesian
closed category. That is why it is called a {\it
formal\/} function space.

Forms have much more expressive power than sketches
as originally defined by Ehresmann, in which only
limits and colimits can be specified.

An $\Eb$-form $\F$ is determined by a freely
adjoined global element $\Fa$ of the limit vertex
of a diagram in $\Eb$, obtaining a category
$\SynCat[\Eb,\F]$.  More details are in
Section~\ref{formdes}. An assertion in this setting
is a potential factorization (PF) (defined precisely in
Section~\ref{pfsec}) of an arrow of
the $\SynCat[\Eb,\F]$ through an arrow into its
codomain.  The assertion is valid if
the PF does indeed factorize in every model of
$\SynCat[\Eb,\F]$.

Instead of the set of rules of deduction of a
traditional theory, we have a set of rules of
construction.  More precisely, we give in
Section~\ref{rulesapp} a system of construction rules
that produce all the objects and arrows of the
categorial theory of a finite-limits sketch. These
rules apply in particular to $\SynCat[\Eb,\F]$, which
is constructed in~\ref{sec613} as such a categorial
theory.

We say that the potential factorization is {\bf
deducible} if there is an actual factorization in
$\SynCat[\Eb,\F]$. Such an arrow must be constructible
by the rules in Section~\ref{rulesapp}. Thus the usual
system of inference is replaced by a system of
construction of arrows in the  finite-limit category
$\SynCat[\Eb,\F]$ (no matter what type of category is
sketched by $\E$). This system is sound and complete
with respect to models (Section~\ref{compsub}).

The fact that we have assumed finite-limit sketches as
given prior to the general definition of $\Eb$-form is
basic to the strategy of this monograph, which is to make
finite-limit logic the  logic for all forms
(Section~\ref{dedrules}) -- but \textit{not} the logic of
the model category. The variation in what can be proved,
for example for finite-product forms ($\Eb$-forms where
$\Eb=\FinProd$ as in Section~\ref{finprodsk}) as contrasted
with cartesian closed forms ($\Eb$ is $\CCC$ as in
Section~\ref{cccsec}) is entirely expressed by the choice
of $\Eb$ and has no effect on the rules of construction.

In \cite{gogmes}, Goguen and Meseguer produced a sound and
complete entailment system for multisorted equational
logic. In Chapter~\ref{EqThChap}, we verify that the
theorems of that logic for a particular signature and
equations all occur as actual factorizations in
$\SynCat[\FinProd,\F]$, where $\F$ is a $\FinProd$ form
induced (in a manner to be described) by the given
signature and equations.  We also compare the expressive
powers of these two systems.

\todo{The converse construction: Show how to create a multisorted equational theory with the same models as a $\FinProd$ form, and use the completeness theorem in Section~\ref{compsub} to show that every actual factorization of the form arises from a deduction in the theory.}

\section{Glossary}

This monograph introduces a large number of structures with
confusingly similar roles. We list the most important
here with a reference to the section in which they are
defined.

\begin{itemize}

\item $\LinTh[L]$, \ref{linsk}.
\item $\FinLimTh[\S]$ , \ref{catthsec}.
\item $\Eb$ ($=\FinLimTh[\E]$), \ref{cspaces}.
\item  $\SynCat[\Eb,\F]$, \ref{formdes}.
\item
$\CatTh[\Eb,\F]$, \ref{catthform}.

\end{itemize}

\chapter{Preliminaries}\label{prelims}

\section{Lists}
Given a set $A$, $\List[A]$ denotes the set of lists of
elements of $A$, including the empty list.  The $k$th entry
in a list $w$ of elements of $A$ is denoted by $w_k$ and
the length of $w$ is denoted by $\ab{w}$. If $f:A\to B$ is
a function, $\List[f]:\List[A]\to\List[B]$ is by definition
$f$ ``mapped over'' $\List[A]$:  If $w$ is  a list of
elements of $A$, then the $k$th entry of $\List[f](w)$ is
by definition $f(w_k)$. This makes $\List$ a functor from
the category of sets to itself.

\section{Graphs}

\defn A \textbf{graph}\mpar{Whole section rewritten and expanded.} $G$ is a mathematical structure consisting of a set $\Nodes[G]$ of \textbf{nodes} and a set $\Arrows[G]$ of \textbf{arrows}, with two functions $\source:\Arrows[G]\to\Nodes[G]$ and $\target:\Arrows[G]\to\Nodes[G]$.\edefn
Graphs may be pictured by drawing dots for nodes and an arrow $a$ from a node $m$ to a node $n$ whenever $\source(a)=m$ and $\target(a)=n$.  These are what category theorists customarily call ``graphs''.  In the graph theory literature they would be called ``directed multigraphs with loops.''

The
underlying graph of a category $\Csc$ is denoted by
$\UndGr[\Csc]$. A subgraph $H$ of a graph $G$ is said to be
{\bf full} if every arrow $f:h_1\to h_2$ of $G$ between
nodes of $H$ is an arrow of $H$.

\rem
Graphs are conceptually more primitive than the well-formed formulas used in string-based logic.  Graphs are given by a linear sketch (see Section~\ref{linsk}), essentially the simplest form of sketch, whereas wff's must be given by a context free grammar (recursive definition) which is equivalent to a finite-product sketch. (See \cite{ctcs}, page 235.)
\erem

\section{Diagrams}\label{diagsec}

\defn Two graph homomorphisms $\delta:I\to G$ and $\delta':I'\to G$ are
said to be \textbf{equivalent} if there is a graph isomorphism
$\phi:I\to I'$ such that

\begin{equation}\label{diagdiag}
\xymatrix{
I \ar[rr]^\phi \ar[dr]_\delta &&
I'\ar[dl]^{\delta'} \\ & G \\ } \end{equation}
commutes.\edefn This
relation is easily seen to be an equivalence
relation on the set of graph homomorphisms into a
graph $G$.

\defn A {\bf diagram} in $G$ is by definition
an equivalence class of
graph homomorphisms $\delta:I\to G$.\edefn

As is the practice when an object is defined to be
an equivalence class, we will refer to a diagram by
any member of the equivalence class.

\defn If $\delta:I\to G$ is a diagram, $I$ is said
to be a {\bf shape graph} of the diagram, denoted
by $\ShpGr[\delta]$, and $G$ is said to be the {\bf
ambient space} of the diagram, denoted by
$\AmbSp[\delta]$.\edefn

Observe that the ambient space of the diagram is
determined absolutely, but the shape graph is
determined only up to an isomorphism that makes
Diagram~(\ref{diagdiag}) commute.

\defn Let $I$ be  a graph and  $\Csc$  a category.
To say that $\delta$ is a {\bf diagram in $\Csc$}
means that $\delta:I\to \UndGr[\Csc]$ is a diagram.
We write $\delta:I\to\Csc$ to denote this
situation.\edefn

Note that $\Csc$ is part of the definition: there
could be another category with the same underlying
graph.

\defn A diagram $\delta:I\to\Csc$ in a category
$\Csc$ {\bf commutes} if whenever
$(f_1,\ldots,f_n)$ and $(g_1,\ldots,g_m)$ are two
paths in $I$ with the same source and target, then
$$\delta(g_m)\oc \delta(g_{m-1})\oc\cdots \oc
\delta(g_1)=\delta(f_n)\oc \delta(f_{n-1})\oc\cdots\oc
\delta(f_1)$$ in $\Csc$.\edefn

Observe that commuting is defined only for diagrams
in a category. More details about this may be found
in \cite{ctcs}, Section 4.1.5.

\section{Convention on drawing diagrams}\label{convention}
It is customary to draw a diagram without naming its shape
graph.  We adopt the following convention:  If a diagram is
represented by a drawing, the shape graph of the diagram is
a graph that has one node for each object shown in the
drawing and one arrow for each arrow shown, with source and
target as shown. Two objects at different locations in the
drawing correspond to two different nodes of the shape
graph, {\it even if the objects have the same label}, and
an analogous remark applies to arrows.  Thus\mpark{added sentence} the traditional presentation of a graph, as in (\ref{d1ex}), reveals the equivalence class of the graph but not precisely which shape graph is used (which is irrelevant in any case).

\exam

The diagram~(\ref{d1ex})

\begin{equation}\label{d1ex}
\xymatrix{
A \ar[r]^f \ar[d]_g & A
\ar[d]^g \\
B \ar[r]_{\Id[B]} & B
}
\end{equation}
called $\delta$, has shape  graph


\begin{equation}\label{shape1}
\xymatrix{
h \ar[r]^t \ar[d]_u & i
\ar[d]^v \\
j \ar[r]_x & k
}
\end{equation}
so that $\delta(h)=A$, $\delta(i)=A$,
$\delta(v)=g$, $\delta(x)=\Id[B]$ and so on.
Diagram~(\ref{ABCD}) below also has shape
graph~(\ref{shape1}) (or one isomorphic to it, of course):


\begin{equation}\label{ABCD}
\xymatrix{
A \ar[r]^f \ar[d]_g & B
\ar[d]^h \\
C \ar[r]_k & D
}
\end{equation}

On the other hand, Diagram~(\ref{AAB}) below

\begin{equation}\label{AAB}
\xymatrix{
A \ar[rr]^f \ar[dr]_g &&
A\ar[dl]^g \\ & B \\ } \end{equation}

{\it is not the same diagram as~(\ref{d1ex})}.
It has shape
graph

\begin{equation} \label{p7un1}
\xymatrix{
i \ar[rr]^u \ar[dr]_v &&
j\ar[dl]^w \\ & k \\ } \end{equation}

\rem The reader should observe that we use ``diagram'' and
``commutative diagram'' both in the object language and the
metalanguage. For example, in Section~\ref{diagsec} we
refer to Diagram~(\ref{diagdiag}), which must commute in
the category of graphs. Note that we did not mention its
shape graph, but according to the principles just
enunciated that shape graph must be (isomorphic to)

\begin{equation}\label{p7un1}
\xymatrix{ i \ar[rr]^u \ar[dr]_v && j\ar[dl]^w \\ & k \\ }
\end{equation}\erem

\rem Diagrams\label{diagrem} are customarily drawn as
planar or nearly planar graphs or as perspective
representations of three-dimensional graphs.  A well-drawn
graph reveals a lot of information quite efficiently to
human beings and at the same time  shows more of the
structure than the formulas of traditional logic commonly
do.  Nevertheless, the details of the representation
(nearness to planarity, symmetry when possible, and so on)
that aid in human understanding are not part of the
abstract structure of the diagram at all.\erem

\section{Cones}\label{conesec} For any graph $G$
and diagram $\delta:I\to G$, a cone
$\Theta:v\coneto(\delta:I\to G)$ (also written
$\Theta:v\coneto\delta$ if the context makes this
clear) has {\bf vertex} $v$ denoted by
$\Vertex[\Theta]$ and {\bf base diagram} $\delta$
denoted by $\BsDiag[\Theta]$. For each node $i$ of
$\ShpGr[\delta]$, the formal projection of the cone
$\Theta$ from $\Vertex[\Theta]$ to $\delta(i)$ is
denoted by $\Proj[\Theta,i]:v\to \delta(i)$. For a
category $\Csc$, a cone
$\Theta:v\coneto(\delta:I\to\Csc)$ is {\bf
commutative} if, for every arrow $f:i\to j$ in $I$,
the diagram


\begin{equation}\label{d5ex}
\xymatrix{
&\Vertex[\Theta]\ar[dl]_{\Proj[\Theta,i]}
\ar[dr]^{\Proj[\Theta,j]}
\\
\delta(i)\ar[rr]_{\delta(f)}&&\delta(j)
\\}
\end{equation}
commutes.

In the following, we are concerned with categories
with specified finite limits.  In such categories,
the specified limit cone of a diagram $\delta$ will
be denoted by
$\LimCone[\delta]:\Lim[\delta]\coneto\delta$.  This
specifically applies to Rule $\exists$LIM
of~\ref{arrrul}.

\section{Fonts}\label{fontsec}
In general, variable objects are given in slant or
script notation and specific objects (given by proper names) are
given in upright notation.  In more detail, we have the
following notational scheme.

\begin{enumerate}

\item Specific data constructors, such as $\List$,
and specific fieldnames for complex objects, such
as $\Nodes[G]$, are given in {\sf sans serif} and
are capitalized as shown.

\item Specific objects and arrows of sketches or
forms are also given in {\sf sans serif}.

\item Specific constructor spaces, such as
$\FinLim$ and $\CCC$, are given in {\sfb bold sans
serif}. We use $\Eb$ to denote a variable
constructor space because of the unavailability of
bold slanted sans serif.

\item Specific categories other than constructor
spaces, such as $\Set$, are given in {\bf
boldface}.

\item Diagrams (specific and variable) are named by
lowercase Greek letters, such as $\delta$ and $\gamma$.

\item Cones (specific and variable) are named by
uppercase Greek letters, such as $\Phi$ and $\Psi$.

\item Models (specific and variable) are given in
uppercase fraktur, for example $\frak M$, $\frak
C$.

\item Variable sketches and forms
are given in {\sfi slanted sans serif}.

\item Variable categories other than constructor
spaces are given in script, for example $\Asc$,
$\Bsc$, $\Csc$.

\item Other variable objects are given in math
italics, such as $a$, $b$, $c$ or (especially
arrows) in lowercase Greek letters. \end{enumerate}

\chapter{Sketches}\label{oldskdefs}

We use a general concept of form described in
Chapter~\ref{skgeneral} that is based on the concept of
finite limit sketch, a particular case of projective sketch
due to Ehresmann.  In this section, we review briefly some
aspects of linear and finite-limit sketches that are
relevant later.

\section{Linear sketches}\label{linsk}
\defn A {\bf linear sketch} $L$ is a
pair $\left(\Graph[L], \Diagrams[L]\right)$ where $\Graph[L]$ is
a graph and $\Diagrams[L]$ is a set of finite diagrams in
$\Graph[L]$.\edefn
\defn
A {\bf model}\mpar{Added definition of model} $\M$ of a linear sketch $\S$ in a
category $\Csc$ is a graph homomorphism
$\M:\Graph[\S]\to\UndGr[\Csc]$ that takes the diagrams in
$\Diagrams[\S]$ to commutative diagrams in $\Csc$.
A {\bf morphism of models} is a natural transformation.\edefn
\defn The {\bf linear theory} generated by a linear sketch $L$ is
the category obtained from the free category generated by
$\Graph[L]$ by imposing the least congruence relation that
makes the diagrams in $\Diagrams[L]$ commute. Here we call the
linear theory $\LinTh[L]$.\edefn

\defn The {\bf universal model} of a linear sketch $L$, denoted
by $\LinUnivMod[L]:L\to\LinTh[L]$,
is the morphism of sketches whose underlying morphism is the
induced graph homomorphism (quotient map) from $\Graph[L]$ to
the graph $\UndGr\Bigl[\LinTh[L]\,\Bigr]$.\edefn
\rem Linear\mpar{Added remark} sketches are so called because the underlying functor from a model in $\Set$ preserves products and coproducts. \erem
\section{Finite-limit sketches}
\defn A {\bf finite-limit sketch}
$\S$ is a triple
$$\left(\Graph[\S],\Diagrams[\S],\Cones[\S]\right)$$ where
$\Graph[\S]$ is a graph, $\Diagrams[\S]$ is a set of finite
diagrams in $\Graph[\S]$, and $\Cones[\S]$ is a set of cones in
$\Graph[\S]$, each to a finite diagram in $\Graph[\S]$ (which
need not be in $\Diagrams[\S]$).\edefn
\defn For finite-limit sketches $S$ and $S'$, a {\bf sketch
morphism} ${\frak m}:\S\to\S'$ is a graph homomorphism ${\frak
m}:\Graph[\S]\to\Graph[\S']$ that takes the diagrams in
$\Diagrams[\S]$ to diagrams in $\Diagrams[\S']$ and the cones in $\Cones[\S]$
to cones in $\Cones[\S']$.\edefn
\defn A {\bf model} $\M$ of a finite-limit sketch $\S$ in a
category $\Csc$ is a graph homomorphism
$\M:\Graph[\S]\to\UndGr[\Csc]$ that takes the diagrams in
$\Diagrams[\S]$ to commutative diagrams in $\Csc$ and the cones
in $\Cones[\S]$ to limit cones in $\Csc$. A {\bf morphism of models} is a natural transformation.\mpar{Added sentence}\edefn

\defn The {\bf forgetful functor} $\UndSk$ from the category of
small categories with finite limits and finite-limit preserving
functors to the category of finite-limit sketches takes a
category $\Csc$ to $(\UndGr[\Csc],D,L)$ where $D$ is the set of
all finite commutative diagrams in $\Csc$ and $L$ is the set of
all limit cones in $\Csc$ to finite diagrams in $\Csc$.\edefn

\section{Categorial theories of finite-limit
sketches}\label{catthsec}
\defn Let $S$ be a finite-limit sketch.
The {\bf finite-limit theory} generated by $\S$,
denoted by $\FinLimTh[\S]$,
is a category with finite limits together with a model
$$\FinLimUnivMod[\S]:\S\to\FinLimTh[\S]$$ called the {\bf universal
model} of the sketch.  It has the following property: For every
model $\M$ of $\S$, there is a finite-limit
preserving
functor $\FinLimTh(\M):\FinLimTh[\S]\to\Set$,
determined uniquely up to natural isomorphism, with the property
that

{\muchwider\muchwider\muchwider\muchwider\taller
\taller\taller\taller
\begin{equation}\label{thiss}
\xymatrix{\S\ar[r]^{\FinLimUnivMod[\S]}\ar[dr]_{\M}
&\UndSk[\FinLimTh[S]]\ar[d]^{\UndSk[\FinLimTh[\M]]}\\
&\UndSk[\Set]}
\end{equation}
}
commutes.\edefn

\rem It follows from the defining
properties that $\FinLimTh[\S]$ is determined up to
equivalence of categories and $\FinLimUnivMod[\S]$ is
determined up to natural isomorphism.\erem

\section{Remark concerning models}\label{modrem1}
Category theorists commonly use the name of a category
to refer indifferently to any equivalent category.  A
related phenomenon occurs with respect to models, as
we discuss. This monograph is concerned with syntax,
and it may be necessary for clarity if not for strict
correctness to distinguish between mathematical
constructions that would be regarded as the ``same''
by many category theorists.

Let $\Catsk$ be the finite-limit sketch for categories
described in detail in Section~\ref{realcatsk}. We may consider the
following three mathematical
entities.\begin{enumerate} \item Some small category
$\Csc$ of one's choice. \item A model $\C$ of the
sketch $\Catsk$ in the category of sets which ``is''
or ``represents'' $\Csc$.  This means that $\C$ is a
morphism of finite-limit sketches from $\Catsk$ to an
underlying finite-limit sketch of $\Set$ with the
property that $\C(\ob)$ is the set of objects of
$\Csc$, $\C(\arr)$ is the set of arrows of $\Csc$,
$\C(\comp)$ is the composition function of $\Csc$ (up
to natural isomorphism), and so on. $\C$ is determined
uniquely by $\Csc$ if one assumes that $\Set$ has
specified finite limits.  It is determined up to
natural isomorphism in any case. \item The model
$\FinLimTh[\C]$ of $\FinLimTh[\Catsk]$ induced by $\C$
according to Section~\ref{catthsec}. This is a
finite-limit-preserving functor from $\FinLimTh[\Catsk]$
to $\Set$.\end{enumerate}

For many category
theorists, $\Csc$ and $\C$ denote the ``same thing''.
Other mathematicians would
disagree, saying $\Csc$ is a category (presumably for them ``category'' has
some meaning other than ``model for the sketch for categories'') and
$\C$ is a morphism of sketches, so how could they be the same thing?
This\mpar{Here you wanted an elaborate explanation of the meaning of $(\Z,+)$.  I don't under\-stand why.  This is standard notation from first-year abstract algebra.} difference in point of view occurs in other situations involving models;
for instance,
is $(\Z,+)$ (the integers with addition as operation) a model of the
axioms for a group, or is it correct to say that there is a model
for the group axioms that ``corresponds'' to $(\Z,+)$?

The distinction between $\C$ and $\FinLimTh[\C]$ is slightly different.  The
first might be described as presentation data for the second.  Since Lawvere,
many category theorists take the view that an algebraic
structure consists of the entire clone of operations rather than
some generating subset of operations.  From that point of view, the category
$\Csc$ is the ``same as'' $\FinLimTh[\C]$ rather than the same as
$\C$. Similarly, from the Lawverean perspective the group
discussed in the previous paragraph is {\it determined\/} by
saying $\Z$ is its underlying set and $+$ is its operation, but
the group is the entire clone. (This must be distinguished from
the relationship between a group and its presentation by
generators and relations, although of course there is an analogy
between those two situations.)

Similar remarks apply to the analogous constructions
obtained  when $\Catsk$ is replaced by an arbitrary finite limit
sketch.

In any case the following three categories are equivalent: \begin{enumerate}
\item The category of small categories and functors. \item
$\Mod[\Catsk,\Set]$, the category of models of the sketch $\Catsk$ in the
category $\Set$ with natural transformations as morphisms. \item The category
of finite-limit-preserving functors from $\FinLimTh[\Catsk]$ to $\Set$ with
natural transformations between them as morphisms. \end{enumerate}

In the sequel, we will distinguish between these
constructions typographically as an aid to reading the
monograph.  We continue the discussion of this section
about similar constructions in
Sections~\ref{moremodrem1}, \ref{modthform}
and~\ref{modrem2}.

\chapter [Finite-limit theories]{Construction
of finite-limit theories}\label{thconst}

\section{Introductory comments}

This chapter provides a construction of the categorial
theory $\FinLimTh[\S]$ of a finite-limit sketch $\S$. It is
essentially a special case of the construction
in~\cite{ehres68}. A related, more general construction is
given by Duval and Reynaud in~\cite{duvrey1}.

The categorial theory $\FinLimTh[\S]$ is in fact the initial
model of $\SynCat[\FinLim,\S]$, where $\FinLim$ is a finite-limit sketch
defined in Section~\ref{skedef}, for categories with finite
limits.  From that point of view, the categorial theory is
a term algebra for $\SynCat[\FinLim,\S]$. A recursive
construction of term algebras for finite-limit sketches is
given in~\cite{ctcs}, Section~9.2.

On the other hand, for any constructor space $\Eb$ and any
$\Eb$-form $\F$, $\SynCat[\Eb,\F]$ is (equivalent as a
category to) a particular example of a finite-limit theory
as described in~\ref{sec613}. Thus these rules of
construction may be used to construct actual factorizations
(when they exist) of potential factorizations (see Section~\ref{pfsec}) for any
$\Eb$-sketch.  It is in this sense that we have reduced
reasoning about arbitrary forms to reasoning (in the sense
of constructing actual factorizations) about finite limits.
The special character of $\Eb$-forms for a particular
constructor space $\Eb$ is encoded in the constructor-space
sketch that generates $\Eb$.

An approach to string-based logic that would be analogous
to this setup would be a system that captured many
different types of string-based logic (but not necessarily
all of them) using uniform rules about manipulating the
strings, with the special behavior for a particular type of
logic $L$, for example first order logic or linear logic,
encoded in a formal description of $L$ that caused the
uniform string manipulation rules to produce the correct
behavior for $L$. As far as we know no such system has been
defined. The one general approach to logic that we know
about, the theory of institutions described in
\cite{gogburst}, is much more abstract than the sort of
system about which we are speculating here; it does not
deal with strings and rules of deduction in a computational
way.

\section{A preliminary construction}\label{indnotsub}
Let $\Csc$ be a category with the property that for some
set (possibly empty) of diagrams in $\Csc$, limit cones
have been specified. In this section, we construct a graph
$G$ and a set of diagrams $D$ in $G$ by
$\Gamma$.\ref{firsto} to $\Gamma$.\ref{lasto} below; $G$
and $D$ are determined by $\Csc$. This construction is the
basis of the inductive construction given in
Section~\ref{indcon}. The definition is deliberately made
elementary (and therefore much more discursive than it
might be), since it forms the basis of our deduction rules
in Section~\ref{dedrules}.

\blist{$\Gamma$} \item\label{firsto} If $x$ is a node of
$\UndGr[\Csc]$, then $x$ is a node of $G$.
\item\label{arro} If $f:x\to y$ is an arrow of
$\UndGr[\Csc]$, then $f:x\to y$ is an arrow of $G$.
\item\label{compo} If the composite $g\oc f$ of two arrows
$f$ and $g$ of $\Csc$ is defined, then the diagram
\begin{equation}\label{p7top}
\xymatrix{ \dom(f) \ar[d]_{f}\ar[dr]^{g\oc f}
&\\
\cod(f) \ar[r]_{g} &\cod(g)
\\
}
\end{equation}
is in $D$.

\item\label{nextolasto} If $\delta :I\to\Csc$ is a diagram
in $\Csc$ that does not have a specified limit, then $G$
contains an object $\Lim[\delta]$ not in $\Csc$ and a cone
$\LimCone[\delta]:\Lim[\delta]\coneto\delta$,  and moreover
for each arrow $u:i\to j$ of $I$, the diagram
\begin{equation}\label{p7bottom}
\xymatrix{ &\Lim(\delta)
\ar[dl]_{\Proj[\LimCone[\delta],i]}
\ar[dr]^{\Proj[\LimCone[\delta],j]}
&\\
\delta(i)\ar[rr]_{\delta(u)} &&
\delta(j)\\
}
\end{equation}
is in $D$.

\item\label{lasto} If $\delta :I\to\Csc$ is a diagram in
$\Csc$ that does not have a specified limit in $\Csc$ and
$\Theta$ is a commutative cone to $\delta $,
then\begin{enumerate} \item\label{filin} $G$ contains an
arrow
$\Fillin[\Theta,\delta]:\Vertex[\Theta]\to\Lim[\delta] $
not in $\Csc$. \item For each node $i$ of $I$, the diagram
\begin{equation}\label{1p11}
\xymatrix{ \Vertex[\Theta] \ar[rr]^{\Fillin[\Theta,\delta]}
\ar[dr]_{\Proj[\Theta,i]} &&\Lim[\delta]
\ar[dl]^{\Proj[\LimCone[\delta],i]}
\\
&\delta(i)&\\
}
\end{equation}
is in $D$. \item If $k:\Vertex[\Theta]\to \Lim[\delta]$ is
in $\Csc$ and, for each node $i$ of $I$, the diagram
\begin{equation}\label{2p11}
\xymatrix{ \Vertex[\Theta] \ar[rr]^{k}
\ar[dr]_{\Proj[\Theta,i]} &&\Lim[\delta]
\ar[dl]^{\Proj[\LimCone[\delta],i]}
\\
&\delta(i)&\\
}
\end{equation}
is in $D$, then the diagram {\spreaddiagramcolumns{1.5cm}
\begin{equation}
\xymatrix{ \Vertex[\Theta]
\artwo {\Fillin[\Theta,\delta]}k &\Lim[\delta]}
\end{equation}\mpark{Reversed $k$}
} is in $D$.
\end{enumerate}
\elist

\section{The construction}\label{indcon}
In this section, we fix an arbitrary finite-limit sketch
$$\S\ceq (\Graph[\S],\Diagrams[\S],\Cones[\S])$$
In Definitions~\ref{seqdef1} through~\ref{seqdeflast}
below, we construct an infinite sequence
\begin{equation}
\xymatrix{ \Csc_0 \ar[r]^{F_0} & \Csc_1 \ar[r]^{F_1} &
\Csc_2 \ar[r]^{F_3} &\Csc_3\ldots
\\
}
\end{equation}
of categories and functors generated by $\S$.

\defn\label{seqdef1} Let
$D'$ be the set of diagrams consisting of all the diagrams
in $\Diagrams[\S]$ and all the diagrams necessary to make
each cone in $\Cones[\S]$ a formally commutative cone,
namely those of the form
\begin{equation}\label{makconcom}
\xymatrix{ &\Vertex[\Theta] \ar[dl]_{\Proj[\Theta,i]}
\ar[dr]^{\Proj[\Theta,j]}
&\\
\delta(i)\ar[rr]_{\delta(u)} &&
\delta(j)\\
}
\end{equation}
for each cone $\Theta:v\coneto\delta$ in $\Cones[\S]$ and
each arrow $u:i\to j$ of the shape graph of the base
diagram $\delta$ of $\Theta$. Then $\Csc_0$ is defined to
be the induced category $\LinTh[G,D']$ as in
Section~\ref{linsk}.

\edefn

\defn\label{seqdef2} Let $\S'\ceq (\Graph[\S'],\Diagrams[\S'])$
be the linear sketch generated by $\Csc_0$ as described in
Section~\ref{indnotsub}. Let $D''$ be the set of diagrams
in $\Graph[\S']$ containing all the diagrams in
$\Diagrams[\S']$ and, for each cone $\Theta:v\coneto\delta$
in $\Cones[\S]$ and each arrow
$k:\Vertex[\Theta]\to\Lim[\delta]$, the diagram
{\spreaddiagramcolumns{1.5cm}
\mpark{Reversed $k$}\begin{equation}
\xymatrix{ \Vertex[\Theta]
\artwo {\Fillin[\Theta,\delta]}k &\Lim[\delta]}
\end{equation}
}%
Then $\Csc_1$ is defined to be the induced category
$\LinTh[\Graph[\S],D'']$ and the image of each cone in
$\Cones[\S]$ is defined to be a specified limit (it is easy
to see that it is indeed a limit cone).\edefn

\defn\label{seqdeflast}
Assume $\Csc_k$ has been defined for $k\geq 1$. Then
$\Csc_{k+1}$ is defined to be $\LinTh[G,D']$ where $G$ and
$D$ are defined from $\Csc_k$ as in
Section~\ref{indnotsub}, and $F_i:\Csc_k\to\Csc_{k+1}$ is
defined to be the functor that takes each object and arrow
of $\Csc$ to its congruence class in $\Csc_{k+1}$. The
specified limits in $\Csc_{k+1}$ are the images of all
those in $\Csc_k$ plus all those constructed by rule
$\Gamma$.\ref{nextolasto}.

\edefn

\thm\label{constthm} The colimit $\T$ of the sequence of categories and functors defined by
Definitions~\ref{seqdef1} through ~\ref{seqdeflast} is a
category with finite limits. The induced graph morphism
from $G$ to  $\UndGr[\T]$ is a universal model of
$\S$.\ethm

It follows from the theorem that $\T$ is the finite-limit
theory of the sketch $\S$.  We denote it by $\FinLimTh[\S]$.
$\FinLimTh[\S]$ is equivalent as a category to
$\CatTh[\FinLim,\S]$ (Section~\ref{skedef}).

\pf Routine using the definition of colimit.\epf

\thl\label{piessens}{\hspace{-.75em}\bf (Piessens)} Every
object of $\FinLimTh[\S]$ is a limit of a diagram built from
the objects and arrows of $\S$. \ethl (See the remarks
in~\ref{forgetful}.)

\pf (Private communication from Frank Piessens.) Via the
Yoneda embedding, $\FinLimTh[\S]\op$ is a full subcategory of
the category $\Func(\Fsc,\Set)$, where $\Fsc$ the free
category generated by $\Graph[\S]$. Every functor from
$\Fsc$ to $\Set$ is a colimit of representables and, hence,
in $\FinLimTh[\S]$, every object is a limit of them. Thus every
object of the theory is a limit of a diagram built from the
objects and arrows of $\S$.\epf

\section{Rules of construction}\label{rulesapp}

These rules construct the objects, arrows and commutative
diagrams of the category $\FinLimTh[\S]$ for a given
finite-limit sketch $\S$.  The rules are given in two lists
in~\ref{arrrul} and~\ref{diagrul} below.

\subsection{Rules that construct objects and
arrows}\label{arrrul}

This list gives all the rules that construct objects and
arrows in the category $\FinLimTh[\S]$. The following
definition forces the distinguished cones of the sketch
$\S$ to become limit cones in the theory.

\defn\label{discondef}
Let $\Theta:v\coneto\delta$ be a distinguished cone of the
sketch $\S$. We define $\LimCone[\delta]\ceq\Theta$ and
$\Lim[\delta]\ceq v$.\edefn

\medskip

$\Rule{}{c}{$\exists$OB}{for every object $c$ of $\S$.}$

\medskip

$\Rule{}{\xymatrix@1{a \ar[r]^f & b}}{$\exists$ARR}%
{for every arrow $f:a\to b$ of $\S$.}$

\medskip

$\Rule{\xymatrix@1{  a\ar[r]^f & b\ar[r]^g & c}}%
{\xymatrix @1{a\ar[r]^{g\oc f} & c}}{$\exists$COMP}%
{for every object $b$ and pair of arrows $f:a\to b$ and
$g:b\to c$ of $\FinLimTh[\S]$.}$

\medskip

$\Rule{c}{ \xymatrix@1{ c\ar[r]^{\Id[c]} & c}}  {
$\exists$ID}{for every object $c$ of $\FinLimTh[\S]$.}$

\medskip

$\Rule{\delta:I\to\FinLimTh[\S]}%
{\LimCone[\delta]:\Lim[\delta]\coneto\delta}%
{$\exists$LIM}{for every diagram $\delta:I\to\FinLimTh[\S]$
that is not the base of a distinguished cone of $\S$.}$

\medskip

$\Rule{\Theta:v\coneto\delta}{\Fillin[\Theta,\delta]:v\to\Lim[\delta]}%
{$\exists$FIA}{for every diagram $\delta$ and every cone
$\Theta:v\coneto\delta$ in $\FinLimTh[\S]$.}$

\rem\label{rulerem} The first two rules are justified by
the inclusion of the sketch $\S$ into $\FinLimTh[\S]$. In rule
$\exists$LIM, $\LimCone[\delta]$ is the specified limit of
$\delta$. The exception in rule $\exists$LIM will force the
distinguished cones of $\S$ to become limit cones in
$\FinLimTh[\S]$: because of Definition~\ref{discondef}, rule
$\exists$FIA applies to those distinguished cones as well
as to all cones constructed by $\exists$LIM.  This remark
also applies to rules CFIA and !FIA in~\ref{diagrul}.

The rules $\exists$COMP corresponds to the arrow $\comp$
and $\exists$ID to  $\unit$ of the sketch for
categories~\ref{realcatsk}. $\exists$LIM and $\exists$FIA
correspond to arrows in $\FinLim$ but not specifically to
arrows of the sketch~\ref{finlimsk}, because an arbitrary
finite limit is constructed from a combination of products
and equalizers.\erem

\rem\label{rem2} The rules just given construct specific
objects and arrows in $\FinLimTh[\S]$.  Rule $\exists$LIM, for
example, constructs a specific limit cone called
$\LimCone[\delta]$, thus providing specified limits for
$\FinLimTh[\S]$.  It is true that there are other limit cones
in general for a given diagram $\delta$, but
$\LimCone[\delta]$ is a specific one.

Of course, in many cases, the entity constructed is the
\textit{unique} entity satisfying some property.  For example, the
arrow $g\oc f$ constructed by $\exists$COMP is (by
definition of commutative diagram) the only one making the
bottom diagram in COMPDIAG commute in $\FinLimTh[\S]$. The
arrow constructed by $\exists$ID is (by an easy theorem of
category theory) the only one making the bottom diagrams in
IDL and IDR commute.  The arrow constructed by $\exists$FIA
is (because of !FIA) the only one making the bottom diagram
in CFIA commute.  In connection with the point that each
rule constructs a specific arrow, these observations are
red herrings: in fact, each rule constructs a specific
arrow with the name given, independently of any uniqueness
properties arising from any other source. This point of
view is contrary to the spirit of category theory. We
follow it here because we are constructing syntax with an
eye toward implementation in a computer language.  This
situation is analogous to the way in which mathematicians
give invariant (basis-free) proofs concerning linear spaces
but use bases for calculation.

Implementing the specific constructions defined above would
be relatively straightforward using a modern
object-oriented language.  Note that we are not asserting
that it would be straightforward to find a confluent and
normalizing form of these rules for automatic theorem
proving, only that there are no obvious difficulties in
implementing them so that they could be applied in an {\it
ad-hoc\/} manner.\erem

\subsection{Rules that construct
formally commutative diagrams}\label{diagrul} The following
rules produce the existence of diagrams that must commute
in $\FinLimTh[\S]$.

\medskip

$\Rule{\xymatrix@1{ a \ar[r]^f &   b\\}} {\xymatrix@1{ a
\artwo ff & b\\}} {REF} {for every arrow $f:a\to b$ of
$\FinLimTh[\S]$}$

\medskip

$\Rule{\xymatrix@1{ a \artwo fg & b}
\quad\xymatrix@1{a \artwo gh & b\\
}}{\xymatrix@1{ a \artwo fh& b\\
}}{TRANS}{for all objects $a$ and $b$ and all arrows
$f,g,h:a\to b$ of $\FinLimTh[\S]$}$

\medskip

$\Rule{}{\UnivMod[\FinLim,\S]\oc\delta:I\to\FinLimTh[\S]}%
{$\exists$DIAG}{}$\\[6pt]
for every diagram $\delta$ in the set $D_{\S}$ of
distinguished diagrams of $\S$.

\medskip

$\Rule{
\xymatrix{ a \ar[r]^f & b \ar[d]^g  \\
& c \\ }}{
\xymatrix{ a \ar[r]^f \ar[dr]_{g\o f} & b \ar[d]^g \\
& c \\ }}{COMPDIAG}%
{for every pair of arrows $f:a\to b$ and $g:b\to c$ of
$\FinLimTh[\S]$.}$

\medskip

$\Rule{\begin{array}{c}
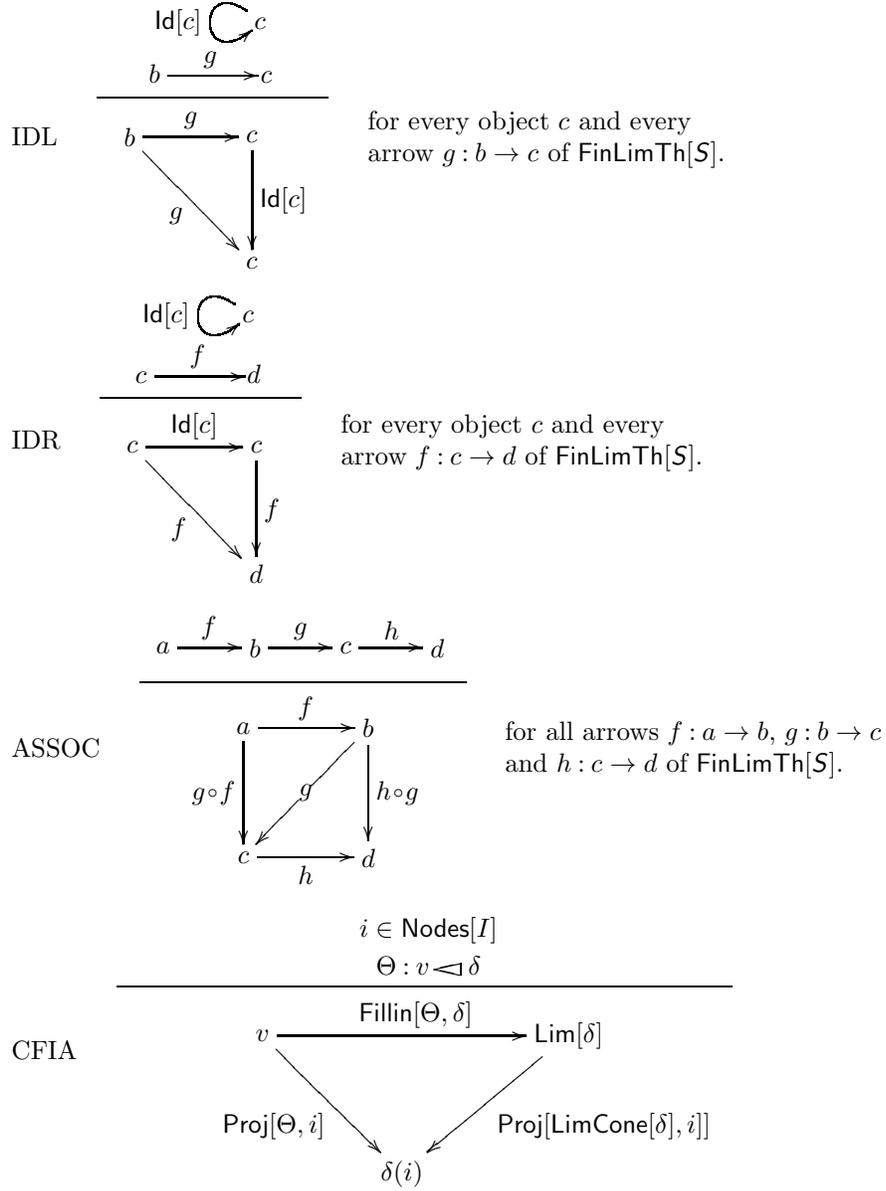
\xymatrix@1{ c\ar@(ul,dl)[]_{\Id[c]}}\\
\xymatrix@1{b \ar[r]^g& c}\end{array}}%
{\xymatrix{ b \ar[r]^g \ar[dr]_g & c \ar[d]^{\Id[c]}\\
& c\\ }
}%
{IDL}{for every object $c$ and every arrow $g:b\to c$ of
$\FinLimTh[\S]$.}$

\medskip

$\Rule{\begin{array}{c}
\xymatrix@1{ c\ar@(ul,dl)[]_{\Id[c]}}\\
\xymatrix@1{c \ar[r]^f & d\\}\end{array}}%
{\xymatrix{ c \ar[r]^{\Id[c]} \ar[dr]_f & c \ar[d]^f \\
& d \\ }
}%
{IDR}{for every object $c$ and every arrow $f:c\to d$ of
$\FinLimTh[\S]$.}$

\medskip

$\Rule{{\narrower\narrower
\xymatrix{ a \ar[r]^f &b\ar[r]^g &c\ar[r]^h& d }}}%
{\xymatrix{ a\ar[r]^f \ar[d]_{g\oc f} & b \ar[d]^{h\oc g}
\ar[dl]|g\\
c \ar[r]_h & d}}{ASSOC}{for all arrows $f:a\to b$, $g:b\to
c$ and $h:c\to d$ of $\FinLimTh[\S]$.}$

\medskip

 $\Rule{\hspace{3em}\begin{array}{c}
i\in \Nodes[I]\\
\Theta:v\coneto\delta
\end{array}\hspace{2.5em}}{\hspace{3em}
\xymatrix{ v \ar[rr]^(.5){\Fillin[ \Theta,\delta]}
\ar[dr]_{\Proj[ \Theta,i] } &&
 \Lim[ \delta]
\ar[dl]^{\hspace{1em}\Proj[ \LimCone[ \delta] ,i]]}\\
&\delta(i) & \\ }}%
{CFIA}{}$\\[6pt]
for every diagram $\delta:I\to\FinLimTh[\S]$, every node $i$ of
$I$, and every cone $\Theta$ with base diagram $\delta$.

\medskip

$\Rule{\begin{array}{c}
\delta:I\to\FinLimTh[\S]\\
\Theta:v\coneto\delta\\
h:v\to \Lim[\delta]   \\
k:v\to \Lim[\delta]\\
\mbox{and each of the following diagrams for each node $i$
of
$I$:}\\
{\xymatrix{ v \ar[rr]^(.5){h} \ar[dr]_{\Proj[ \Theta,i] }
&&
 \Lim[ \delta]
\ar[dl]^{\hspace{1em}\Proj[ \LimCone[ \delta] ,i]]}\\
&\delta(i) & \\
}}\\[10pt]
{\xymatrix{ v \ar[rr]^(.5){k} \ar[dr]_{\Proj[ \Theta,i] }
&&
 \Lim[ \delta]
\ar[dl]^{\hspace{1em}\Proj[ \LimCone[ \delta] ,i]]}\\
&\delta(i) & \\
}}
\end{array}}%
{\hspace{-3em}\xymatrix{ v \ar@<1ex>[r]^(.4){h} \ar@<-1ex>[r]_(.4){k} & \Lim[\delta]\\
}
}%
{!FIA}{}$\\[6pt]
for every diagram $\delta:I\to\FinLimTh[\S]$, every cone
$\Theta$ in $\FinLimTh[\S]$ with base diagram $\delta$, and
every pair of arrows $h,k:v\to\Lim[\delta]$.

\rem Note that we do not need a rule of the form

$$\Rule{\xymatrix{ a \ar[r]<1ex>^f \ar[r]<-1ex>_g & b\\
}}{\xymatrix{ a \ar[r]<1ex>^g \ar[r]<-1ex>_f & b\\
}}{SYM}{}$$ since the two diagrams exhibited are actually
the same diagram (see~\ref{diagsec}). \erem

\section{A specific choice of $\FinLimTh[\S]$}\label{specch}
In this monograph, for a given finite-limit sketch $\S$, we
assume given a particular instance of $\FinLimTh[\S]$: that
constructed in this chapter. It has the following
properties (which are not preserved by equivalence of
categories): \blist{T} \item $\FinLimTh[\S]$ is a category with
specified finite limits. (The construction explicitly
produces the specified limits.) \item Every arrow of
$\FinLimTh[\S]$ is a composite of projections from specified
limits, fill-in arrows and arrows of the form
$\FinLimUnivMod[\S](\f)$ for arrows $\f$ of the graph of $\S$.
\elist

The following proposition is significant for this monograph when
$\FinLimTh[\S]$ is taken to be $\SynCat[\Eb,\F]$ (defined as a
particular finite-limit theory in~\ref{termin} below),
where $\F$ is a form.

\thp\label{conobar} For a given sketch $\S$, every object
and every arrow of $\FinLimTh[\S]$ is constructible by repeated
applications of the constructions of Section~\ref{rulesapp}
to the objects and arrows of the sketch $\S$.\ethp

\pf This proof depends on the specific choice of
$\FinLimTh[\S]$ defined in~\ref{specch}.  It is clearly closed
under all the constructions of Section~\ref{rulesapp}. The
properties listed in~\ref{specch} imply that $\FinLimTh[\S]$ is
minimal with respect to the constructions of that section,
so that in fact those constructions can be taken as an
recursive definition of $\FinLimTh[\S]$.\epf

\rem It is also true by Lemma~\ref{piessens} that every
object of $\FinLimTh[\S]$ is the limit (not necessarily
specified) of a diagram of the form
$\FinLimUnivMod[\S]\oc\delta$ where $\delta$ is a diagram in the
graph of $\S$. This latter property, of course, is
preserved by equivalences of categories that commute with
the universal model.\erem

\notat It is clear that $\FinLimTh$ is a functor from the
category of finite-limit sketches and sketch morphisms to
the category of finite-limit categories and finite-limit
preserving functors. For any finite-limit sketch $\S$ and
morphism $\eta:\S\to\T$ of sketches, the induced functor
between categorial theories will be denoted by
$$\FinLimTh[\eta]:
\FinLimTh[\S]\to\FinLimTh[\T]$$ \enotat
\chapter{Limits of diagrams}\label{limdiagchap}

In this section, we develop some techniques for dealing
with limits of diagrams that are used extensively in the
example proofs in Sections~\ref{exam1} and \ref{exam2}.

\section{Morphisms of Diagrams}

\defn\label{dimordef} A {\bf morphism of diagrams}\mpar{You suggested using the double arrow for natural transformations and  calling it a 2-cell, but that is not right: it is a 2-cell in the 2-category of sets but it is \textit{not} a 2-cell in the category of diagrams, and that is what the word ``morphism'' refers to.}
$(\psi,\alpha):(\delta:I\to G)\to(\delta':I'\to G)$
is a graph morphism $\psi:I\to I'$ together with a
natural transformation
$\alpha:\delta'\oc\psi\to\delta$. \edefn

\rem It is easy to see that this definition of
morphism of diagrams is compatible with the
equivalence relation that defines diagrams. It was
first given by Eilenberg and Mac Lane
\shortcite{eilmaclane} and studied further in
\cite{kocklimmon,guit74,guit77}. It is not the same
notion of morphism of diagrams as that
of~\cite{barrexact}, page~52, studied
in~\cite{tholtozz}.\erem

\defn\label{strictmordef} A {\bf strict morphism of
diagrams} $\psi:(\delta:I\to G)\to(\delta':I'\to
G)$ is a graph morphism $\psi:I\to I'$ for which
the diagram


\begin{equation}\label{diagdiag2}
\xymatrix{
I \ar[rr]^\psi \ar[dr]_\delta &&
I'\ar[dl]^{\delta'} \\ & G \\ } \end{equation}
commutes.\edefn

\rem A strict morphism of diagrams is a special case
of morphism of diagrams (set $\alpha$ to be the
inclusion).  All applications in this monograph use
strict morphisms only.\erem

\thp\label{inducedprop} Let
$(\psi,\alpha):(\delta:I\to \Csc)\to(\delta':I'\to
\Csc)$ be a morphism of diagrams in a category
$\Csc$. Given a commutative cone
$$\Theta':\Vertex[\Theta']\coneto\bigl(\delta':I'\to \Csc\bigr)$$
there is a commutative cone
$$\Theta:\Vertex[\Theta]\coneto\bigl(\delta:I\to
\Csc\bigr)$$
with the following properties:
\alist
\item
$\Vertex[\Theta] = \Vertex[\Theta']$
\item For every $i\in \Nodes[I]$, $\Proj[\Theta,i]=\alpha i\oc
\Proj[\Theta',\psi i]$.
\ealist
\ethp


\pf That $\Theta$ is commutative follows from the
fact that $\Vertex[\Theta] = \Vertex[\Theta']$ and
the fact that for every $f:i\to j$ in $I$, the
following diagram commutes because $\Theta'$ is
commutative and $\alpha$ is natural.

{\shorter\wider\wider

{\shorter\wider\wider
\begin{equation}
\xymatrix{
&\delta'\psi i\ar[r]^{\alpha i}\ar[dd]^{\delta'\psi
f}&\delta i\ar[dd]^{\delta f}\\
\Vertex[\Theta']\ar[ur]^{\Proj[\Theta',\psi i]}
 \ar[dr]_
{\Proj[\Theta',\psi j]}\\
& \delta'\psi j\ar[r]^{\alpha
j}&\delta j\\
}
\end{equation}
}
\thc\label{morcor} Let $(\psi,\alpha):(\delta:I\to
G)\to(\delta':I'\to G)$ be a morphism of
diagrams. Then there is a unique arrow
$\phi:\Lim[\delta']\to\Lim[\delta] $
for which for all nodes $i$ of $I$,\mpark{\rm{Repaired bottom line of diagram}}
{\wider\wider\taller\taller

\begin{equation}
\xymatrix{
\Lim[\delta'] \ar[r]^\phi
\ar[d]_{\Proj[\Lim[\delta'],\psi i]} & \Lim[\delta]
\ar[d]^{\Proj[\Lim[\delta],i]} \\ \delta'\psi i
\ar[r]_{\alpha i} & \delta i } \end{equation}
}\
 \ethc

\pf
This follows from
Proposition~\ref{inducedprop} by letting
$\Theta'\ceq\LimCone[\delta']$ and then setting
$\phi$ to be the fill-in arrow from $\Theta$ to
$\Lim[\delta]$, where $\Theta$ is the cone defined
in Proposition~\ref{inducedprop}. \epf

\rem
When the target of the diagrams is a category
$\Csc$ with finite limits, the preceding
constructions make $\Lim$ a contravariant functor
from the category of diagrams to $\Csc$.\erem

\section{Restrictions of diagrams}

\defn
Let
$\delta:I\to\Csc$ be a diagram and
$\Incl[J\includedin I]:J\to I$
an inclusion
of graphs.
The {\bf restriction} of $\delta$ to $J$, denoted
by $\delta\|_J$, is $\delta\oc\Incl[J\includedin
I]:J\to \Csc$.  $\delta\|_J$ is called a {\bf
subdiagram} of $\delta$.\edefn

\rem $\Incl[J\includedin I]$ is a strict
morphism of diagrams from $\delta\|_J$ to
$\delta$.\erem

\defn Let $\delta:I\to\Csc$ be a diagram and
$\Incl[J\includedin I]:J\to I$ an inclusion of
graphs. Let $\Theta:v\coneto\delta$ be a cone.  The
{\bf base-restriction of} $\Theta$ {\bf to}  $J$ is
defined to be the cone
$\Theta\|^J:v\coneto(\delta\oc \Incl[J\includedin
I])$ with vertex $v$ and projections defined by
$\Proj[\Theta\|^J,j]\ceq\Proj[\Theta,j]:v\to
\delta(j)$ for all nodes $j$ of $J$. In this case,
we also say that $\Theta$ is a {\bf base-extension}
of $\Theta\|^J$. \edefn \rem  If $\Theta$ is
commutative, then so is $\Theta\|^J$.\erem

\defn  Let $\delta:I\to\Csc$ be a diagram and $J$
a subgraph of $I$. Then the subdiagram $\delta\|_J$
is said to {\bf dominate} $\delta$, or to be {\bf
dominant in} $\delta$, if every commutative cone
$\Theta:v\coneto(\delta\|_J)$ in $\Csc$ has a
unique base extension to a commutative cone
$\Theta':v\coneto\delta$ with the same vertex.
\edefn

\rem Tholen and Tozzi~\shortcite{tholtozz} give a
condition (``confinality'') on $I$ and $J$ such
that any diagram based on $I$ is dominated by its
restriction to $J$.  One type of dominance that
their condition does not cover is the case in which
$\delta$ is obtained from $\delta\|_J$ by adjoining
a limit cone over a subdiagram of $\delta\|_J$ (see
Section~\ref{speccas}.) \erem


\section{Limits of subdiagrams}\label{limsubdiag}

\rem Some\mpark{Moved remark to here from 5.2} of the definitions and lemmas in
this section have variants in which one
has graph homomorphisms rather than inclusions.  We
shall not, however, need these.\erem
\thl\label{lemma1} Let $\delta:I\to \Csc$ be a
diagram and let $J$ be a subgraph of $I$ with
inclusion $\Incl[J\includedin I]$. Let
$\gamma=\delta\|_J$.  Then there is a unique arrow
$\phi:\Lim[\delta]\to\Lim[\gamma]$ such that for
all nodes $j$ of $J$,

\begin{equation}
\xymatrix{
\Lim[\delta] \ar[rr]^\phi
\ar[dr]_{\Proj[\LimCone[\delta],j]} &&
\Lim[\gamma]\ar[dl]^{\Proj[\LimCone[\delta],j]}
\\ & \delta(j) \\ }
\end{equation}
Furthermore,\mpark{\rm{Added sentence}} if $J$ dominates $I$, then $\phi$ is an isomorphism.\ethl
\pf The existence and
uniqueness of $\phi$ is a special case of
Corollary~\ref{morcor}.

Now assume that $\gamma$ dominates $\delta$. Let
$\Psi:\Lim[\gamma]\coneto\delta$ be the unique
extension of $\LimCone[\gamma]$ to $\delta$.  Using
$\exists$FIA of Section~\ref{rulesapp}, we define
$$\psi\ceq\Fillin[\Psi,\delta]:\Lim[\gamma]\to\Lim[\delta]$$
It follows from  !FIA of Section~\ref{rulesapp} that
$\psi$ is the only arrow from $\Lim[\gamma]$ to
$\Lim[\delta]$ that makes all diagrams of the form\mpark{Changed top arrow to $\psi$}

\begin{equation}
\xymatrix{
\Lim[\gamma] \ar[rr]^\psi
\ar[dr]_{\Proj[\psi],j]} &&
\Lim[\delta]\ar[dl]^{\Proj[\LimCone[\delta],j]}
\\ & \gamma(j) \\ }
\end{equation}
commute for each node $j$ of $J$.
Since $\phi\oc\psi:\Lim[\gamma]\to\Lim[\gamma]$ and
$\Id[\Lim[\gamma]]:\Lim[\gamma]\to\Lim[\gamma]$ both commute
with all the projections to nodes of $J$,
it follows from !FIA that
$\phi\oc\psi=\Id[\Lim[\gamma]]$, A similar argument shows that
$\psi\oc\phi=\Id[\Lim[\delta]]$, so that $\phi$ is an
isomorphism. \epf

\section{Special cases of extending diagrams}\label{speccas}

Here we define some special cases of dominance that are easy to
recognize.

\defn\label{addarrow}
Let  graphs $I$ and $J$ be given such that $J\includedin I$ and $I$ and $J$
have the same nodes, and suppose that $I$ has exactly one arrow $a:j\to k$ not
in $J$. Let $\delta:I\to\Csc$ be a diagram with the property that for all
nodes $j'$ of $J$ and all arrows $f:j\to j'$ and $g:j'\to k$,


\begin{equation}\label{adjcompdef}
\xymatrix{
\delta(j) \ar[rr]^{\delta(a)}
\ar[dr]_{\delta(f)} && {\delta(k)} \\
& \delta(j')\ar[ur]_{\delta(g)} \\ }
\end{equation} 
commutes  in $\Csc$. Then we say $\delta$ {\bf
extends} $\delta\|_J$ {\bf by adjoining a composite}.\edefn

\defn\label{addcocone}
Let $I$ and $J$ be graphs with the following properties:
\blist{ACC}
\item $J\includedin I$.
\item $I$ has exactly one node $v$ not in $J$.
\item $I$ has at least one arrow not in $J$.\mpark{Should ACC-3 be omitted?}
\item Every arrow in $I$ not in $J$ has target $v$.
\elist Suppose that $\delta:I \to\Csc$ is a diagram with the
property that if $a:i\to v$, $b:j\to v$ and $f:i\to j$ are
arrows of $I$, then


\begin{equation}
\xymatrix{
\delta(i) \ar[rr]^{\delta(a)}
\ar[dr]_{\delta(f)} && {\delta(v)} \\
& \delta(j)\ar[ur]_{\delta(b)} \\ }
\end{equation} 
commutes.  Then we say $\delta$ {\bf extends} $\delta\|_J$ by
{\bf adjoining a commutative cocone}.\edefn

\defn\label{adjlim}
Let $I$, $J$ and $J'$ be graphs with $J'\includedin J\includedin
I$, such that $J'$ is full in $J$, $I$ contains exactly one
node $v$ not in $J$, and for each node $j$ of $J'$, $I$ contains exactly
one arrow $p_j:v\to j$ and no other arrows not in $J$.
Let $\delta:I\to\Csc$ be a diagram, and suppose further that
$\delta$ extends $\delta\|_J$
in such a way that
$\delta(v)$ and the
arrows $\delta(p_j)$ constitute a limit cone to $\delta\|_{J'}$.
Then we say that $\delta$ {\bf extends} $\delta\|_{J}$ {\bf by
adjoining a limit.}\edefn

\defn\label{addfillin}
Let $I$ be a graph
and let $\delta:I\to\Csc$ be a diagram.  Let $J'$ be a
nonempty subgraph of $I$ and let $\Theta:v\coneto\delta\|_{J'}$
and $\Psi:w\coneto\delta\|_{J'}$ be commutative cones for
which\alist \item $\Theta$ is a limit cone. \item
Each projection $\Proj[\Theta,i]$ and
$\Proj[\Psi,i]$ is a composite of arrows in the image of
$\delta$ (it follows that $v$ and $w$ are in the image of
$\delta$.)

\ealist Let $\phi:w\to v$ be the unique
fill-in arrow given by the definition of limit, and suppose $f$
is an arrow of $I$
for which $\delta(f)=\phi$. Let $J$ be the subdiagram of $I$
obtained by omitting $f$. Then $\delta$ {\bf extends}
$\delta\|_J$ {\bf by adjoining a fill-in arrow}.\edefn

\thl\label{lemma2}
Suppose that $\delta':I\to\Csc$ extends $\delta:J\to\Csc$ by adjoining a
composite, a commutative cocone, a limit or a fill-in arrow.
Then
$$\Fillin\left[\LimCone[\delta'\|_J],\delta\right]:\Lim[\delta']\to\Lim[\delta]$$
is an isomorphism.
\ethl

\pf We will
show in each case that
$\delta\|_J$ dominates $\delta$.

In the case of adjoining a composite, it follows from
the fact that all the diagrams~(\ref{adjcompdef})
must commute that a commutative cone over $\delta\|_J$ is already
a commutative cone over $\delta$.

If $\delta$ extends $\delta\|_J$ by adjoining a commutative
cocone, then in the notation of Definition~\ref{addcocone} any
cone $\Theta:u\coneto\delta\|_J$ extends uniquely to a cone
$\Theta':u\coneto\delta$ by defining $\Proj[\Theta',v]\ceq
\delta(f)\oc\Proj[\Theta,i]$, where $f:i\to v$ is an arrow of $I$
not in $J$.

If $\delta$ extends $\delta\|_J$ by adjoining a limit, then
in the notation of Definition~\ref{adjlim},\mpark{Restored missing sentence}
$\Theta:u\coneto\delta\|_J$ extends uniquely to
$\Theta':u\coneto\delta$ by defining
$\Proj[\Theta',v]\ceq
\Fillin\left[\Theta',\delta\|_{J'}\right]$.

Finally, suppose $\delta$ extends
$\delta\|_J$ by adjoining a fill-in arrow.
By repeatedly adjoining composites we can assume $\delta$ has
the property that every projection arrow
$\Proj[\Theta,i]$ and $\Proj[\Psi,i]$ (notation as in
Definition~\ref{addfillin}) is in the image of $\delta$.
(We are using the fact that dominance is transitive, which is
easy to show.)
Now let $\Phi:x\coneto\delta\|_J$ be a
commutative cone and let $\delta(m)=v$, $\delta(n)=w$.
It is necessary and sufficient to show that the diagram

\begin{equation}\label{tring}
\xymatrix{
&w\ar[dd]^{\delta f}\\
x\ar[ur]^{\Proj[\Phi,n]}
 \ar[dr]_
{\Proj[\Phi,m]}\\
& v\\
}
\end{equation} 
}

commutes.


For every arrow $g:j\to j'$ of $J'$, we have a diagram


{\Wider\Wider\Wider\taller\taller
\begin{equation}\label{xdia}
\xymatrix{
&w\ar[r]^{\Proj[\Psi,j]}
\ar[ddr]|(.25){\Proj[\Psi,k]}
\ar[dd]^{\delta f}
&\delta (j)
\ar[dd]^{\delta(g)}\\
x\ar[ur]^{\Proj[\Phi,n]}
\ar[dr]_{\Proj[\Phi,m]}\\
&v\ar[uur]|(.25){\Proj[\Theta,j]}\ar[r]_{\Proj[\Theta,k]}
&\delta (k)\\}
\end{equation}
}
Let the cone $\Psi':x\coneto\delta\|_{J'}$ be defined
by requiring that
$$\Proj[\Psi',j]=\Proj[\Psi,j]\oc\Proj[\Phi,n]$$ for
every node $j$ of $J'$. The upper right triangle of
Diagram~(\ref{xdia}) commutes because $\Psi$ is a
commutative cone.  It follows that $\Psi'$ is a
commutative cone.  We now prove that both
$\Proj[\Phi,m]$ and $\delta f\oc\Proj[\Phi,n]$ satisfy
the requirements of $\Fillin[\Psi',\delta\|_{J'}]$ in
the notation of of Section~\ref{rulesapp}.  It will
follow from rule !FIA in that section that
Diagram~(\ref{tring}) commutes, as required.

\alist\item We must show that
for all nodes $j$ of $J'$,
$$\Proj[\Psi,j]\oc\Proj[\Phi,n]=\Proj[\Theta,j]\oc\Proj[\Phi,m]$$
This follows from the fact that $\Phi$ is a commutative cone to
$J$ and $J'\includedin J$.

\item We must show that for all nodes $j$ of $J'$,
$$\Proj[\Theta,j]\oc\delta
f\oc\Proj[\Phi,n]=\Proj[\Psi,j]\oc\Proj[\Phi,n]$$
This follows from the fact that the upper left triangle inside the
rectangle in Diagram~(\ref{xdia}) commutes because $\delta f$ is
a fill-in arrow.
\ealist
\epf


\chapter{Forms}\label{skgeneral}

We provide here a definition of ``form'' (generalized
sketch) based on \cite{gensk}. Some of the terminology has
been changed. A form is a generalization of the concept of
sketch invented by Charles Ehresmann and described in
\cite{bastehres} or in \cite{ctcs}. The definitions below
presuppose the concept of finite-limit sketch (see
Section~\ref{oldskdefs}).

\section{Constructor spaces}\label{cspaces}
We will assume given a fixed finite-limit sketch
$\Catsk$ whose category of models is the category of small
categories and functors. A specific such sketch is given
explicitly in Section~\ref{realcatsk}. The presentation
that follows has $\Catsk$ as an implicit parameter.

\defn\label{cssdef}
A finite-limit sketch $\E$ together with a morphism $\eta:\Catsk\to\E$ of
sketches
is called a {\bf constructor space sketch}
provided that every object in
$\FinLimTh[\E]$ is the limit of a finite diagram whose nodes are of
the form  $\FinLimTh[\eta](\n)$,
where $\n$ is a node of
$\Catsk$. The morphism $\eta$ is denoted
by $\CatStruc[\E]:\Catsk\to\E$. \edefn

\rems The notation ``$\CatStruc$'' abbreviates ``categorial
structure''.

Definition~\ref{cssdef} is more general than
Definition~4.1.2 in~\cite{gensk} in that
$\CatStruc[\E]$ need not be an inclusion. However, in
all the examples in this monograph, $\FinLimTh[\CatStruc[\E]]$
is injective on objects.\erems

\defn\label{cspacedef}
A category of the form $\FinLimTh[\E]$ for some constructor
space sketch $\E$ is called a {\bf constructor
space}.\edefn \notat\label{bfnot} We will normally denote
the constructor space $\FinLimTh[\E]$ by $\Eb$ (note the
difference in fonts). In particular, we have the
constructor space $\Cat$ corresponding to the constructor
space sketch $\Catsk$ given in~\ref{realcatsk}.  For this
example, $\CatStruc[\Catsk]:\Catsk\to\Catsk$ is the
identity functor. For the constructor spaces $\FinProd$,
$\FinLim$ and $\CCC$ constructed in Chapter~\ref{conspsk},
the structure map is in each case inclusion.

\edefn \defn\label{ecatdef} A model in
$\Set$ of a constructor space $\Eb$ is called an
{\bf$\Eb$-category}, and a morphism of such models is called an
{\bf$\Eb$-functor} (see Section~\ref{moremodrem1}). \edefn

\rems Recall that a model of $\Eb$ is a
finite-limit preserving functor from $\Eb$ to $\Set$, and a
morphism of models is a natural transformation from one such
functor to another.

Observe that $\CatStruc$ induces an underlying functor from
the category of $\Eb$-categories to the category of
categories.

Definitions~\ref{cssdef},
\ref{cspacedef}
and~\ref{ecatdef} are essentially the same as those in \cite{gensk}, and are a
special case (where all 2-cells are identities) of the two-dimensional
version given in~\cite{powerwells}.
\erems

\section{Further remarks concerning models}\label{moremodrem1}
We continue the discussion about models begun in
Section~\ref{modrem1}. Constructor spaces $\FinProd$ (for
categories with specified finite products), $\FinLim$ (for
categories with specified finite limits) and $\CCC$ (for
Cartesian closed categories with specified structure) are
given in Chapter~\ref{conspsk}.  The remarks concerning
models of $\Catsk$ in Section~\ref{modrem1} apply equally
well to models of these and other constructor spaces.

For\mpar{amplified} example, each
model of $\CCC$  is a functor, but it
corresponds to a certain Cartesian closed category {\it
with specified structure\/} whose objects, arrows, sources
and targets, composition, binary product structure and
closed structure are all determined by the values (in the
model under consideration) of certain nodes and arrows of
the sketch $\CCC$. Morphisms of models are Cartesian closed functors that preserve all this specified structure on the nose.  Cartesian closed categories in the usual sense form a large category isomorphic to the category of models of $\CCC$.

We will identify Cartesian closed
categories with models in $\Set$ of $\CCC$ in the sequel,
and similarly for other constructor spaces $\Eb$.  In
particular a $\FinProd$ category is a category with
specified finite products and a $\FinLim$ category is a
category with certain specified finite limits.

The value in a model $\M$ of an object $\vsfi$ in a constructor space
is the set of all examples of a particular construction that is
possible in the $\Eb$-category $\M$.  Section~\ref{consnot} gives an extended example of this.
Thus each object of $\Eb$ represents a {\it type of construction\/}
possible in an $\Eb$-category; hence the name ``constructor
space''.


\section{Notation for diagrams in a constructor
space}\label{consnot}
The object of $\FinLimTh[\Catsk]$
whose value in a model is the set of all not
necessarily commutative diagrams of the form

\begin{equation}\label{sqdiag1}
\xymatrix{
A \ar[r]^h \ar[d]_f & B
\ar[d]^k \\
C \ar[r]_g\ar[ur]^x & D
}
\end{equation}

is the limit of the diagram
{\wider\wider\taller
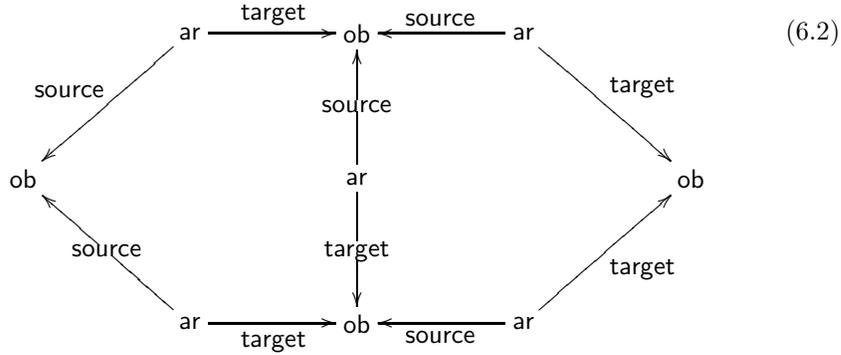
\begin{equation}\label{basicunlab}
\xymatrix{
&\arr\ar[r]^{\target}\ar[dl]_{\source}
&\ob
&\arr\ar[l]_{\source}\ar[dr]^{\target}
\\
\ob
&&\arr\ar[u]|{\source}\ar[d]|{\target}
&&\ob
\\
&\arr\ar[r]_{\target}\ar[ul]|{\source} &\ob
&\arr\ar[l]^{\source}\ar[ur]_{\target}
\\
}
\end{equation}
}
Observe that $\FinLimTh[\Catsk]$
is the constructor space for unrestricted categories, so that
Diagram~(\ref{basicunlab}) (more precisely, its image under
$\CatStruc[E]$) occurs in any constructor space $\Eb$.

We now describe this diagram in more detail and introduce some
notation that makes the discussion of such diagrams easier to
follow. We use the notation $D(n)$ to refer to the diagram shown
herein with label $(n)$, and $I(n)$ for its shape graph. For
example, the limit of the diagram above is
$\Lim[D(\ref{basicunlab})]$.

Every node of $D(\ref{basicunlab})$ is
either the object $\ob$ (the object that becomes the set of
objects in a model) or the object $\arr$ (the
object that becomes the set of arrows in a model)
of $\FinLimTh[\Catsk]$. For a model $\C$ of
$\FinLimTh[\Catsk]$ in $\Set$, an element of
$\C(\Lim[D(\ref{basicunlab})]$ is a diagram in
$\C$, not
necessarily commutative, of the form of Diagram~(\ref{sqdiag1}).

In order to make the relation between Diagrams~(\ref{sqdiag1})
and~(\ref{basicunlab}) clear, we give the
shape graph of~(\ref{basicunlab}):

{\wider\wider\taller
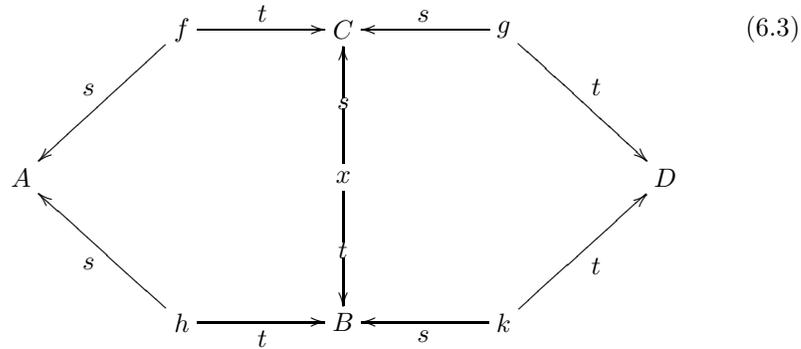
\begin{equation}\label{basicshape}
\xymatrix{
&f\ar[r]^{t}\ar[dl]_{s}
&C
&g\ar[l]_{s}\ar[dr]^{t}
\\
A
&&x\ar[u]|{s}\ar[d]|{t}
&&D
\\
&h\ar[r]_{t}\ar[ul]^{s}
&B
&k\ar[l]^{s}\ar[ur]_{t}
\\
}
\end{equation}
}

We have labeled the nodes of Diagram~(\ref{basicshape}) by the objects and
arrows that occur in Diagram~(\ref{sqdiag1}) in such a way that the node
named by an object or arrow of Diagram~(\ref{sqdiag1}) will inhabit the
value of that node in the model $\C$.  For example, the object
$A$ of $\C$ is at the upper left corner of Diagram~(\ref{sqdiag1}) and
the projection arrow from $\C(\Lim[D(\ref{basicunlab})])$
to $\C(\ob)$ determined by the node labeled $A$ of the shape
graph~(\ref{basicshape})
is a function from the set of diagrams
in $\C$ of the form of Diagram~(\ref{sqdiag1})
to the set of objects of $\C$ that takes a diagram to the object
in its upper
left corner.
The  arrows of Diagram~(\ref{basicshape}) are
labeled in accordance to their values in
Diagram~(\ref{basicunlab}).  {\it It is
important to understand that each distinct arrow in
Diagram~(\ref{basicshape}) is a different arrow of the shape
graph}, whether they have different labels or not.

We will combine diagrams such as
Diagram~(\ref{basicunlab})
and their shape graphs into one graph by labeling the nodes of
the diagram by superscripts naming the corresponding node of the
shape graph.  In the case of Diagram~(\ref{basicunlab}), doing this
gives the
following {\bf annotated diagram}:

{\wider\wider\taller
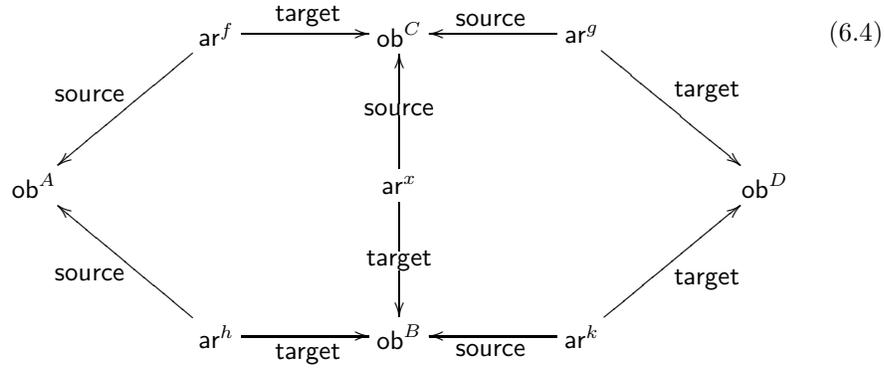
\begin{equation}\label{basic}
\xymatrix{
&\arr^{f}\ar[r]^{\target}\ar[dl]_{\source}
&\ob^{C}
&\arr^{g}\ar[l]_{\source}\ar[dr]^{\target}
\\
\ob^{A}
&&\arr^{x}\ar[u]|{\source}\ar[d]|{\target}
&&\ob^{D}
\\
&\arr^{h}\ar[r]_{\target}\ar[ul]^{\source}
&\ob^{B}
&\arr^{k}\ar[l]^{\source}\ar[ur]_{\target}
\\
}
\end{equation}
}

Here, the superscript $A$ on the leftmost node indicates that
the corresponding node of the shape graph is labeled $A$.
Formally,
the expression $\ob^A$ is used as the label for
the node $\delta(A)$ and its use signifies that
$\delta(A)=\ob$.
That device
helps the reader to see that Diagram~(\ref{sqdiag1})
is indeed an element of $\C(\Lim[D(\ref{basicunlab})])$.

For example, the particular
arrow $h$ of Diagram~(\ref{sqdiag1}) is an element
of $\arr$, and the label $\arr^h$ in
Diagram~(\ref{basic}) helps one see that it is that
node that projects to $h$ in the model $\C$ and
that the source of $h$ is $A$ and that the target
is $B$.

{\it It is important to understand that an annotated diagram
such as~(\ref{basic}) denotes precisely the same diagram
as~(\ref{basicunlab}).}  The fact that one node is
labeled $\ob^A$ and another $\ob^B$ does not change the fact
that {\it both nodes\/} are $\ob$.  The superscript merely
gives information about the relation between
Diagram~(\ref{basicunlab}) and Diagram~(\ref{sqdiag1}).

Diagram~(\ref{basic})
could also be drawn as the base of a limit cone
$\Theta$
with limit $\Lim[D(\ref{basic})]$ (which of course is the
same as $\Lim[D(\ref{basicunlab})]$) as
follows.
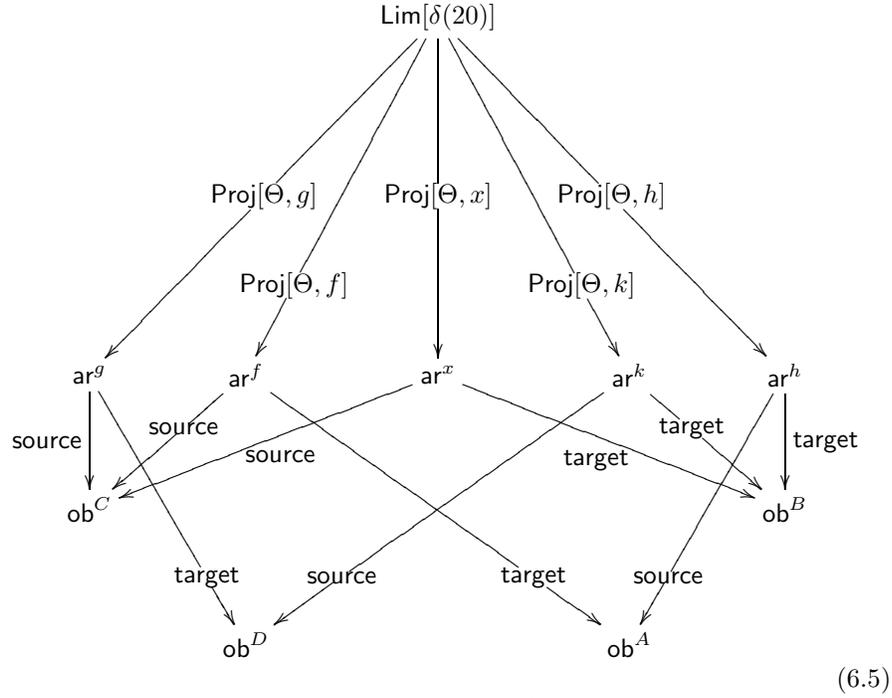
\begin{equation}\label{olimcone}
\xymatrix{
&&\Lim[\delta(20)]
\ar[dddll]|{\Proj[\Theta,g]}
\ar[dddl]|(.75){\Proj[\Theta,f]}
\ar[ddd]|{\Proj[\Theta,x]}
\ar[dddr]|(.75){\Proj[\Theta,k]}
\ar[dddrr]|{\Proj[\Theta,h]}
&&\\
&&&&\\
&&&&\\
\arr^{g}\ar[d]_{\source}\ar[ddr]|(.75){\target}
&\arr^{f}\ar[dl]|(.4){\source}\ar[ddrr]|(.75){\target}
&\arr^{x}\ar[dll]^{\source}\ar[drr]_{\target}
&\arr^{k}\ar[dr]|(.4){\target}\ar[ddll]|(.75){\source}
&\arr^{h}\ar[d]^{\target}\ar[ddl]|(.75){\source}
\\
\ob^{C}
&&&&\ob^{B}
\\
&\ob^{D}
&&\ob^{A}
&
\\
}
\end{equation}


Because
of the typographical complexity of doing this for diagrams more
complicated than Diagram~(\ref{basic}), we will usually give diagrams
whose limits we discuss in the form of Diagram~(\ref{basic}),
without showing the cone, instead of in the form of
Diagram~(\ref{olimcone}).

Showing the cone explicitly as in
Diagram~(\ref{olimcone})
nevertheless has an advantage.
It makes it clear
that many of the projection arrows from
$\Lim[D(\ref{basic})]$ are induced by others; in the
particular case of Diagram~(\ref{olimcone}),
all the arrows to nodes labeled
$\ob$ are induced by composing arrows to some node
labeled $\arr$ with $\source$ or $\target$.
Diagram~(\ref{basic}) does not make this property
as easy to discover as Diagram~(\ref{olimcone})
does.

A systematic method of translating from graphical
expressions such as Diagram~(\ref{olimcone}) to a string-based
expression could presumably be based on this, following the
notation
introduced in \cite{ttt}, page 38. In the case of
Diagram~(\ref{olimcone}), the string-based expression would be
something like this: $$ \begin{array}{l}
\bigl[<g,f,x,k,h>\mid\sr(g)=C,\ \tar(g)=D,\ \sr(f)=A,\\
\tar(f)=C,\ \sr(x)=C,\ \tar(x)=B,\ \sr(k)=B,\\ \quad \tar(k)=D,\
\sr(h)=A,\ \tar(h)=B\bigr]\end{array}$$ or in more familiar terms,
$$\left[<g,f,x,k,h>\mid g:C\to D, f:A\to C, x:C\to B, k:B\to D,
h:A\to B\right]$$

\section{Forms}\label{skedef}\label{formdes}

In this section, we outline  those facts about forms
that are needed in this monograph. More complete treatments are in
\cite{gensk} and \cite{powerwells}.

Let $\E$ be a constructor-space sketch, $\Eb$ (which is
$\FinLimTh[\E]$) the constructor space it generates, and
$\delta:I\to\Graph[\E]$ a diagram. We may freely adjoin a
global element $\phi:1\to\Lim[\delta]$ to obtain a
finite-limit category, denoted by $\Eb[\phi]$ in the
literature and called a {\bf polynomial category}.

\defn\label{formdef} In the notation of the
preceding paragraph, the {\bf $\Eb$-form} $\F$ {\bf
determined by $\delta$} is the value $\I(\phi)$ of a freely
adjoined global element $\phi:1\to\Lim[\delta]$, where $\I$
is the initial model of $\Eb[\phi]$ in $\Set$. \edefn

\notat\label{termin} If $\F$ is an $\Eb$-form determined by
$\delta$ as in the definition, we write $\Fa$ for $\phi$.
The diagram $\delta$ is called the {\bf description} of
$\F$. We denote $\Eb[\Fa]$ by $\SynCat[\Eb,\F]$ and call it
the {\bf syntactic category} of $\F$.\enotat

\rem The ``$\Eb$'' in the notation $\SynCat[\Eb,\F]$ is
redundant, but helpful as a reminder of which constructor
space we are using.\erem

\section{Constructing $\SynCat[\Eb,\F]$}\label{sec613} One
way of constructing $\SynCat[\Eb,\F]$ is as follows: First
adjoin $\Fa$ to $\Graph[\E]$ to obtain a graph $G$. Then
define the finite-limit sketch $\S\ceq
(G,\Diagrams[\E],\Cones[\E])$, and finally set
$\SynCat[\E,F]\ceq\FinLimTh[\S]$. The inclusion of $\Graph[\E]$
into $G$ is a sketch map from $\E$ to $\S$ and so generates
a finite-limit\mpar{``finite limits preserving'' is not idiomatic English.  After all, you don't say ``fruits harvesting machine''.} preserving functor
$\Constants[\F]:\Eb\to\SynCat[\Eb,\F]$. This construction
can cause considerable collapsing, for example if one
adjoins a global element of an initial object.

\rems\label{forgetful} A model of $\Eb$ in $\Set$ (a
finite-limit pre\-serv\-ing functor from $\Eb$ to $\Set$)
is an $\Eb$-category. A model $\Ff$ of $\SynCat[\Eb,\F]$ is
an $\Eb$ category together with a chosen element of
$\Ff\left(\Lim[\delta]\right)$, where $\delta$ is the
description of $\F$ as in~\ref{termin}.

The functor $\Constants$ induces a forgetful functor (it
forgets the chosen element of
$\Ff\left(\Lim[\delta]\right)$) from set-valued (or more
general) models of $\SynCat[\Eb,\F]$ to models of $\Eb$.

Each object of $\SynCat[\Eb,\F]$ is the limit of a diagram
in $\E[\Fa]$ by Lemma~\ref{piessens}, and the diagram is in
some sense a description of a possible construction in any
model of the form $\F$.\erem

\exam\label{skexam1} As an example, consider the
finite-limit sketch $S$ with graph
\begin{equation}
\xymatrix{ A\ar@<1ex>[r]^{v} &B\ar@<1ex>[l]^{u}\ar[r]^{f}
&C
\\
}
\end{equation}
one diagram
\begin{equation}
\xymatrix{ B\ar[r]^{u}\ar[d]_{Id[B]} &A\ar[dl]^{v}
\\
B
\\
}
\end{equation}
and one cone
\begin{equation}
\xymatrix{ &B\ar[dl]_{u}\ar[dr]^{f} &
\\
A &&C
\\
}
\end{equation}
This is in fact a finite-product sketch, but any such
sketch is also a finite-limit sketch. Then one way to
capture the information in the sketch is to take $\delta$
to be the following diagram in $\FinLim$ and define the
$\FinLim$ form $\F$ determined by a freely adjoined
constant $\Fa:1\to\Lim[\delta]$.

{\wider\wider
\begin{equation}
\xymatrix{ &&\ob^{A}\x \ob^{C}\ar[d]^{\prod} &&
\\
&\arr^{u}\ar[dl]_{\target}
&\cone\ar[l]_{\lproj}\ar[r]^{\rproj}
&\arr^{f}\ar[r]^{\target}\ar[d]^{\source} &\ob^{C}
\\
\ob^{A} &\arr_{2}\ar[u]|{\rfac}\ar[r]|{\comp}\ar[d]|{\lfac}
&\arr &\ob^{B}\ar[l]|{\unit} &
\\
&\arr^{v}\ar[ul]^{\source}\ar[urr]|{\target} &&&
\\
}
\end{equation}
} \eexam

\section{Theories and models of forms}\label{modthform}\label{catthform}
For completeness, we define models of forms and their
morphisms, and theories of forms, using the notation of
this monograph, but only briefly since these ideas are not used
in this monograph. More detail and examples may be found in
\cite{gensk} or \cite{powerwells}.

\rem We are in this section identifying a model of $\Eb$
with an actual category with structure imposed by $\Eb$,
and similarly for models of $\SynCat[\Eb,\F]$.  This is an
example of the phenomenon mentioned in
Sections~\ref{modrem1} and~\ref{moremodrem1}. See also
Section~\ref{modrem2}.\erem

Now let $\F$ be an $\Eb$-form with description
$\delta:1\to\Eb$, so that it is named by
$\Fa:1\to\Lim[\delta]$. Then $\SynCat[\Eb,\F]$ is a
finite-limit theory and so has an initial model.

\defn\label{catthform} The initial model of $\SynCat[\Eb,\F]$ is
called the {\bf $\Eb$-theory} of $\F$, denoted by
$\CatTh[\Eb,\F]$.\edefn

\rem The form $F$ is an element of the value of
$\Lim[\delta]$ in $\CatTh[\Eb,\F]$.\erem

\rem\label{finlimflth} Once a finite-limit sketch $S$ is
captured as a form $\F$ as described in
Example~\ref{skexam1}, it follows that $\FinLimTh[\S]$ is
naturally equivalent as a category to
$\CatTh[\FinLim,\F]$.\erem

\defn\label{fmoddef} A {\bf model of $\F$ in an $\Eb$-category
$\Csc$} is defined to be a model of $\SynCat[\Eb,\F]$ with
underlying $\Eb$-category $\Csc$.\edefn

This means that $\Csc$ is the value of the functor from
models of $\SynCat[\Eb,\F]$ to models of $\Eb$ that forgets
the element corresponding to $\F$ ---  see
Section~\ref{forgetful}.

Let $\M$ be a model of $\F$ with underlying $\Eb$-category
$\Csc$.  Since $\CatTh[\Eb,\F]$ is the initial model of
$\F$, there is a unique $\Eb$-functor
$\phi:\CatTh[\Eb,\F]\to\Csc$  that takes
$\CatTh[\Eb,\F](\Fa)$ to $\M(\Fa)$. In the case of familiar
sketches, say finite-limit or finite-product sketches
(corresponding to $\Eb=\FinLim$ and $\Eb=\FinProd$
respectively), that functor $\phi$ is what would usually be
called the functor from the theory induced by a model of
the sketch.  To define for forms the entities that
correspond in those cases to the actual sketch and its
models involves complications and is carried out in two
different ways in \cite{gensk} and \cite{powerwells}.

Finally, a {\bf morphism of models} of a form $\F$ in a
category $\Csc$ is simply a natural transformation between
the functors $\phi$ and $\phi'$ corresponding as described
in the previous paragraph to models $\M$ and $\M'$ in
$\Csc$.

Any work making extensive use of the entities constructed
in this section will probably need to introduce more
elaborate terminology and notation (for example for $\phi$)
than is used here.  Such refinements, however, are peripheral to our
concerns and we hope that this discussion is sufficiently
detailed to obviate confusion.

\section{Relationship between forms and sketches}\label{modrem2}
This section continues the discussion begun in
Sections~\ref{modrem1} and~\ref{moremodrem1}. One
constructor space is $\FinLim$, defined in
Section~\ref{finlimsk}. A finite-limit sketch $S$ in the
traditional sense (a graph with diagrams and cones)
corresponds to a $\FinLim$-form  in the construction  using
the methods of Example~\ref{skexam1} and has the same
models. The traditional finite-limit sketch $S$ is an
element of the value in the initial model of a certain node
$\vsfi$ of $\FinLim$ (not uniquely determined by $S$) which
is the limit of a generally large and complicated diagram
$\delta$ (not uniquely determined either by $S$ or by
$\vsfi$) in $\FinLim$ that specifies the graph, diagrams
and cones of $S$.

These remarks apply to other types of Ehresmann sketches by
replacing $\FinLim$ by the suitable constructor space.


\chapter{Examples of sketches for constructor
spaces}\label{conspsk}

Here we present constructor space sketches for certain
types of categories. In each case the models are
categories of the sort described and the morphisms of
models are functors that preserve the structure on the
nose. It is an old result that such categories can be
sketched. See \cite{burr}, \cite{burr70},
\cite{macst2}, and \cite{coppeylair88}, for example.

The embedding $\CatStruc[\Eb]$ of Section~\ref{cssdef}
will in each case be inclusion.

\section{Notation}\label{convia} We denote the $i$th
projection in a product diagram of the form
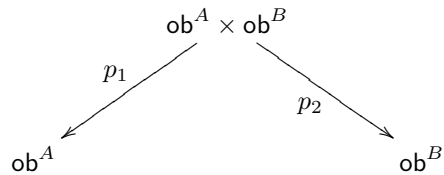
\begin{equation*}
\xymatrix{ &\ob^{A}\x \ob^{B} \ar[dl]_{p_1}
\ar[dr]_{p_2}
&\\
\ob^{A} & &\ob^{B}
\\
}
\end{equation*}
as $p_i$, or $p_i^{A\xord B}$ if the source or target
is not shown. We use a similar device for the product
of three copies of $\ob$.

%
%

\section{The sketch $\protect\Catsk$ for categories}\label{realcatsk}
This version of the sketch for categories is based on
\cite{ctcs}. Another version is given in
\cite{coppeylair88}, page~64. The first versions were
done by Ehresmann~\shortcite{ehrkan}, \shortcite{csg}
and~\shortcite{ehres68}.

\subsection{The graph of $\Catsk$}\label{grphfl}

The graph of the sketch for categories contains nodes
as follows.
\begin{enumerate}
\item $\one$, the formal terminal object. \item $\ob$,
the formal set of objects. \item $\arr$, the formal
set of arrows. \item $\arr_2$, the formal set of
composable pairs of arrows. \item $\arr_3$, the formal
set of composable triples of arrows. \end{enumerate}
The arrows for the sketch for categories are
\begin{enumerate}
\item $\unit:\ob\to\arr$ that formally picks out the
identity arrow of an object. \item
$\source,\target:\arr\to\ob$ that formally pick out
the source and target of an arrow. \item
$\comp:\arr_2\to\arr$ that picks out the composite of
a composable pair. \item $\lfac,\rfac:\arr_2\to\arr$
that pick out the left and right factors in a
composable pair. \item
$\lfac,\mfac,\rfac:\arr_3\to\arr$ that pick out the
left, middle and right factors in a composable triple
of arrows. \item $\lass,\rass:\arr_3\to\arr_2$:
$\lass$ formally takes $<h,g,f>$ to $<h\oc g,f>$ and
$\rass$ takes it to $<h,g\oc f>$. \item
$\lunit,\runit:\arr\to\arr_2$: $\lunit$ takes an arrow
$f:A\to B$ to $<\Id[B],f>$ and $\runit$ takes it to
$<f,\Id[A]>$. \item Arrows $\id:x\to x$ as needed.
\end{enumerate}
Observe that $\id$, $\lfac$ and $\rfac$, like $p_1$
and $p_2$, are overloaded. We will observe the same
care with these arrows as with $p_1$ and $p_2$ as
mentioned in Section~\ref{convia}.

\subsection{Cones of $\Catsk$}

{\wider\taller

$\arr_2$ and $\arr_3$ are defined by these cones:


$$\begin{array}{c@{\hspace{2em}}c}
\xymatrix{ & \arr_{2}\ar[dl]_{\lfac}\ar[dr]^{\rfac}
&\\
\arr\ar[dr]_{\source} & &\arr\ar[dl]^{\target}
\\
&\ob
&\\
} & \xymatrix{ &
\arr_{3}\ar[dl]_{\lfac}\ar[dr]^{\rfac}\ar[d]|{\mfac}
&\\
\arr\ar[d]_{\source}
&\arr\ar[dr]|{\source}\ar[dl]|{\target}
&\arr\ar[d]^{\target}
\\
\ob& &\ob
\\
}
\end{array}$$
}

\subsection{Diagrams of $\Catsk$}

{\wider\taller\taller

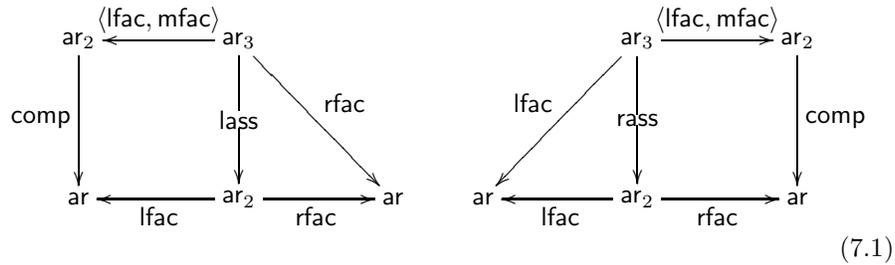
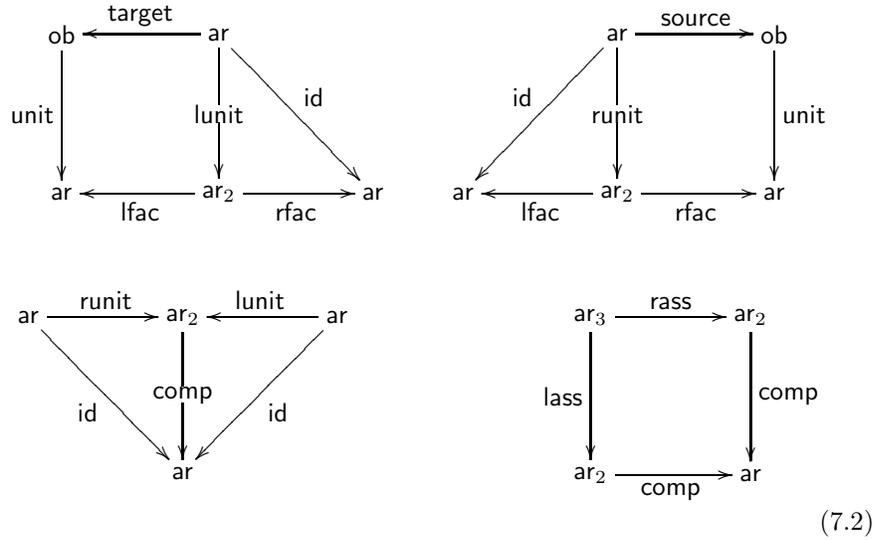
\begin{equation}\label{lass&rass}
\begin{array}{c@{\hspace{2em}}c}
\xymatrix{ \arr_{2} \ar[d]_{\comp}
&\arr_{3}\ar[l]_{<\lfac,\mfac>} \ar[dr]^{\rfac}
\ar[d]|{\lass}
&\\
\arr &\arr_{2} \ar[l]^{\lfac} \ar[r]_{\rfac} &\arr
\\
} & \xymatrix{ &\arr_{3} \ar[r]^{<\lfac,\mfac>}
\ar[dl]_{\lfac} \ar[d]|{\rass} &\arr_{2}\ar[d]^{\comp}
\\
\arr &\arr_{2} \ar[l]^{\lfac} \ar[r]_{\rfac} &\arr
\\
}
\end{array}
\end{equation}

$$\begin{array}{c@{\hspace{2em}}c}
\xymatrix{ \ob \ar[d]_{\unit} &\arr\ar[l]_{\target}
\ar[dr]^{\id} \ar[d]|{\lunit}
&\\
\arr &\arr_{2} \ar[l]^{\lfac} \ar[r]_{\rfac} &\arr
\\
} & \xymatrix{ &\arr \ar[r]^{\source} \ar[dl]_{\id}
\ar[d]|{\runit} &\ob\ar[d]^{\unit}
\\
\arr &\arr_{2} \ar[l]^{\lfac} \ar[r]_{\rfac} &\arr
\\
}
\end{array}
$$

\begin{equation}\label{u&ass}
\begin{array}{c@{\hspace{2em}}c}
\xymatrix{ \arr \ar[r]^{\runit} \ar[dr]_{\id}
&\arr_{2} \ar[d]|{\comp} &\arr \ar[l]_{\lunit}
\ar[dl]^{\id}
\\
&\arr & &
\\}
& \xymatrix{ \arr_{3} \ar[r]^{\rass} \ar[d]_{\lass}
&\arr_{2} \ar[d]^{\comp}
\\
\arr_{2} \ar[r]_{\comp} &\arr
\\
}
\end{array}
\end{equation}
}

\section{The sketch for the constructor space
$\FinProd$}\label{finprodsk}

To get the sketch for categories with finite products,
we must add the following nodes and arrows to the
sketch for categories:

Nodes:
\begin{enumerate}
\item $\ta$, the formal set of terminal arrows. \item
$\cone$, the formal set of cones of the form
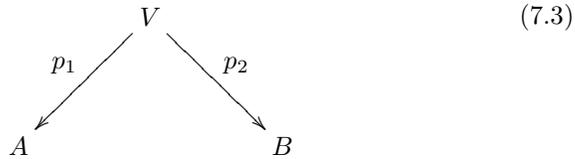
\begin{equation}
\xymatrix{ &V \ar[dl]_{p_{1}} \ar[dr]^{p_{2}} &
\\
A & &B
\\
}
\end{equation}
\item $\fid$, the formal set of fill-in diagrams
(``sawhorses'') of the form
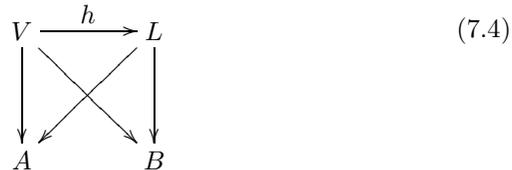
\begin{equation}\label{sawhorse}
\xymatrix{ V \ar[r]^{h} \ar[dr] \ar[d] &L \ar[dl]
\ar[d]
\\
A &B
\\
}
\end{equation}
where $h$ commutes with the cone projections.
\end{enumerate}
Arrows:\begin{enumerate} \item $\ter:1\to\ob$, that
formally picks out a particular terminal object. \item
$!:\ob\to\ta$, that picks out the arrow from an object
to the terminal object. \item $\inc:\ta\to\arr$, the
formal inclusion of the set of terminal arrows into
the set of arrows. \item $\prod:\ob\x\ob\to\cone$,
that picks out the product cone over a pair of
objects. \item $\soco:\fid\to\cone$, that picks out
the source cone of a fill-in arrow. \item
$\taco:\fid\to\cone$, that picks out the target cone
of a fill-in arrow. \item $\ufid:\cone\to\fid$, that
takes a cone to the unique fill-in diagram that has
the cone as source cone. \item $\fia:\fid\to\arr$ that
formally picks out the fill-in arrow in a fill-in
diagram.
\end{enumerate}

\subsection{Cones for $\FinProd$}\label{bpcones}
$\FinProd$ has four cones in addition to those of the
sketch for categories. One is the cone  $$1$$
over the empty diagram.  The  one below says that $\ta$
is the formal set of arrows to the terminal object:
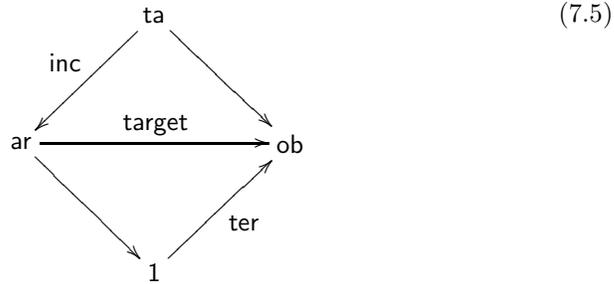
\begin{equation}\label{tacone}
\xymatrix{
& \ta \ar[dl]_\inc \ar[dr] &\\
\arr \ar[rr]^\target \ar[dr] && \ob\\
& \one \ar[ur]_\ter\\
}
\end{equation}
Note that in giving this cone, we are not only saying
that $\ta$ is the limit of the diagram
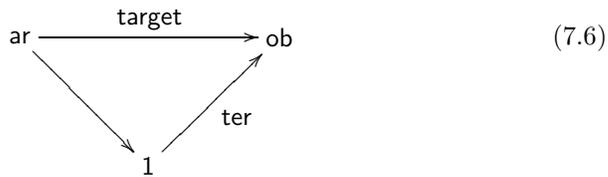
\begin{equation}
\xymatrix{
\arr \ar[rr]^\target \ar[dr] && \ob\\
& \one \ar[ur]_\ter\\
}
\end{equation}
but also that $\inc$ is one of the projection arrows.
(Indeed, this is the only projection arrow that
matters, since the other two are induced.)

The following cone makes $\cone$ the formal object of
cones to a discrete diagrams consisting of a pair of
objects.
 {\wider\wider\taller
\begin{equation}\label{bindi} \xymatrix{
& \cone \ar[dl]|{\lproj} \ar[dr]|{\rproj} \\
\arr \ar[r]_{\source} & \ob & \ar[l]^{\source} }
\end{equation}}

Finally, there must be a cone with vertex $\fid$ over
Diagram~\ref{bpbigdiag} below, which is annotated to
refer to Diagram~(\ref{sawhorse}), in which $\Gamma$
is the cone with vertex $V$, $\Lambda$ is the cone
with vertex $L$ and is a limit cone, and $h$ is the
fill-in arrow. In this case, the projection arrows of
the cone are not shown. {\taller\taller\wider\wider 
\begin{equation}\label{bpbigdiag}
\xymatrix{ \arr_2 \ar[d]|{\comp} \ar `u[r]
`[rrrr]|{\rfac} [rrrrdd] \ar `l[d] `[ddd]|(0.6){\lfac}
[ddd]
&&\ob^V &&\\
\arr \ar[d]|{\target} \ar[urr]|{\source} &
\cone^{\Gamma} \ar[l]^{\lproj} \ar[r]_{\rproj}
& \arr \ar[d]|{\target} \ar[u]|{\source} &&\\
\ob^A & \ob^A\xord\ob^B \ar[l]_{p_1} \ar[r]^{p_2}
\ar[d]|{\prod} & \ob^B & \arr_2 \ar[ul]|(0.4){\comp}
\ar[dl]|(.4){\lfac} \ar[r]^{\rfac}
& \arr^h \ar[uull]|{\source} \ar[ddll]|{\target}\\
\arr \ar[u]|{\target} \ar[drr]|{\source} &
\cone^{\Lambda} \ar[l]_{\lproj}
\ar[r]^{\rproj} & \arr \ar[u]|{\target} \ar[d]|{\source}\\
 & &\ob^L\\
}
\end{equation}}%
In addition, we require: \sqlist \item The projection
to $\cone^\Gamma$ must be $\soco$. \item The
projection to $\cone^{\Lambda}$ must be $\taco$. \item
The projection to $\arr^h$ must be $\fia$. \esqlist

\subsection{Diagrams for $\FinProd$}
The following two diagrams make the arrow to the
terminal object have the correct source and target.
{\wider\wider\taller
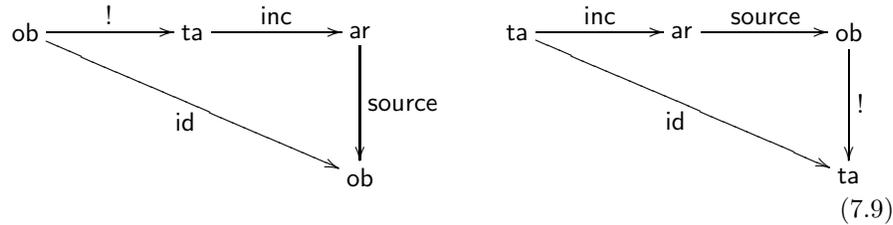
\begin{equation}\label{terdiags}\begin{array}%
{c@{\hspace{2em}}c} \xymatrix{ {\ob} \ar[r]^{!}
\ar[rrd]_{\id} & {\ta} \ar[r]^{\inc} &
{\arr} \ar[d]^{\source}  \\
&& {\ob} \\ }  & \xymatrix{ {\ta} \ar[r]^{\inc}
\ar[drr]_{\id} & {\arr} \ar[r]^{\source} &
{\ob} \ar[d]^{!}\\
&&{\ta} \\ }
\end{array}\end{equation}
}
The diagram below makes the fill-in arrow to a product
unique.
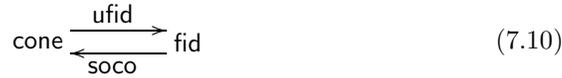
\begin{equation}
\xymatrix{
\cone \ar[r]<1ex>^{\ufid} & \fid \ar[l]<1ex>^{\soco} \\
}\end{equation} The  diagram below forces the product
cone projections to have the correct targets.
{\taller\wider\wider
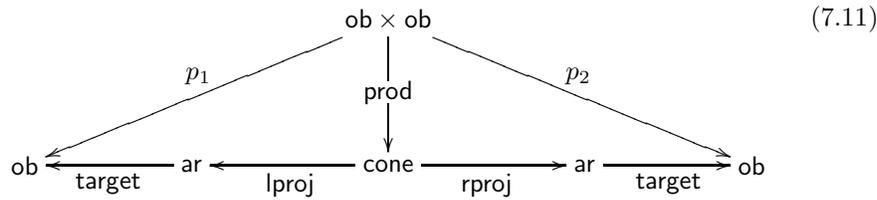
\begin{equation}\label{prodtars}
\xymatrix{
&& \ob\x\ob \ar[dll]_{p_1} \ar[drr]^{p_2} \ar[d]|{\prod}  \\
\ob & \arr \ar[l]^{\target} & \cone \ar[l]^{\lproj}
\ar[r]_{\rproj} & \arr \ar[r]_{\target} & \ob\\
}
\end{equation}}

\section{Modules}\label{modsubsec}
As we proceed to sketch more complicated
constructions, we will need to use some device to
communicate the nature of the necessary diagrams,
which become too large to comprehend easily. Here we
introduce the first of several {\bf modules}: diagrams
that occur frequently as subdiagrams  because they are
needed to force the value of a node in a model to
contain certain types of constructions.

Modules are a well-understood part\mpark{Added paragraph}
of programming language methodology.  We believe that the concept
called ``module'' here can be made explicit enough to become part of a programming
language based on the techniques of this monograph, but that work is yet to be accomplished.

\subsection{The module for the product of objects}
Every occurrence of $\ob$ that is annotated $M\x N$
must be part of a subdiagram of the following form:
{\wider\taller
\begin{equation}\label{MxNdiag}\xymatrix{
\ob^M & \arr^{p_1^{M\xord N}} \ar[l]_{\target} \ar[dr]^{\source}\\
\ob^M\x\ob^N \ar[r]^<>(.5){\prod} \ar[u]^{p_1}
\ar[d]_{p_2}
& \cone \ar[u]|{\lproj} \ar[d]|{\rproj} & \ob^{M\xord N}\\
\ob^N & \arr^{p_2^{M\xord N}}\ar[l]^{\target}
\ar[ur]_{\source} }
\end{equation}
}%
Henceforth, an occurrence of $\ob$ annotated $A\x B$
(for example) will be taken to \emph{imply} the
existence of a subdiagram of the form of
Diagram~(\ref{MxNdiag}) with $M$ replaced with $A$ and
$N$ replaced with $B$.  The subdiagram will not
necessarily be shown. If this is part of a diagram
$\delta$, the diagram can be reconstructed by taking
the union of the shape graph of the
module~(\ref{MxNdiag}) and the shape graph of the part
of $\delta$ that is shown on the page, and defining
the diagram based on the resulting graph as the
pushout of the diagram shown and the module.  This is
illustrated in Diagrams~(\ref{crossdiag})
and~(\ref{crossdiag2}) in the next section.

\subsection{The module for the product of
arrows}\label{parrmod}

In the commutative diagram
\begin{equation}\label{ucrossa}
\xymatrix{ K \ar[dd]_u & K\x N \ar[l]_(0.5){p_1}
\ar[dd] \ar[dr]^{p_2}\\
&&N\\
M & M\x N \ar[l]_(0.5){p_2} \ar[ur]_{p_2} }
\end{equation}
the unlabeled arrow is necessarily $u\x \Id[N]:K\x
N\to M\x N$. Such a diagram must be an element in a
model of the value of Diagram~(\ref{crossdiag}) below,
which is therefore a module for the product of an
arrow and an identity arrow.  In this diagram,
$\phi\ceq <p_2^{M\x N}, u\x\Id[N]>$.
{\narrower\narrower\narrower\taller\taller
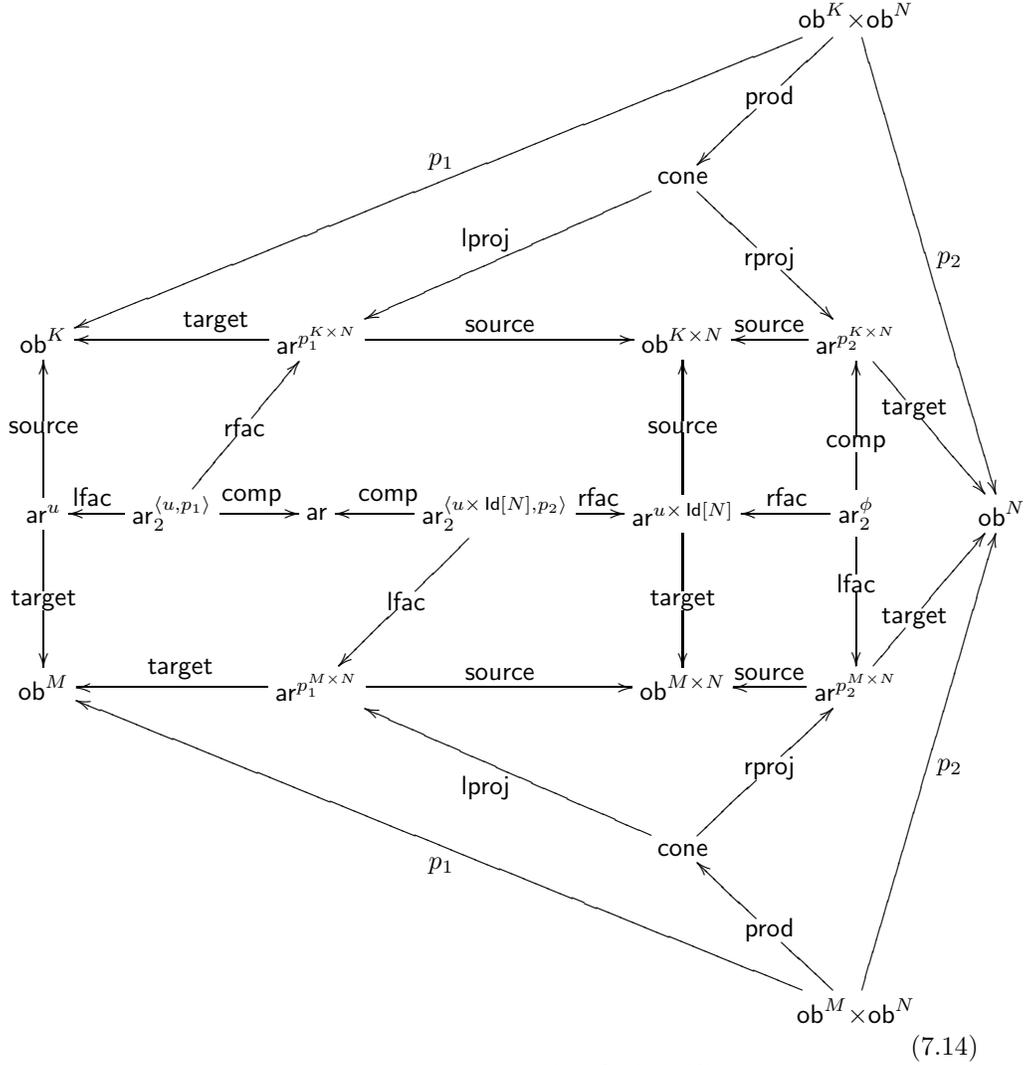
\begin{equation}\label{crossdiag}
\hspace{-.3in}\xymatrix{ &&&&& \ob^K\xord\ob^N
\ar[2,-5]_{p_1} \ar[dl]|{\prod}
\ar[3,1]^{p_2}\\
&&&& \cone \ar[dll]_{\lproj}
\ar[dr]|{\rproj}\\
\ob^K && \arr^{p_1^{K\x N}} \ar[ll]_(0.37){\target}
\ar[rr]^(0.5){\source} && \ob^{K\x N} & \arr^{p_2^{K\x
N}} \ar[l]_{\source} \ar[dr]|(0.4){\target}
\\
\arr^u \ar[u]|{\source} \ar[d]|{\target} &
\arr_2^{<u,p_1>} \ar[l]_(0.62){\lfac} \ar[ur]|{\rfac}
\ar[r]^(0.55){\comp} & \arr &
\arr_2^{<u\xord\Id[N],p_2>} \ar[l]_(0.6){\comp}
\ar[r]^(0.55){\rfac} \ar[dl]|{\lfac} &
\arr^{u\xord\Id[N]} \ar[d]|{\target} \ar[u]|{\source}
& \arr_2^{\phi} \ar[l]_(0.4){\rfac}
\ar[u]|(0.4){\comp} \ar[d]|(0.4){\lfac} & \ob^N
\\
\ob^M && \arr^{p_1^{M\x N}} \ar[ll]_{\target}
\ar[rr]^{\source} && \ob^{M\x N} & \arr^{p_2^{M\x N}}
\ar[l]_(0.5){\source} \ar[ur]|(0.4){\target}
\\
&&&& \cone \ar[ull]^{\lproj} \ar[ur]|{\rproj}
\\
&&&&& \ob^M\xord\ob^N \ar[-2,-5]^{p_1} \ar[ul]|{\prod}
\ar[-3,1]_{p_2}
\\
}
\end{equation}
}%
More precisely, let $x$ be an element of
$\M(\Lim[D(\ref{crossdiag})])$, for some category $\M$
with finite limits. Then if
$\Proj[\Lim[D(\ref{crossdiag})],h](x)=h$, then
$$\Proj[\Lim[D(\ref{crossdiag})],h\xord\Id[A]](x)=h\xord\Id[A]$$ as suggested by the
notation. This will be used in Section~\ref{cccsec}
below.

Diagram~(\ref{crossdiag}) contains two copies of
Diagram~(\ref{MxNdiag}), the module for the product of
two objects.  The copy at the bottom is precisely
Diagram~(\ref{MxNdiag}), and the copy at the top is
Diagram~(\ref{MxNdiag}) with $M$ replaced with $K$. In
the sequel, a diagram such as
Diagram~(\ref{crossdiag}) will be drawn without the
modules, as shown below.
{\taller\taller\narrower\narrower\narrower\narrower
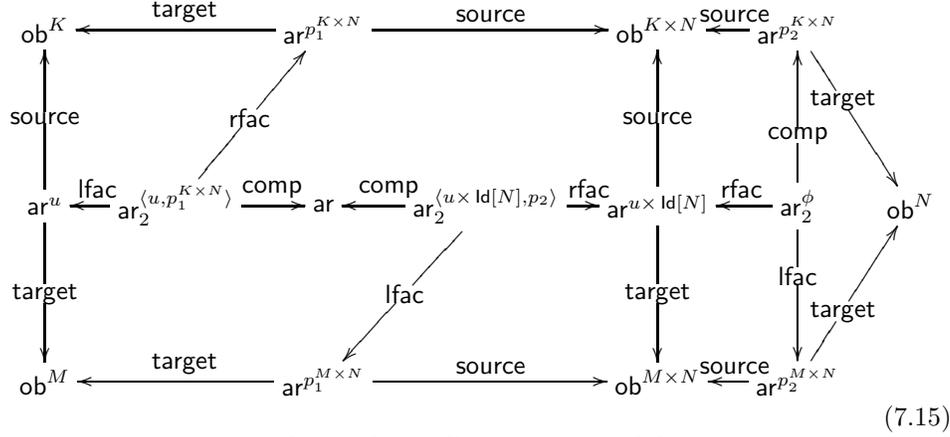
\begin{equation}\label{crossdiag2}
\hspace*{-.2in}\xymatrix{ \ob^K && \arr^{p_1^{K\xord
N}} \ar[ll]_{\target} \ar[rr]^(0.5){\source} &&
\ob^{K\xord N} & \arr^{p_2^{K\xord N}}
\ar[l]_(0.45){\source} \ar[dr]|(0.4){\target}
\\
\arr^u \ar[u]|{\source} \ar[d]|{\target} &
\arr_2^{<u,p_1^{K\xord N}>} \ar[l]_(0.6){\lfac}
\ar[ur]|{\rfac} \ar[r]^(0.65){\comp} & \arr &
\arr_2^{<u\xord\Id[N],p_2>} \ar[l]_(0.6){\comp}
\ar[r]^(0.6){\rfac} \ar[dl]|{\lfac} &
\arr^{u\xord\Id[N]} \ar[d]|{\target} \ar[u]|{\source}
& \arr_2^{\phi} \ar[l]_(0.4){\rfac}
\ar[u]|(0.4){\comp} \ar[d]|(0.4){\lfac} & \ob^N
\\
\ob^M && \arr^{p_1^{M\x N}} \ar[ll]_{\target}
\ar[rr]^{\source} && \ob^{M\x N} & \arr^{p_2^{M\x N}}
\ar[l]_(0.45){\source} \ar[ur]|(0.4){\target}
\\
}
\end{equation}
}%
Diagram~(\ref{crossdiag}) may be mechanically
reconstructed from Diagram~(\ref{crossdiag2}) and the
annotations that include the symbols $M\x N$ and $K \x
N$ (three of each).  The shape diagram of
Diagram~(\ref{crossdiag}) is the pushout of the shape
diagram of Diagram~(\ref{crossdiag2}) and the shape
diagrams of the modules Diagram~(\ref{MxNdiag})  and
Diagram~(\ref{MxNdiag}) with $M\leftarrow K$.  Each of
the latter two have four annotated nodes and six
annotated arrows in common with
Diagram~(\ref{crossdiag2}), and the values of any two
of the three smaller diagrams at a given common node
or arrow is of course the same, so that
Diagram~(\ref{crossdiag}) is the union of
Diagram~(\ref{crossdiag2}) and the two modules.

\section{The sketch
for the constructor space $\FinLim$}\label{finlimsk}

We sketch the constructor space $\FinLim$ by adding
data to the sketch for $\FinProd$ that ensure that a
$\FinLim$-category has equalizers of pairs of arrows.
The sketch has the following nodes:
\begin{enumerate}
\item $\ppair$ is the formal set of parallel pairs of
the form
\begin{equation}\label{ppa}\xymatrix{
A \ar[r]<1ex>^f \ar[r]<-1ex>_g  & B\\
}\end{equation} \item $\econe$ is the formal set of
diagrams\mpar{I couldn't read what you wrote on the bottom of page 44.}
\begin{equation}\label{econes}\xymatrix{
E \ar[r]^u & A \ar[r]<1ex>^f \ar[r]<-1ex>_g  & B\\
}\end{equation} in which $f\oc u= g\oc u$.  Of course, a
cone to Diagram~(\ref{ppa}) also has a projection to
$B$, but that is forced and need not be included in
the data for the cone. \item $\efid$ is the set of
fill-in diagrams {\taller\begin{equation}\xymatrix{
& X \ar[dl]_v \ar[d]^u\\
E \ar[r]^e & A \ar[r]<.7ex>^f \ar[r]<-1ex>_g  & B\\
}\end{equation}} in which $f\oc e=g\oc e$ and $u=e\oc v$.
\end{enumerate}
The arrows of the sketch include:
\begin{enumerate}
\item $\equ:\ppair\to\econe$, that formally picks out
the equalizer of the parallel pair. \item
$\top,\bot:\ppair\to\arr$, that pick out $f$ and $g$
in Diagram~(\ref{ppa}). \item
$\etop,\ebot:\econe\to\arr$, that pick out $f$ and $g$
in Diagram~(\ref{econes}). \item
$\esoco,\etaco:\efid\to\econe$ that pick out the
source and target cones of the fill-in arrow. \item
$\eufid:\econe\to\efid$ that takes a diagram of the
form of Diagram~(\ref{econes}) to the unique fill-in
diagram that has this diagram as source cone. \item
$\efia:\efid\to\arr$ that picks out the fill-in arrow
in a fill-in diagram. \end{enumerate}

\subsection{Cones for $\FinLim$}\label{finlimcones}
$\ppair$ is the limit of the diagram
{\wider\taller\taller\begin{equation}\xymatrix{ \ar^f
\ar[r]^{\target} \ar[d]|{\source}
& \ob^B\\
\ob^A
& \ar^g \ar[l]_{\source} \ar[u]|{\target}\\
}\end{equation}}%
The following projections from $\ppair$ have names:
$\soob:\ppair\to \ob^A$, $\top:\ppair\to\arr^f$, and
$\bot:\ppair\to\arr^g$.

$\econe$ is the limit of
{\wider\wider\taller\begin{equation}\label{econe2}
\xymatrix{ & \arr_2 \ar[d]|{\rfac} \ar[r]^{\lfac}
\ar[dl]_{\comp}
& \arr^f \ar[d]|{\source} \ar[dr]^{\target}\\
\ar & \arr^e \ar[r]^{\target} & \ob^A
& \ob^B\\
& \arr_2 \ar[ul]^{\comp} \ar[r]_{\lfac} \ar[u]|{\rfac}
& \arr^g \ar[u]|{\source}
\ar[ur]_{\target}\\
}
\end{equation}}%
Two projections have names: $\etop:\econe\to\arr^f$
and $\ebot:\econe\to\arr^g$.

$\efid$ is the limit of the pushout of
Diagram~(\ref{econe2}) and the following diagram. Note
that the common part of the two diagrams is
$$\xymatrix{ \arr^e \ar[r]^t & \ob^A}$$
\rem We could
have presented Diagram~(\ref{bpbigdiag}) as a pushout
in much the same way (the common part would describe
the arrow $h:V\to L$).  We have deliberately varied
the way we present the data in this monograph because
we are not sure ourselves which approach communicates
best. {\wider\wider\taller\begin{equation} \xymatrix{
\arr^v\ar[dd]_{\target}\ar[r]^{\source} & \ob^X &
\arr^u \ar[l]_{\source}
\ar[dd]^{\target}\\
& \arr_2 \ar[ul]|{\rfac}
\ar[ur]|{\comp} \ar[d]|{\lfac} \\
\ob^E & \arr^e \ar[l]_{\source} \ar[r]^{\target}
& \ob^A\\
}
\end{equation}}%
The named projections are $\esoco:\efid\to\arr^u$,
$\etaco:\efid\to\arr^e$ and $\efia:\efid\to\arr^v$.\erem

\subsection{Diagrams for $\FinLim$}\label{finlimdiags}
The following diagram makes the fill-in arrow unique.
{\wider\wider
\begin{equation}
\xymatrix{
\econe \ar[r]<1ex>^{\eufid} & \efid \ar[l]<1ex>^{\esoco} \\
}\end{equation}
These two diagrams ensure that the equalizer cone be a
cone to the correct diagram.
\begin{equation}
\begin{array}{cc}
\xymatrix{ \ppair \ar[r]^{\equ} \ar[dr]_{\top}
& \econe \ar[d]^{\etop}\\
& \arr\\
} & \xymatrix{ \ppair \ar[r]^{\equ} \ar[dr]_{\bot}
& \econe \ar[d]^{\ebot}\\
& \arr \\
}
\end{array}
\end{equation}

\section{The sketch for the constructor space
$\CCC$}\label{cccsec}

\defn\label{cccdefn}
A {\bf Cartesian closed category} is a category $\Csc$
with the following structure: \blist{CCC} \item $\Csc$
has binary products. \item For each pair of objects
$A$ and $B$ of $\Csc$, there is an object $B^A$ and an
arrow $\eval:B^A\x A\to B$. \item For each triple of
objects $A$, $B$ and $C$ of $\Csc$, there is a map
\begin{equation}\label{lambdadef}
\lambda:\Hom(B\x A, C)\to \Hom(B,C^A)\end{equation}
such that for every arrow $f:B\x A\to C$,
{\wider\wider\wider\taller
\begin{equation}\label{evaleq}
\xymatrix{ B\xord A \ar[r]^{\lambda f\xord\Id[A]}
\ar[dr]_f
& C^A \xord A \ar[d]|{\eval}\\
& C\\
}
\end{equation}}%
commutes. \item\label{fourthreq} For any arrow $g:B\to
C^A$, $\lambda(\eval\oc (g\xord\Id[A]))=g$.
\elist\edefn Using this definition, the sketch for the
constructor space for Cartesian closed categories may
be built on the sketch for $\FinProd$ by adding the
following nodes and arrows.

The nodes are:
\begin{enumerate}
\item $\twovf$, the formal set of ``functions of two
variables'', that is, arrows of the form $B\x A\to C$.
\item $\curry$, the formal set of ``curried
functions'' $B\to C^A$.
\end{enumerate}
The sketch for $\CCC$ has arrows
\begin{enumerate}
\item $\fs:\ob^B\x\ob^A\to\ob^{B^A}$ that picks out
the function space $B^A$ of two objects $B$ and $A$.
\item $\ev:\ob^B\x\ob^A\to\arr$ that picks out the
arrow $\eval:B^A\x A\to B$. \item\label{lamitem}
$\lam:\arr\to\arr$, the formal version of the mapping
$\lambda$ of Diagram~(\ref{lambdadef}). \item
$\tsource:\twovf\to\ob^{B\xord A}$, that picks out the
source of a function $f:A\x B\to C$. \item
$\ttarget:\twovf\to\ob^C$, that picks out the target
of a function $f:A\x B\to C$. \item
$\arrow:\twovf\to\arr^f$, that picks out the arrow $f$
itself. \item $\csource:\curry\to\ob^B$, that picks
out the source of a curried function $g:B\to C^A$.
\item $\ctarget:\curry\to\ob^{C^A}$, that picks out
the target of a curried function $f:B\to C^A$. \item
$\arrow:\curry\to\arr^g$, that picks out the arrow $g$
itself. \end{enumerate}

\subsection{Cones for $\CCC$}\label{cccc}
The constructor space $\CCC$ must have two cones
{\wider\taller
\begin{equation}\label{twovfcone}
\xymatrix{ & \twovf \ar[dl]_{\tsource} \ar[d]|{\arrow}
\ar[dr]^{\ttarget}
\\
\ob^{B\x A} & \arr^f \ar[l]^{\source} \ar[r]_{\source}
& \ob^C
\\
}\end{equation}
}%
{\taller\taller\taller\Narrower
\begin{equation}\label{currycone}\xymatrix{
&&& \curry \ar[1,-3] \ar[dl]|{\ctarget}
\ar[dr]|{\arrow} \ar[1,3]^{\csource}
\\
\ob^C\xord \ob^A \ar[rr]_{\fs} && \ob^{C^A} && \arr^g
\ar[ll]^{\target} \ar[rr]_{\source} && \ob^B
\\
}
\end{equation}
}

\subsection{The module for function
spaces}\label{fsmodule} Henceforth, we will assume the
module
\begin{equation}\label{fsmodulediag}\xymatrix{
\ob^M & \ob^M\xord \ob^N \ar[l]_{p_1} \ar[r]^{\fs}
\ar[d]^{p_2} & \ob^{M^N}
\\
& \ob^N\\
}\end{equation} is attached whenever an occurrence of
$\ob$ is annotated by $M^N$.  Note that this occurred
in Diagram~(\ref{currycone}).

\subsection{Diagrams for $\CCC$}

Diagram~(\ref{evaldiag}) below forces $\eval$ to have
the correct domain and codomain. {\wider\wider
\begin{equation}\label{evaldiag}
\xymatrix{ \ob^{B^A\xord A} & \arr^{\eval}
\ar[l]_{\source} \ar[r]^{\target} & \ob^B }
\end{equation}
}%

Expanding this diagram using the required modules is a
two stage process, giving {\narrower\begin{equation}
\xymatrix{ \ob^A & \ob^{B^A}\xord\ob^A
\ar[l]_(0.5){p_2} \ar[r]^(0.5){p_1} \ar[d]|{\prod} &
\ob^{B^A}
\\
\arr \ar[u]|{\rproj} \ar[dr]_{\source} & \cone
\ar[l]_{\lproj} \ar[r]^{\rproj} & \arr \ar[u]|{\lproj}
\ar[dl]^{\source} & \ob^B\xord\ob^A \ar[ul]_{\fs}
\ar[r]^(0.5){p_1} \ar[d]^{p_2} & \ob^A
\\
& \ob^{B^A\xord A} & \arr^{\eval} \ar[l]^{\source}
\ar[r]_{\target} &
\ob^B\\
}\end{equation}}

Diagram~(\ref{lfst}) below forces $\lambda f$ to have
the correct domain and codomain. {\wider\wider\taller
\begin{equation}\label{lfst}\xymatrix{
& \ob^C &
\ob^{C^A}\\
\ob^{B\xord A} & \twovf^f \ar[l]_{\tsource}
\ar[r]^{\lam} \ar[u]|{\ttarget} & \curry^{\lambda f}
\ar[u]|{\csource} \ar[d]|{\ctarget}
\\
&& \ob^B }\end{equation} }

Diagram~(\ref{lamcom}) below forces
Diagram~(\ref{evaleq}) to commute. {\wider
\begin{equation}\label{lamcom}
\xymatrix{ \ob^{B\xord A} & \arr^{\lambda f\xord
\Id[A]} \ar[l]_{\source} \ar[r]^{\target}
&\ob^{C^A\xord A}
\\
& \arr_2 \ar[u]|{\rfac} \ar[dl]_{\comp}
\ar[dr]^{\lfac}
\\
\arr^f \ar[r]_{\target} \ar[uu]|{\source} & \ob^C &
\arr^{\eval} \ar[l]^{\target} \ar[uu]|{\source} }
\end{equation}
}

Diagram~(\ref{otherlam}) below ensures that
requirement CCC--\ref{fourthreq} holds.
{\Narrower\taller
\begin{equation}\label{otherlam}
\xymatrix{ \ob^{C^A\xord A} & \arr^{\eval}
\ar[l]_{\source} \ar[r]^{\target} & \ob^C & \ob^B
\\
\arr^{g\xord\Id[A]} \ar[u]|{\target} \ar[d]|{\source}
& \arr_2^{<\eval,g\xord\Id[A]>} \ar[u]|{\lfac}
\ar[l]_{\rfac} \ar[dr]|{\comp} && \arr^g
\ar[u]|{\source} \ar[r]^{\target} & \ob^{C^A}
\\
\ob^{B\xord A} && \arr^{\eval\o(g\xord\Id[A])}
\ar[ll]_{\source} \ar[uu]|{\target} & \curry
\ar[u]|{\arrow}
\\
&&&\twovf \ar[-1,-3]|{\tsource} \ar[-3,-1]|{\ttarget}
\ar[ul]|{\arrow} \ar[u]|{\lam}
\\
}
\end{equation}
}

\subsection{Invertibility of $\lambda$}\label{laminvr}

It follows from CCC-\ref{fourthreq} that if $\Csc$ is
any Cartesian closed category corresponding to a model
$\C$, then $\C(\lambda)$ is a bijection.  The
Completeness Theorems~\ref{complthm}
and~\ref{compcomm} then imply that there is an arrow
$\lambda\inv$ in $\CCC$ that, as its name suggests, is
a formal inverse to $\lambda$.

\chapter{Graph-based logic}\label{gbl}
\section{Assertions in graph-based logic}\label{pfsec}

\defn
Let $\F$ be an $\Eb$-form.  A diagram of the kind
\begin{equation}\label{PFShow}
\xymatrix{
&D'\ar[d]^{f'}
\\
D''\ar[r]_{f''}
&D
\\
}
\end{equation}
in $\SynCat[\Eb,\F]$
is called a {\bf potential
factorization} or {\bf PF}.\edefn

\defn Suppose there are morphisms of diagrams
\begin{equation}\label{diagdiag3}
\xymatrix{
&\delta'\ar[d]^{\phi'}
\\
\delta''\ar[r]_{\phi''}
&\delta
\\
}
\end{equation}
for which
\alist
\item
$\delta$, $\delta'$ and
$\delta''$ are all diagrams in $\SynCat[\Eb,\F]$.
\item $D=\Lim[\delta]$, $D'=\Lim[\delta']$ and
$D''=\Lim[\delta'']$. \item $f'$ is the fill-in arrow induced by
$\phi'$ and $f''$ is the fill-in arrow induced by
$\phi''$.\ealist
Then Diagram~(\ref{diagdiag3}) is a {\bf description} of the
potential factorization~(\ref{PFShow}).\edefn

By Lemma~\ref{piessens}, the description of a
potential factorization can be taken to lie in the graph of the
constructor space sketch $\E$ that generates $\Eb$. There are in
general many descriptions of a PF.
The description is not part of the
structure.  We do not know of an example of a diagram of the
form of~(\ref{PFShow}) that we can prove does not have a
description.


\subsection{Notation}
We will use suggestive notation for a potential
factorization, as exhibited in the diagram below.
{\wider\wider
\begin{equation}\label{pfdiag1}
\xymatrix{
&\claim\ar[d]^{\claimcon}
\\
\hyp\ar[r]_{\hypcon}
&\wksp
\\
}
\end{equation}
}%
The names $\hypcon$, $\claimcon$ and $\wksp$
respectively abbreviate ``hypothesis construction'',
``claim construction'' and ``workspace''. The reason
for the names of the arrows and nodes is discussed
in~\ref{disc}. In many examples, including all those
in this monograph, $\claimcon$ is monic and
corresponds to a formal selection of a subset of those
objects formally denoted by $\wksp$.

\subsection{Actual factorizations}\label{afs}
If an actual factorization arrow $\verif:\hyp\to\claim$ can
be constructed using the rules of Chapter~\ref{thconst} in
$\SynCat[\Eb,\F]$ that makes
{\wider\wider
\begin{equation}\label{p17diag}
\xymatrix{
&\claim\ar[d]^{\claimcon}
\\
\hyp\ar[r]_{\hypcon}\ar[ur]^{\verif}
&\wksp
\\
}
\end{equation}
}%
commute, then we say that the potential factorization
(\ref{PFShow}) is {\bf deducible}. We remind the reader
that by~\ref{sec613}, $\SynCat[\Eb,\F]$ is
$\FinLimTh\left[\E,\Name[\F]\right]$, so that the rules of
Chapter~\ref{thconst} apply.

If for some model $\M$ of
$\SynCat[\Eb,\F]$ there is an arrow $\xi$ of $\Set$ that makes
{\wider\wider\wider\wider
\begin{equation}\label{newdiag}
\xymatrix{
&\M(claim)\ar[d]^{\M(claimcon)}
\\
\M(hyp)\ar[r]_{\M(hypcon)}\ar[ur]^{\xi}
&\M(wksp)
\\
}
\end{equation}}%
commute, then we say the model $\M$ {\bf satisfies} the
potential factorization. If for {\it every\/} model $\M$
there is such an arrow $\xi$ then we say that the potential
factorization is {\bf valid}.

We give examples of potential
factorizations in Sections~\ref{exam1} and~\ref{exam2}.
Section~\ref{disc} discusses the general concept of
potential factorization.

\section{Soundness and completeness}\label{compsub}

\thm\label{complthm}
In any syntactic category $\SynCat[\Eb,\F]$, a potential factorization is
deducible if and only if it is valid.
\ethm

\pf
That deducibility implies validity follows from the fact that functors
preserve factorizations.

For the converse, let
{\wider\wider
\begin{equation}\label{pfdiag1}
\xymatrix{
&\claim\ar[d]^{\claimcon}
\\
\hyp\ar[r]_{\hypcon}
&\wksp
\\
}
\end{equation}}
be a potential factorization in $\SynCat[\Eb,\F]$.  Suppose it is
valid.
Because the functor $\Hom(\hyp,-)$ is a model, the hypothesis
of the theorem implies that we may choose
an arrow $\xi$ of $\Set$ such that the diagram
{\wider\wider\wider\wider\wider\taller
\begin{equation}\label{homdi}
\xymatrix{
&\Hom(\hyp,\claim)\ar[d]^{\M(\hyp,\claimcon)}
\\
\Hom(\hyp,\hyp)\ar[r]_{\M(\hyp,\hypcon)}\ar[ur]^{\xi}
&\M(\hyp,\wksp)
\\
}
\end{equation}}%
commutes.
Define $\verif\ceq\xi(\id_{\hyp})$.  Then $\claimcon\oc\verif=\hypcon$,
so that the potential
factorization is deducible.
\epf

\rem Any model $\M$ is a functor that preserves finite limits, so
that if $\claimcon$ is monic, so is $\M(\claimcon)$.
In that case, $\M(\verif)$ necessarily equals $\xi$.

Section~\ref{laminvr} gives an example of how completeness can
be used.
\erem

\rems Proofs of soundness in string-based logic contain an induction
which is missing in the preceding argument. In the present
system, a theorem can be identified with an arrow of
$\SynCat[\Eb,\F]$. Chapter~\ref{thconst} describes the
recursive construction of arrows in the finite-limit theory
of a sketch and so is the analog of the inductive part of
string-based proofs of completeness. We repeat once more
that $\SynCat[\Eb,\F]$ is indeed the finite-limit theory of
a sketch, namely $\FinLimTh\bigl[\E[\Fa]\bigr]$
(see~\ref{sec613}) and hence the constructions in
Chapter~\ref{thconst} do apply in this case. \erems

\thp\label{compcomm}
Let $\delta:I\to\SynCat[\Eb,\F]$ be a diagram.  Suppose for every model
$\M$ of $\SynCat[\Eb,\F]$, $\M\oc\delta$ commutes.  Then $\delta$
commutes. \ethp

\pf Suppose
\begin{equation}\label{adiag}
\xymatrix{
A\ar[r]^{f}\ar[dr]_{h}
&B\ar[d]^{g}
\\
&C
\\
}
\end{equation}
is a diagram in $\SynCat[\Eb,\F]$.  Then, because $\Hom(A,-)$ is
a model,
{\wider\wider\wider\wider\taller
\begin{equation*}
\xymatrix{
\Hom(A,A)\ar[r]^{\Hom(A,f)}\ar[dr]_{\Hom(A,h)}
&\Hom(A,B)\ar[d]^{\Hom(A,g)}
\\
&\Hom(A,C)
\\
}
\end{equation*}}%
commutes.
By chasing $\Id[A]$ around the diagram both ways, we get
$g\oc f=h$,
so  Diagram~(\ref{adiag}) commutes as well.  The
general result follows because every diagram can
be triangulated.\epf

\section{Example: A fact about diagrams in any category}\label{exam1}

The following proposition holds in any category.

\thp\label{examthm} In any category, given the following diagram
\begin{equation}\label{sqdiag}
\xymatrix{
A \ar[r]^h \ar[d]_f & B
\ar[d]^k \\
C \ar[r]_g\ar[ur]^x & D
}
\end{equation}
if the two triangles commute then so does the outside square.
\ethp
\pf $k\oc h=k\oc(x\oc f)=(k\oc x)\oc f=g\oc f$.\epf

Let $\F$ be the $\Cat$-form ($\Cat=\FinLimTh[\Catsk]$)
with not necessarily commutative diagrams of the
form~(\ref{sqdiag}) as models.  Such a form can be realized by
specifying that $\Fa$ be a constant whose type is the limit of
Diagram~(\ref{basic}).

We construct here the potential factorization in
$\SynCat[\Cat,\F]$ that corresponds
to Proposition~\ref{examthm}. The construction takes place
entirely in $\FinLimTh[\Catsk]$; no reference to the constant $\Fa$
is made, since we are working directly with the description of
the type (codomain) of $\Fa$, which is Diagram~(\ref{basic}).

As we pointed out in Section~\ref{consnot},
an element of the value in a model of $\Lim[D(\ref{basic})]$ is
a diagram such as Diagram~(\ref{sqdiag}).  However,
Diagram~(\ref{basic})
carries only the information as to the source and target of the
arrows in Diagram~(\ref{sqdiag}).

The structure we must actually work with should include the
information as to which pairs are composable.  There are four
composable pairs in Diagram~(\ref{sqdiag}), and each one
inhabits the value of $\arr_2$, the node of formal
composable pairs in a category (see Appendix~\ref{realcatsk}).

The following diagram thus contains the basic information
about sources, targets and composability that are required for
stating  Proposition~\ref{examthm},
and so its limit is suitable for being the node
$\wksp$ of the proof corresponding to the Proposition.
{\wider\wider\taller
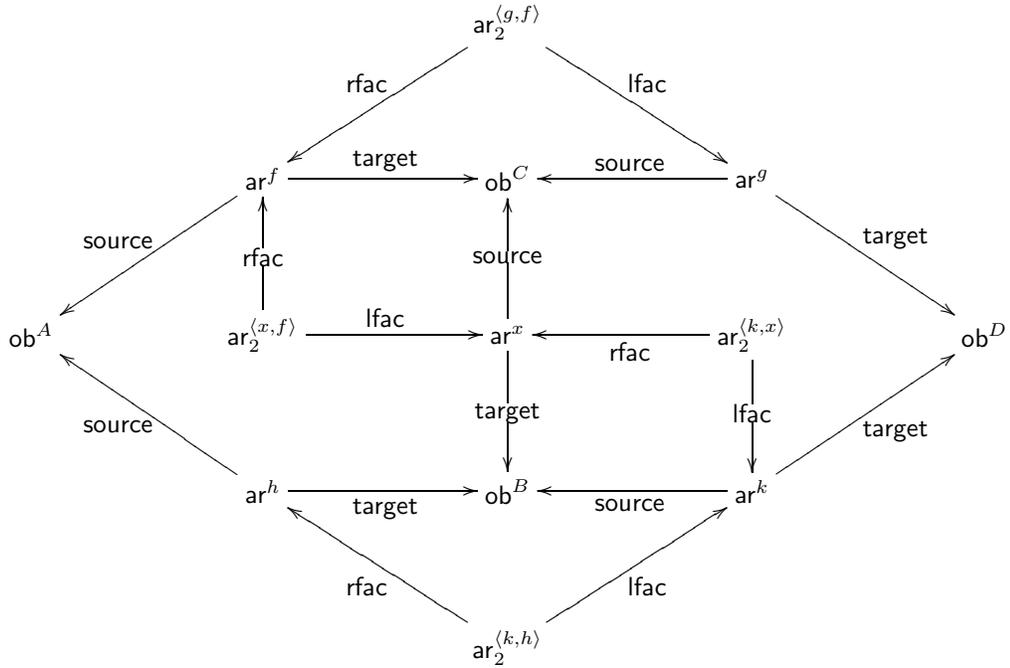
\begin{equation}\label{compinfo}
\xymatrix{
&&\arr_{2}^{<g,f>}\ar[dl]_{\rfac}\ar[dr]^{\lfac}
&&
\\
&\arr^{f}\ar[dl]_{\source}\ar[r]^{\target}
&\ob^{C}
&\arr^{g}\ar[l]_{\source}\ar[dr]^{\target}
\\
\ob^{A}
&\arr_{2}^{<x,f>}\ar[u]|{\rfac}\ar[r]^{\lfac}
&\arr^{x}\ar[u]|{\source}\ar[d]|{\target}
&\arr_{2}^{<k,x>}\ar[d]|{\lfac}\ar[l]^{\rfac}
&\ob^{D}
\\
&\arr^{h}\ar[lu]^{\source}\ar[r]_{\target}
&\ob^{B}
&\arr^{k}\ar[l]^{\source}\ar[ru]_{\target}
\\
&&\arr_{2}^{<k,h>}\ar[ul]^{\rfac}\ar[ur]_{\lfac}
&&
\\
}
\end{equation}
}%
The statement
that the two triangles in Diagram~(\ref{sqdiag}) commute is:
$x\oc f=h$ and $k\oc x=g$.  Using the composition arrow
$\comp:\arr_2\to \arr$ of $\SynCat[\Catsk]$, this
statement amounts to saying that
Diagram~(\ref{sqdiag}) is a member of
$\C\bigl(\Lim[D(\ref{hypoth})]\bigr)$ (below)
so we take $\hyp=\Lim[D(\ref{hypoth})]$.
{\wider\taller
\begin{equation}\label{hypoth}
\xymatrix{
&&\arr_{2}^{<g,f>}\ar[dl]_{\rfac}\ar[dr]^{\lfac}
&&
\\
&\arr^{f}\ar[dl]_{\source}\ar[r]^{\target}
&\ob^{C}
&\arr^{g}\ar[l]_{\source}\ar[dr]^{\target}
\\
\ob^{A}
&\arr_{2}^{<x,f>}\ar[u]|{\rfac}\ar[r]^{\lfac}
\ar[d]|{\comp}
&\arr^{x}\ar[u]|{\source}\ar[d]|{\target}
&\arr_{2}^{<k,x>}\ar[d]|{\lfac}\ar[l]^{\rfac}
\ar[u]|{\comp}
&\ob^{D}
\\
&\arr^{h}\ar[lu]^{\source}\ar[r]_{\target}
&\ob^{B}
&\arr^{k}\ar[l]^{\source}\ar[ru]_{\target}
\\
&&\arr_{2}^{<k,h>}\ar[ul]^{\rfac}\ar[ur]_{\lfac}
&&
\\
}
\end{equation}}%
The statement that the outside of Diagram~(\ref{sqdiag})
commutes is that the diagram is
a member of
$\C\bigl(\Lim[\delta(\ref{conc})]\bigr)$, so we
take $\claim=\Lim[D(\ref{conc})]$:
{\narrower\begin{equation}\label{conc} \xymatrix{
&&\arr_{2}^{<g,f>}\ar[dl]_{\rfac}\ar[dr]_{\lfac}
\ar[ddrrr]^{\comp} &&
\\
&\arr^{f}\ar[dl]_{\source}\ar[r]_{\target} &\ob^{C}
&\arr^{g}\ar[l]^{\source}\ar[dr]|{\target}
\\
\ob^{A}
&\arr_{2}^{<x,f>}\ar[u]|{\rfac}\ar[r]^{\lfac}
&\arr^{x}\ar[u]|{\source}\ar[d]|{\target}
&\arr_{2}^{<k,x>}\ar[d]|{\lfac}\ar[l]^{\rfac}
&\ob^{D}
&\arr^{k\oc g\oc f}
\\
&\arr^{h}\ar[lu]^{\source}\ar[r]^{\target} &\ob^{B}
&\arr^{k}\ar[l]_{\source}\ar[ru]|{\target}
\\
&&\arr_{2}^{<k,h>}\ar[ul]^{\rfac}\ar[ur]^{\lfac}
\ar[uurrr]_{\comp} &&
\\
}
\end{equation}}
Diagram~(\ref{compinfo})
is a restriction of both Diagram~(\ref{hypoth}) and
Diagram~(\ref{conc}). By Lemma~\ref{lemma1}, this inclusion induces
arrows
$$\claimcon:\Lim[D(\ref{conc})]\to\Lim[D(\ref{compinfo})]$$
and
$$\hypcon:\Lim[D(\ref{hypoth})]\to\Lim[D(\ref{compinfo})]$$
producing a potential factorization in $\FinLimTh[\Catsk]$
(hence in $\SynCat[\Cat,\F]$, which contains $\FinLimTh[\Catsk]$ as
a subcategory):
{\wider\wider\wider
\begin{equation}\label{pottoprove}
\xymatrix{
&\Lim[D(\ref{conc})]\ar[d]^{\claimcon}
\\
\Lim[D(\ref{hypoth})]\ar[r]_{\hypcon}
&\Lim[D(\ref{compinfo})]
\\
}
\end{equation}
}

This potential factorization expresses the content of
Proposition~\ref{examthm}
in diagrammatic form.  It should be clear that the node
$\Lim[D(\ref{compinfo})]$
could have been replaced by
$\Lim[D(\ref{basic})]$.

In Section~\ref{proofexamthm} we construct an arrow $\verif$ making the diagram
{\begin{equation}\label{pottoprovefulfilled}
\xymatrix{
&\Lim[D(\ref{conc})]\ar[d]^{\claimcon}
\\
\Lim[D(\ref{hypoth})]\ar[r]_{\hypcon}\ar[ur]^{\verif}
&\Lim[D(\ref{compinfo})]
\\
}
\end{equation}} an actual factorization (see Section \ref{afs}).

\section{Example: A fact about Cartesian closed
categories}\label{exam2} Proposition~\ref{examthm}
holds in any category. We now discuss a theorem of
Cartesian closed categories, to show how the system
presented in this monograph handles structure that
cannot be expressed using Ehresmann sketches.  The
latter are equivalent in expressive power to ordinary
first order logic. (An excellent presentation of the
details of this fact may be found in~\cite{arbook}.)
Thus this example in a certain sense requires
higher-order logic.

\thp\label{exam2thm}
In any Cartesian closed category, if
\begin{equation*}
\xymatrix{
A\x B\ar[r]^{g}\ar[dr]_{h}
&C\ar[d]^{f}
\\
&D
\\
}
\end{equation*}%
commutes, then so does
\begin{equation*}\label{mustcom}
\xymatrix{
A\ar[r]^{\lambda g}\ar[dr]_{\lambda h}
&C^{B}\ar[d]^{f^{B}}
\\
&D^{B}
\\
}
\end{equation*}%
\ethp
See Appendix~\ref{cccsec} for
notation.
The arrow $f^B$ is defined by
{\wider
\begin{equation}\label{fbdef}
\xymatrix{
f^{B}:=\lambda(C^{B}\x B\ar[r]^{\eval}
&C\ar[r]^{f}
&D):C^{B}\ar[r]
&D^{B}
\\
}
\end{equation}
}

The proof follows from the fact that $\lambda$ is invertible and
the calculation
$$\mld \eval\oc\left((f^B\oc\lambda g)\x \Id[B]\right) &=
\eval\oc\left((f^B\x \Id[B]) \oc (\lambda g\x\Id[B])\right)\\
=\left(\eval\oc(f^B\x\Id[B])\right)\oc(\lambda g\x \Id[B])\\
= (f\oc\eval)\oc(\lambda g\x\Id[B])\\
= f\oc(\eval\oc(\lambda g\x\Id[B]))\\
=f\oc g
= h= \eval\oc(\lambda h\x\Id[B])$$
The first equality is
based on an assertion true in all
categories that can be handled in our system in a manner
similar to (but more complicated than) that of
module~(\ref{MxNdiag}) in Appendix~\ref{modsubsec}.
The second and fourth equalities are
associativity of composition and are proven using
Figure~(\ref{u&ass})  of Appendix~\ref{realcatsk}.
The sixth equality is a hypothesis.  The other three equalities
are all based on Diagram~(\ref{lamcom}) in
Appendix~\ref{cccsec}.

We present here the potential factorization corresponding
to the third equality, which is the most complicated of
those based on Diagram~(\ref{lamcom}).  In this
presentation, unlike that of~\ref{exam1}, we will use the
modules developed in Chapter~\ref{conspsk} to simplify the
figures.  The actual factorization corresponding to this
potential factorization is given in~\ref{exam2sec}.

The fact under discussion is that the diagram
{\wider\wider\wider\wider\wider
\begin{equation}\label{mustcom}
\xymatrix{
A\x B\ar[r]^{\lambda g\x\Id[B]}
&C^{B}\x B\ar[r]^{f^{B}\x\Id[B]}\ar[d]_{\eval}
&D^{B}\x B\ar[d]^{\eval}
\\
&C\ar[r]_{f}
&D
\\
}
\end{equation}}%
commutes.

Thus $\wksp$ will be the limit of the following diagram,
which describes the objects and arrows in Diagram~(\ref{mustcom})
but has no requirements on its commutativity.
{\wider\wider
\begin{equation}\label{wk1}
\xymatrix{\arr^{\lambda g\x\Id[B]}
\ar[r]^{\target}\ar[d]_{\source}
&\ob^{C^{B}\x B}
&\arr^{f^{B} g\x\Id[B]}
\ar[r]^{\target}\ar[l]_{\source}
&\ob^{D^{B}\x B}
\\
\ob^{A\x B}
&\arr^{eval}
\ar[d]_{\target}\ar[u]^{\source}
&&\arr^{eval}
\ar[d]^{\target}\ar[u]_{\source}
\\
&\ob^{C}
&\arr^{f}
\ar[r]_{\target}\ar[l]^{\source}
&\ob^{D}
}
\end{equation}
}

We define $\hyp$ to be
the limit of the following diagram, in which
$\phi=(f\oc\eval)\oc(\lambda g\xo \Id[B])$.
{\taller\taller\narrower
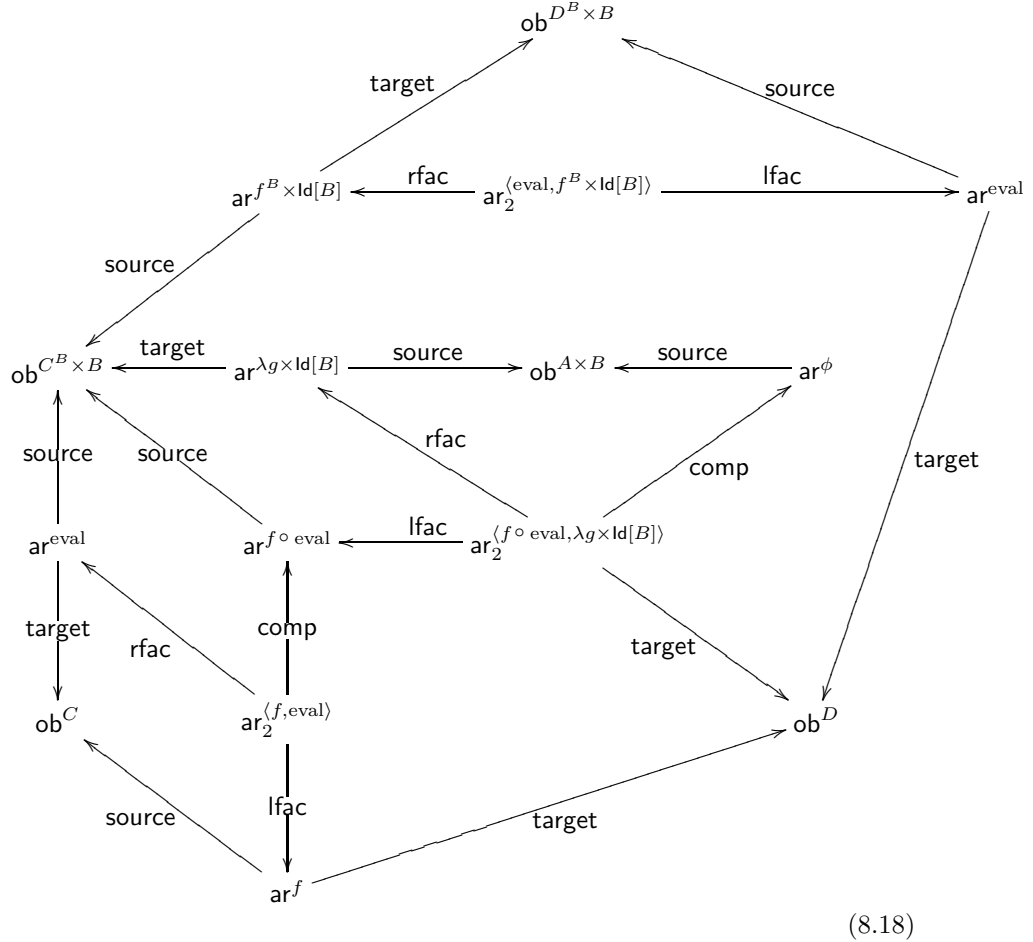
\begin{equation}\label{wk2}
\xymatrix{
&&\ob^{D^{B}\x B}
&&\\
&\arr^{f^{B}\x \Id[B]}
\ar[ur]^{\target}\ar[dl]_{\source}
&\arr_{2}^{<\eval,f^{B}\x \Id[B]>}
\ar[rr]^{\lfac}\ar[l]_{\rfac} &&\arr^{\eval}
\ar[llu]_{\source}\ar[dddl]^{\target}
\\
\ob^{C^{B}\x B}
&\arr^{\lambda g\x \Id[B]}
\ar[l]_{\target}\ar[r]^{\source}
&\ob^{A\x B}
&\arr^{\phi}
\ar[l]_{\source}
&\\
\arr^{\eval} \ar[d]|{\target}\ar[u]|{\source}
&\arr^{f\oc \eval} \ar[lu]|{\source} &\arr_{2}^{<f\oc
\eval,\lambda g\x \Id[B]>}
\ar[l]_{\lfac}\ar[lu]_{\rfac}\ar[ru]_{\comp}
\ar[dr]_{\target}
&&\\
\ob^{C} &\arr_{2}^{<f, \eval>}
\ar[lu]^{\rfac}\ar[u]|{\comp}\ar[d]|{\lfac} & &\ob^{D}
\\
&\arr^{f}
\ar[lu]^{\source}\ar[rru]_{\target}
&&&\\
}
\end{equation}}%
Then Diagram~(\ref{wk1}) is a subdiagram of
Diagram~(\ref{wk2}) (the big rectangle
in Diagram~(\ref{wk1}) is the perimeter of Diagram~(\ref{wk2}))
and we define $\hypcon$ to be
the induced arrow from $\Lim[D(\ref{wk2})]$ to
$\Lim[D(\ref{wk1})]$.

The object $\claim$ is the limit of
a diagram we shall refer to as Diagram~(\ref{wk2}$'$),
obtained from Diagram~(\ref{wk2}) by adjoining an arrow
labeled $\comp$ from
$\arr_2^{<\eval,f^B\x \Id[B]>}$ to
$\arr^{f\oc\eval}$. Diagram~(\ref{wk2}$'$) includes
Diagram~(\ref{wk2}) and hence Diagram~(\ref{wk1}), and we take
the arrow $\claimcon$ to be the arrow from
Diagram~(\ref{wk2}$'$) to Diagram~(\ref{wk1}) induced by this
inclusion.

We now have a potential factorization
{\wider\wider\wider
\begin{equation}\label{pfex2}
\xymatrix{
&\Lim[D(\ref{wk2}')]\ar[d]^{\claimcon}
\\
\Lim[D(\ref{wk2})]\ar[r]_{\hypcon}
&\Lim[D(\ref{wk1})]
\\
}
\end{equation}}

This potential factorization expresses the content of the
third equality in the calculation in the proof of
Proposition~\ref{exam2thm}. We provide an actual factorization

{\wider\wider\wider
\begin{equation}\label{pfex2}
\xymatrix{
&\Lim[\ref{wk2} ']\ar[d]^{\claimcon}
\\
\Lim[D(\ref{wk2})]\ar[r]_{\hypcon}\ar[ur]^{\verif}
&\Lim[D(\ref{wk1})]\\
}
\end{equation}}in Section \ref{exam2sec}.

\rem Diagram~(\ref{pfex2})
is a diagram in $\CCC$.  Let the $\CCC$-form $\F$
be determined by requiring that $\Fa$ be a freely adjoined
global element with target $\wksp$.  Then via the embedding
$\CatTh[\FinLim,\CCC]$ into $\SynCat[\CCC,\F]$,
Diagram~(\ref{pfex2}) is also a diagram in $\SynCat[\CCC,\F]$,
and if it has an actual factorization in $\CatTh[\FinLim,\CCC]$
then it also has an actual factorization in $\SynCat[\CCC,\F]$.\erem

\section{Discussion of the examples}\label{disc}

\subsection{General discussion}
In a potential
factorization
{\wider\wider
\begin{equation}\label{PF2}
\xymatrix{
&\claim\ar[d]^{\claimcon}
\\
\hyp\ar[r]_{\hypcon}
&\wksp
\\
}
\end{equation}
}%
each of $\hyp$, $\claim$ and $\wksp$ represents
a type of entity that can be
constructed in an $\Eb$-category, specifically
in an arbitrary category for
Example~\ref{exam1} and in any Cartesian closed category for
Example~\ref{exam2}. In this section, we discuss a way of
thinking about these nodes that exhibits how they could
represent a theorem about $\Eb$-categories.

\begin{enumerate}

\item The
node $\wksp$ (for ``workspace'') represents the data involved in
the union of the hypothesis and the conclusion.  For a given theorem,
the choice of what is actually included in $\wksp$ may be
somewhat arbitrary (see Item~\ref{exam2item1} in Section~\ref{discex2}
below).

\item The node $\hyp$
represents possible additional properties that are part of the
assumptions
in the theorem being represented.

\item The node $\claim$ represents the
properties that the theorem asserts must hold given the
assumptions.

\item The arrow $\hypcon$ represents the
selection or construction necessary to see the
hypothesis as part of
the workspace.  In both our examples, $\hypcon$ represents a
simple forgetting of properties.

\item The arrow $\claimcon$ represents the selection or
construction necessary to see the claim as part of the
workspace.

\item The arrow $\verif$ in an actual factorization
represents a specific way, {\it uniform in any model},
that any entity of type $\hyp$ can be transformed into, or
recognized as, an entity of type $\claim$.
\end{enumerate}

\subsection{Discussion of Example~\protect\ref{exam1}}
In Example~\ref{exam1},
\begin{enumerate}
\item $\wksp$ represent squares of the form
\begin{equation}\label{sqd2}
\xymatrix{
A \ar[r]^h \ar[d]_f & B
\ar[d]^k \\
C \ar[r]_g\ar[ur]^x & D
}
\end{equation}
with no commutativity conditions.
\item $\hyp$ represents diagrams of the form of
Diagram~(\ref{sqd2})
in which the two triangles commute;
\item $\claim$ represents diagram of the form of Diagram~(\ref{sqd2})
in which the outside square commutes.
\item $\hypcon$ represents forgetting that the two triangles commute.
Because it represents forgetting a property in this case,
$\hypcon$
is monic, but in general it need not be.
\item
$\claimcon$ represents forgetting that the outside square commutes.
\end{enumerate}

\subsection{Discussion of Example~\protect\ref{exam2}}
In Example~\ref{exam2}\label{discex2},
\begin{enumerate}
\item\label{exam2item1} $\wksp$ represent diagrams of the form of
Diagram~(\ref{mustcom}) with no commutativity conditions
and no recognition that any sequence of arrows is composable.
The phrase ``The form of Diagram~(\ref{mustcom})'' refers to the
source and target commonalities of the arrows in the diagram.
Obviously, some sequences compose but we have not represented
that in $\wksp$, although we could have.

\item $\hyp$ represents diagrams of the form
of Diagram~(\ref{mustcom}),
recognizing the composite $f\oc\eval$ and the fact that
$f\oc\eval,\,\lambda g\x\Id[B]$
and $\eval,\,f^B\x\Id[B]$ are composable pairs.

\item $\claim$ represents diagrams of the form of
Diagram~(\ref{mustcom}) that commute.

\item $\hypcon$ represents forgetting the information concerning
composition in $\hyp$.

\item
$\claimcon$ represents forgetting the information concerning
composition in $\claim$.

\end{enumerate}

We discuss the meaning of the actual factorization arrows
$\verif$ in~\ref{pfdis}.

\subsection{Explicit description instead of pattern recognition}\label{vfdisc} The representation of facts such as
those of Example~\ref{exam1} and~\ref{exam2} as  potential
factorizations is {\it variable free\/} in the sense that
in each statement, one does not refer to a particular
diagram such as Diagram~(\ref{sqdiag}) or
Diagram~(\ref{mustcom}) which stands as a pattern for all
such diagrams. Propositions~\ref{examthm}
and~\ref{exam2thm} state the fact in question using those
diagrams as  patterns, and understanding their meaning
calls on the reader's ability to recognize patterns.  Our
description of the fact in the examples as a potential
factorization is much more complicated because the diagrams
involved in the potential factorization are essentially
explicit descriptions of the relations between the nodes
and arrows of Diagram~(\ref{sqdiag}) and
Diagram~(\ref{mustcom}) respectively, relations which a
knowledgeable reader grasps from seeing the diagrams
without having them indicated explicitly.

Thus\mpark{Added paragraph} at the price of considerably more complexity we have substituted
{\it explicit description\/} of the structure {\it implied\/} by the diagrams in the assertions of
Propositions~\ref{examthm}
and~\ref{exam2thm}.  The relation between the explicit description and the written assertions reminds us of the relation between a program in a high-level programming language and the assembly code (or at least lower-level code) of a compiled form
of the assertions.  We believe this is an important step in the process of implementing on a computer the machinery described in this monograph.

\section{The rules of graph-based
logic}\label{ruled}\label{dedrules}

In first-order logic, rules in the\mpark{Added phrase} form of a context-free grammar are given for constructing terms and
formulas, and further rules (rules of inference) are given for
deriving formulas from formulas. These rules are intended to
preserve truth.

In Section~\ref{rulesapp} we gave rules for
constructing all the objects and arrows of $\FinLimTh[\S]$
(hence in any syntactic category -- see~\ref{sec613})
for an arbitrary finite-limit sketch $\S$, and for
constructing a basis for all the commutative diagrams
in $\FinLimTh[\S]$. These rules correspond to both the
term and formula construction rules and the rules of
inference of string-based logic.  The tools of a
typical string-based logic include constant symbols,
variables, function symbols, logical operators and
quantifiers. Here we have nodes, arrows and
commutative diagrams. (See Remark~\ref{whatarediagrams}.) What corresponds to a sentence
is  a potential factorization as in
Diagram~(\ref{PF2}), and what corresponds to the
satisfiability of the sentence in a model $\M$ is the
existence of an arrow $\xi$ in that model for which
{\wider\wider\wider\wider\wider\taller
\begin{equation}\label{newdiag2}
\xymatrix{ &\M(\claim)\ar[d]^{\M(\claimcon)}
\\
\M(\hyp)\ar[r]_{\M(\hypcon)}\ar[ur]^{\xi} &\M(\wksp)
\\
}
\end{equation}}%
commutes. A demonstration of the deducibility of the sentence corresponds to
the construction of an arrow $\verif:\hyp\to \claim$
in $\FinLimTh[\S]$ such that
{\wider\wider
\begin{equation}
\xymatrix{
&\claim\ar[d]^{\claimcon}
\\
\hyp\ar[r]_{\hypcon}\ar[ur]^{\verif}
&\wksp
\\
}
\end{equation}}%
commutes.  Thus the same rules suffice for {\it constructing\/}
the sentence (the potential factorization) and for {\it proving\/}
it (constructing the arrow that makes it an actual
factorization). See Section~\ref{disc} for further discussion of
these points.

\rems Each rule in Section~\ref{rulesapp} is actually
a {\it rule scheme\/}, and each instance of the scheme is an
assertion that, given certain arrows and commutative
diagrams (see Remark~\ref{remm} below) in
$\SynCat[\Eb,\S]$, other arrows or commutative
diagrams exist in $\SynCat[\Eb,\S]$. For example, if
$\delta$ is a diagram with shape graph
\begin{equation}
\xymatrix{
&i\ar[d]^{u}
\\
j\ar[r]_{v}
&k
\\
}
\end{equation}

then the following
rule is an instance of $\exists$LIM:
{\wider\wider\wider\wider\wider\wider\wider\taller\taller
\begin{equation}
\Rule
{
\xymatrix{
&\delta(i)\ar[d]^{\delta(u)}
\\
\delta(j)\ar[r]_{\delta(u)}
&\delta(k)
\\
}}
{\xymatrix{
\Lim[\delta]
\ar[r]^{\Proj[\Lim[\delta],i]}
\ar[d]_{\Proj[\Lim[\delta],j]}
\ar[dr]|{\Proj[\Lim[\delta],k]}
& \delta(i)
\ar[d]
\\
\delta(j)
\ar[r]
& \delta(k)
}}
{}
{}
\end{equation}}
\erems

\rem\label{remm} Rules $\exists$FIA, !FIA and CFIA assume
the existence of commutative cones, but a commutative cone
is a collection of interrelated commutative diagrams, so
the statement above that each scheme assumes the existence
of certain arrows and commutative diagrams is correct. Thus
given a cone $\Theta:v\coneto(\delta:I\to \Csc)$, an
instance of $\exists$FIA is this rule:
{\wider\wider\wider\wider\wider\wider
\wider\wider\wider\taller\taller
\begin{equation}
\Rule
{\xymatrix{
\Vertex[\Theta]
\ar[r]^{\Proj[\Lim[\Theta],i]}
\ar[d]_{\Proj[\Lim[\Theta],j]}
\ar[dr]|{\Proj[\Lim[\Theta],k]}
&\delta(i)\ar[d]
\\
\delta(j)\ar[r]
&\delta(k)
}}
{\xymatrix{
\Vertex[\Theta]
\ar[r]_\Fillin[\Theta'\delta]
&\Lim[\delta]
\\}}
{}
{}
\end{equation}}%

The point of this remark is that $\exists$FIA is a rule
with a diagram as hypothesis and an arrow as conclusion.
The hypothesis is the \emph{cone itself}, not the string
``$\Theta:v\coneto(\delta:I\to \Csc)$'' or any other
description of it.\erem

\rem\label{whatarediagrams}
In\mpark{Added Remark} the second paragraph of this section
we mentioned the tools of string-based logic -- strings made from symbols with an implicit structure given by a grammar.  The tools of graph-based logic are diagrams made from nodes and arrows.  We have been criticized for not giving a rigorous definition of concepts such as nodes, arrows and diagrams.  We would like to point out that logic texts generally do not spell out what symbols are (a decidedly subtle question) or what strings are.  The difference between the two situations is that string of symbols are more familiar to most logicians than diagrams are.  We claim that that is the \emph{only} difference between them.

\erem

\section{Example: Proof of Proposition~\protect\ref{examthm}}\label{proofexamthm}
We
continue Example~\ref{exam1} by constructing and thereby
deducing
the existence of an arrow
$\verif:\Lim[D(\ref{hypoth})]\to\Lim[D(\ref{conc})]$ making

\begin{equation}\label{pottoprovefulfilled}
\xymatrix{
&\Lim[D(\ref{conc})]\ar[d]^{\claimcon}
\\
\Lim[D(\ref{hypoth})]\ar[r]_{\hypcon}\ar[ur]^{\verif}
&\Lim[D(\ref{compinfo})]
\\
}
\end{equation}
commute.

We first construct
Diagram~(\ref{hypoth}$'$) (not shown)
by adjoining
$\arr_3$ to Diagram~(\ref{hypoth}), along with
arrows \begin{verse} $\lfac:\arr_3\to\arr^k$\\
$\mfac:\arr_3\to\arr^x$\\
$\rfac:\arr_3\to\arr^f$\\
\end{verse}
We further construct Diagram~(\ref{hypoth}$''$) by
adjoining
\begin{verse} $<\lfac,\mfac>:\arr_3\to\arr_2^{<k,x>}$\\
$<\mfac,\rfac>:\arr_3\to\arr_2^{<x,f>}$\\
$\lass:\arr_3\to\arr^{<g,f>}$\\
$\rass:\arr_3\to\arr_2^{<k,h>}$
\end{verse}
to Diagram~(\ref{hypoth}$'$).
These arrows are defined in
Appendix~\ref{realcatsk}.

Diagram~(\ref{compinfo}) is a base restriction of
each of Diagram~(\ref{hypoth}$'$)
and Diagram~(\ref{hypoth}$''$),
so, we may, using Lemma~\ref{lemma1},
choose
arrows
$$\phi_1:\Lim[D(\ref{hypoth}')]
\to\Lim[D(\ref{compinfo})]$$
and
$$\phi_2:\Lim[D(\ref{hypoth}'')]
\to\Lim[D(\ref{compinfo})]$$
Diagram~(\ref{hypoth}) is a dominant subdiagram of
Diagram~(\ref{hypoth}$'$) since the latter is obtained from
the former by adjoining a limit of a subdiagram
together with their projection arrows
($\lfac$, $\mfac$ and $\rfac$).
Therefore, using Lemma~\ref{lemma2}, we may choose an
isomorphism $\psi_1$ making
\begin{equation}\label{comover}
\xymatrix{
\Lim[D(38)] \ar[rr]^{\psi_{1}} \ar[dr]_{\hypcon}
&&
\Lim[38']\ar[dl]^{\phi_{1}} \\
&\Lim[D(37)] \\ }
\end{equation}
commute.

Similarly, Diagram~(\ref{hypoth}$'$) is a dominant subdiagram of
Diagram~(\ref{hypoth}$''$) since the latter is obtained from
the former by adjoining four fill-in arrows.
Therefore by Lemma~\ref{lemma2} we may choose an
isomorphism $\psi_2$ making
\begin{equation}\label{comover2}
\xymatrix{
\Lim[D(38')] \ar[rr]^{\psi_{1}} \ar[dr]_{\phi_{1}}
&&
\Lim[38'']\ar[dl]^{\phi_{2}} \\
&\Lim[D(37)] \\ }
\end{equation}
commute.
We then construct Diagram~(\ref{hypoth}$'''$)
by adjoining arrows
$$\comp:
\arr_2^{<g,f>}\to\arr^{k\oc x\oc f}$$
and
$$\comp:
\arr_2^{<k,h>}\to\arr^{k\oc x\oc f}$$
where $\arr^{k\oc x\oc f}$ is a new node.


Because of associativity (the right diagram
in Figure~(\ref{u&ass})  of Appendix~\ref{realcatsk}),
Diagram~(\ref{hypoth}$'''$) extends Diagram~(\ref{hypoth}$''$)
by adjoining a commutative cocone,
so we may choose an isomorphism
$\psi_3:\Lim[D(\ref{hypoth}'')]\to\Lim[D(\ref{hypoth}''')]$
and an arrow
$\phi_3:\Lim[D(\ref{hypoth}''')]\to\Lim[D(\ref{compinfo})]$
making
\begin{equation}\label{comover4}
\xymatrix{
\Lim[D(38'')] \ar[rr]^{\psi_{3}} \ar[dr]_{\phi_{2}}
&&
\Lim[38''']\ar[dl]^{\phi_{3}} \\
&\Lim[D(37)] \\ }
\end{equation}
commute.

Finally, by Lemma~\ref{lemma1}, we may choose arrow
$\psi_4:\Lim[D(\ref{hypoth}''')]\to\Lim[D(\ref{conc})]$ making
\begin{equation}\label{comover5}
\xymatrix{
\Lim[D(38''')] \ar[rr]^{\psi_{4}} \ar[dr]_{\phi_{3}}
&&
\Lim[39]\ar[dl]^{\claim} \\
&\Lim[D(37)] \\ }
\end{equation}
commute.

We next set
\begin{equation}\label{xidef}
\verif\ceq\psi_4\oc\psi_3\oc\psi_2\oc\psi_1
\end{equation}
whence the theorem follows.


\section{Example: Proof of Theorem~\protect\ref{exam2thm}}\label{exam2sec}
In this section, we provide a factorization of the potential
factorization described in Section~\ref{exam2}.

Diagram~(\ref{wk2}$'$) contains the following as a
subdiagram, in which, using Diagram~(\ref{fbdef}), $$\theta=
<\eval,f^B\x\Id[B]>=
<\eval,\lambda(f\oc\eval)\x\Id[B]>$$
{\wider\taller
\begin{equation}
\xymatrix{
\ob^{C^{B}\x B}
&\arr^{f^{B}\x \Id[B]}
\ar[r]^{\target}\ar[l]_{\source}
&\ob^{D^{B}\x B}
\\
&\arr_{2}^{\theta}
\ar[u]|{\rfac}
\ar[dl]|{\comp}
\ar[dr]|{\lfac}
&
\\
\arr^{f\oc \eval}
\ar[uu]^{\source}
\ar[r]_{\target}
&\ob^{D}
&\arr^{\eval}
\ar[uu]_{\source}
\ar[l]^{\target}
\\
}
\end{equation}}%
This diagram is an instance of Diagram~(\ref{lamcom}), so it
commutes.
The arrow $\comp$ satisfies Definition~\ref{addarrow}, so
Lemma~\ref{lemma2} implies that there is an isomorphism
$\verif:\hyp\to\claim$.  Now the inclusion of
Diagram~(\ref{wk1}) into Diagram~(\ref{wk2}) followed by the
inclusion of Diagram~(\ref{wk2}) into Diagram~(\ref{wk2}$'$) is
precisely the inclusion of
Diagram~(\ref{wk1}) into Diagram~(\ref{wk2}$'$). It follows that
$\hypcon\oc\verif\inv=\claimcon$, so that
$\claimcon\oc\verif=\hypcon$ as required.

\NoCompileMatrices

\section{Discussion of the proofs.}\label{pfdis}

The factorization $\verif:\Lim[D(\ref{hypoth})]\to\Lim[D(\ref{conc})]$
of Diagram~(\ref{pottoprove}) given in Equation~(\ref{xidef})
constitutes the recognition that if the two triangles commute,
then so does the outside square.
The fact that $\verif$  makes
Diagram~(\ref{pottoprove}) commute is a codification of the fact that
if the two triangles commute then so does the outside square
{\it of the same diagram}.  In general, the reason we require
that actual factorizations be an arrow in the comma
category
$(\SynCat[\Eb,\F]\downarrow\wksp)$ instead of merely an arrow
from one node to another is to allow us to assert hypotheses and
conclusions that share data (in this case the data in
Diagram~(\ref{sqd2})).

The factorization $\verif:\hyp\to \claim$ constructed
in~\ref{exam2sec}
constitutes recognition that $\eval\oc (f^B\x\Id[B])=f\oc\eval$ via
the arrow
$\comp:\arr_2^{<\eval,f^B\x\Id[B]>}\to\arr^{f\oc\eval}$ in
Diagram~(\ref{wk2}$'$).
Because of this, the node
$\arr_2^{<f\oc\eval,\lambda g\x\Id[B]>}$ could also
be labeled $<\eval\oc f^B\x \Id[B],\lambda
g\x\Id[B]>$.  Thus the factorization also exhibits the fact that
$$\eval\oc (f^B\x\Id[B])\oc(\lambda g\x\Id[B])=
f\oc\eval\oc(\lambda g\x\Id[B])$$

It is clear that there are many alternative formulations of
Proposition~\ref{exam2thm}.  For example, instead of
first constructing $\arr_2^{f\oc\eval}$ as in
Diagram~(\ref{wk2}) (which is $\hyp$ in this case),
we could have constructed a node $\arr_2^{\eval\oc
(f^B\x\Id[B])}$ and an arrow
$$\comp:\arr_2^{<\eval,(f^B\x\Id[B]>)}\to
\arr_2^{\eval\oc (f^B\x\Id[B])}$$
Then the construction of an arrow
$$\comp:\arr_2^{<f,\eval>}\to\arr^{\eval\oc (f^B\x\Id[B])}$$
would have proved the theorem.

\section{Discussion of graph-based logic}
This\mpar{I am not sure I reworded this the way you wanted.} formalism, which uses diagrams and mappings
between diagrams instead of strings of symbols, perhaps
seems unusual from the point of view of symbolic logic. It
contrasts with the usual string-based formalism in two
ways.  On the one hand, our formalism exhibits explicitly
much more detail than the string-based approach about the relationships between different
parts of the structure.  On the other, our formalism is
very close to the way it would be represented in a modern
object-oriented computer language as compared to
string-based formulas. The nodes become objects and the
arrows become methods.

Thus the arrow $\lam$ of CCC--~\ref{lamitem} (in
Section~\ref{cccsec}) can be directly represented in a
program object as the method that yields the exponential
adjoint of an arrow in a Cartesian closed category.  In
contrast, a formula in first order logic requires a rather
sophisticated parser to translate it into a data structure
on which a program can operate. Parsing is well-understood,
but it results in a computer representation (for example as
a tree or as reordered tokens on a stack) that is very
different from the formula before it is parsed. What must be represented in the computer is the formation-tree of a formula or term, not the string of symbols that
is usually thought of as the formula or term.

In this
sense, the appearance that string based logic
is simpler than graph-based logic is an illusion: For people to understand the structure of an expression represented as a string requires them to have
sophisticated pattern-recognition
abilities.  For a computer program to operate with such expressions requires an elaborate parser.

Thus the approach via diagrams has some of the advantages
(for example, transparent translation into a programming
object) and some of the disadvantages (for example, more of
the structure is explicit) of assembly language versus high
level languages. (See Remarks~\ref{diagrem}
a)
\chapter{Equational Theories}

\section{Signatures}

\subsection{Expressions and terms}\label{termexp}

In the description that follows of the terms and
equations for a signature, we use a notation that
specifies the variables of a term or equation
explicitly. In particular, one may specify variables
that do not appear in the expression. For this reason,
the formalism we introduce in the definitions below
distinguishes an \textit{expression\/} such as
$f(x,g(y,x),z)$ from  a \textit{term\/}, which is an
expression together with a specified set of typed
variables;  in this case that set could be for example
$\{x,y,z,w\}$.  This formalism is equivalent to that
of \cite{gogmes}.

\defn\label{sigf}

A pair $(\Sigma,\Omega)$ of sets together with two
functions $\Inp:\Omega\to\List[\Sigma]$ and
$\Outp:\Omega\to\Sigma$ is called a
\textbf{signature}. Given a signature
$\Ssc\ceq(\Sigma,\Omega)$, elements of $\Sigma$ are
called the \textbf{types} of $\Ssc$ and the elements
of $\Omega$ are called the \textbf{operations} of
$\Ssc$.\edefn

\notat\label{inptlist} Given a signature
$\Ssc=(\Sigma,\Omega)$, we will denote the set
$\Sigma$ of types by $\Types[\Ssc]$ and the set
$\Omega$ of operations by $\Oprns[\Ssc]$. For any
$f\in\Omega$, the list $\Inp[f]$ is called the
\textbf{input type list} of $f$ and the type
$\Outp[f]$ is the \textbf{output type} of $f$. \enotat

\rem\label{inptlrem} The input type list of $f$ is
usually called the \textbf{arity} of $f$, and the
output type of $f$ is usually called simply the
\textbf{type} of $f$.\erem
\defn\label{constdef}
An operation $f$ of a signature $\Ssc$ is called a
\textbf{constant} if and only if $\Inp[f]$ is the
empty list.\edefn

\defn\label{inhabdef} A type $\typevar$ of a signature $\Ssc$ is
said to be \textbf{inhabited} if and only if either
\alist \item there is a constant of output type
$\typevar$ in $\Ssc$, or \item there is an operation
$f$ of output type $\typevar$ for which every type in
$\Inp[f]$ is inhabited. \ealist The type $\typevar$ is
said to be \textbf{empty} if and only if it is not
inhabited. \edefn

\section{Terms and equations}

In this section, we define terms and equations in a given
signature.

\subsection{Assumptions}\label{assu}
In these definitions, we make the following
assumptions, useful for bookkeeping purposes.
\blist{A}\item\label{ordtyp} We assume that we are
given a signature $\Ssc$ for which
$\Types[\Ssc]=\sb{\sigma^i}{i\in I}$ for some ordinal
$I$. \item\label{doubleind} For each $i\in I$, we
assume there is an indexed set $\Vbl[\sigma^i]\ceq
\sb{x_j^i}{j\in\omega}$ whose elements are by
definition \textbf{variables of type $\sigma^i$}. In
this setting, $x_j^i$ is the $j$th variable of type
$\sigma^i$.

\item\label{orddd} The set of variables is ordered by
defining $$x_j^i\lt x_l^k \ceqv
\begin{cases}\text{either} & i\lt k\\ \text{ or }&i=k
\text{ and } j\lt l\end{cases}$$ \elist

We also define $\Vbl[\Ssc]\ceq
\union_{i\in\omega}\Vbl[\sigma^i]$.

\defn For any type
$\tau$, an \textbf{expression of type} $\tau$ is
defined recursively as follows. \blist{Expr} \item A
variable of type $\tau$ is an expression of type
$\tau$. \item If $f$ is an operation with
$\Inp[f]=(\gamma^i\mid i\in 1\twodots n)$ and
$\Outp[f]=\tau$, and $\lstsb{e}{1}{n}$ is a list of
expressions for which each $e_i$ is of type
$\gamma^i$, then $f\lstsb{e}{1}{n}$ is an expression
of type $\tau$.\elist\edefn

\defn\label{typenot}
The type of a variable $x$ is denoted by $\Type[x]$,
so that in the notation of~\ref{assu},
$\Type[x_j^i]=\sigma^i$. This notation will be
extended to include lists and sets of variables as
follows:
\begin{enumerate}
\item If $L=<x^1_2, x^1_2, x^1_3, x^2_2>$, then
$\Type[L]\ceq \sigma^1\xord\sigma^1\xord\sigma^1\xord\sigma^2$.

\item If $W=\{x^1_2,x^1_3,x^2_1,x^3_2\}$, then
$\Type[W]\ceq\sigma^1\xord\sigma^1\xord\sigma^2\xord\sigma^3$.
(Note that this depends on the ordering given by
A.\ref{ordtyp}.)

\end{enumerate}
The type of an expression $e$ will be denoted by
$\Type[e]$.\edefn

\rem Thus the function $\Type$ is overloaded: it may
be applied to variables, lists or sets of variables,
or expressions, and will in the following be applied
to terms and equations as well.

In every case, $\Type[z]$ denotes a single type, never
a list or set of types. In contrast, $\Types[z]$,
defined in Section~\ref{inptlist}\mpar{Fixed reference}, denotes a set of types  and
$\TypeList[z]$, defined below in
Definition~\ref{def234}, denotes a list of types.\erem

\defn\label{def234} For a given expression $e$,
$\VarList[e]$ is defined recursively by requiring that
\blist{VL} \item If $x$ is a variable of type
$\sigma$, $\VarList[x]\ceq(\sigma)$.

\item If $e=f(e_1,\ldots,e_n)$, then
$\VarList[e]\ceq(\VarList[e_1])\cdots(\VarList[e_n])$,
the concatenate of the lists $\VarList[e_i]$. \elist
\edefn \rem For a given expression $e$, $\VarList[e]$
is the list of variables in $e$, in order of
appearance in $e$ from left to right, counting
repetitions. \erem

\defn $\Rng\bigl[\VarList[e]\,\bigr]$, the set of
distinct variables occurring in $e$, is called the
\textbf{variable set} of $e$ and denoted by
$\VarSet[e]$. The list
$(\List[\Type])\bigl[\VarList[e]\,\bigr]$ is called
the \textbf{type list} of $e$, denoted by
$\TypeList[e]$. \edefn \rem If the $k$th entry of
$\VarList[e]$ is $x_j^i$, then the $k$th entry of
$\TypeList[e]$ is $\sigma^i$.\erem

\exam\label{expl} Let $e$ be the expression
$f(x,g(y,x),z)$. If $x$ and $y$ are variables of type
$\typevar$ and $z$ is of type $\tau$, then the
variable list of $e$ is $(x,y,x,z)$, the variable set
is $\{x,y,z\}$ and the type list is
$(\typevar,\typevar,\typevar,\tau)$. Using the
notation of A.\ref{doubleind} and supposing
$\typevar=\sigma^1$, $\tau=\sigma^2$, $x=x_1^1$,
$y=x^1_2$ and $z=x^2_1$, we have
$e=f(x_1^1,g(x_2^1,x_1^1),x_1^2)$ and the following
statements hold:\mpar{added $\Type[e]$}

\begin{align*} \VarList[e]&= (x_1^1,x_2^1,x_1^1,x_1^2)  \\
\VarSet[e]&=\{x_1^1,x_2^1, x_1^2\} \\
\TypeList[e]&= (\sigma^1, \sigma^1,\sigma^1,\sigma^2)
\\
\Type[e] &= \sigma^1\xord\sigma^1\xord\sigma^1\xord\sigma^2
\end{align*}
\eexam

\defn\label{termdeff}
A \textbf{term} $t$ for a signature $\Ssc$ is
determined by the following: \blist{TD}\item A set
$\Var[t]$ of typed variables.  (It is a set, not a
list, but it is ordered by the ordering of
A.\ref{orddd} in~\ref{assu}.) \item An expression
$\Expr[t]$. \item A type
$\Type[t]\in\Types[\Ssc]$.\elist These data must
satisfy the following requirements:\blist{TR} \item
$\VarSet\bigl[\Expr[t]\,\bigr]\includedin \Var[t]$.
\item $\Type[t] = \Type\bigl[\Expr[t]\,\bigr]$.\elist
\edefn \notat A given term $t$ will be represented as
the list $$(\Expr[t],\Var[t],\Type[t])$$ \enotat

\defn\label{InputTypes}
Let $t$ be a term. The list $\InputTypes[t]$ is
defined to be the list whose $i$th entry is the type
of the $i$th variable in $\Var[t]$ using the ordering
given by A.\ref{orddd} in~\ref{assu}. Thus if the
$k$th entry of $\Var[t]$ is $x_j^i$, then the $k$th
entry of $\InputTypes[t]$ is $\sigma^i$. Observe that
there are no repetitions in $\Var[t]$ but there may
well be repetitions in $\InputTypes[t]$. \edefn

\rem It follows immediately from
Definitions~\ref{typenot} and~\ref{InputTypes} that if
$t=(e,V,\tau)$ then
$$\prd\InputTypes[t]=\Type[V]$$\erem

\exam Let $e=f\left(x_1^1,g(x_2^1,x_1^1),x_1^2\right)$
as in Example~\ref{expl}, and suppose
$\Outp[f]=\sigma^5$. Then there are many terms $t$
with $\Expr[t]=e$, for example
$$t_1\ceq\left(e,\{x_1^1,x_2^1,x_1^2\},\sigma^5\right)$$ and
$$t_2\ceq\left(e,\{x_1^1,x_2^1,x_3^1,x_1^2,x_5^7\},\sigma^5\right)$$
We have $\Type[t_1] = \Type[t_2] = \sigma^5$ and (for
example)
$$\InputTypes[t_1]=(\sigma^1,\sigma^1,
\sigma^2)$$ and
$$\InputTypes[t_2]=(\sigma^1,\sigma^1,
\sigma^1,\sigma^2,\sigma^7)$$\eexam

\defn An \textbf{equation} $E$ is determined by
a set $\Var[E]$ of typed variables (ordered by our
convention) and two expressions $\Left[E]$,
$\Right[E]$, for which \blist{ER} \item\label{firster}
$\Type\bigl[\Left[E]\,\bigr]=\Type\bigl[\Right[E]\,\bigr]$.
\item $\VarSet\bigl[\Left[E]\,\bigr]
\union\VarSet\bigl[\Right[E]\,\bigr]\includedin
\Var[E]$.\elist \edefn \notat We will write $e=_Ve'$
to denote an equation $E$ with $V=\Var[E]$,
$e=\Left[E]$ and $e'=\Right[E]$. The notation
$\Type[E]$ will denote the common type of $\Left[E]$
and $\Right[E]$. \enotat

\exam\label{eqex2} Let $e$ be the expression
$f\left(x_1^1,g(x_2^1,x_1^1),x_1^2\right)$ of
Example~\ref{expl}.  Let $e'\ceq g(x_2^1,x_3^1)$.
Then there are many equations with $e$ and $e'$ as
left and right sides, for example:
\begin{equation}\label{mcee}
E_1 \ceq f\left(x_1^1,g(x_2^1,x_1^1),x_1^2\right)
=_{\{x_1^1,x^1_2,x_3^1,x_1^2,x_1^3\}} g(x_2^1,x_3^1)
\end{equation} and
\begin{equation}
E_2 \ceq f\left(x_1^1,g(x_2^1,x_1^1),x_1^2\right)
=_{\{x_1^1,x^1_2,x_3^1,x_1^2,x_2^5\}} g(x_2^1,x_3^1)
\end{equation}
\eexam

For later use, we need the following definition:

\defn\label{concdef} The \textbf{most concrete term} associated with an
expression $e$ is defined to be the unique term $t$
with the properties that $\Expr[t]=e$ and
$\Var[t]=\VarSet[e]$. The \textbf{most concrete
equation} associated with two expressions $e$ and $e'$
is defined to be the unique equation $E$ such that
$\Left[E]=e$, $\Right[E]=e'$,  and
$\Var[E]=\VarSet[e]\union\VarSet[e']$. \edefn

\exam We continue Example~\ref{eqex2}. The most
concrete equation associated with the expressions
$f\left(x_1^1,g(x_2^1,x_1^1),x_1^2\right)$ and
$g(x_2^1,x_3^1)$ is
\begin{equation}\label{referto}f\left(x_1^1,g(x_2^1,x_1^1),x_1^2\right)
=_{\{x_1^1,x_2^1,x_3^1,x_1^2\}}
g(x_2^1,x_3^1)\end{equation} The most concrete term
associated with
$f\left(x_1^1,g(x_2^1,x_1^1),x_1^2\right)$ is
$$\left(f\left(x_1^1,g(x_2^1,x_1^1),x_1^2\right),\{x_1^1,x_2^1,x_1^2\},
\sigma^5\right)$$ The most concrete term associated
with $g(x_2^1,x_3^1)$ is
$$\left(g(x_2^1,x_3^1),\{x_2^1,x_3^1\},\sigma^5\right)$$
assuming $\Outp[g]=\sigma^5$ (this \textit{must} be true if
Equation~(\ref{referto}) is true.) \eexam

\section{Equational theories}\mpark{We had failed to define
equational theory in the original paper.  In fact we
confused theory and signature.}

\defn\label{eqthdef}
An \textbf{equational theory} $(S,E)$ is a signature $S$
together with a set of equations $E$ in $S$.\mpar{I said we should
give an example of an equational theory and you agreed, but I think we can post this preliminary version without it.  When we do it, we ought to give the signature of the example in 9.1.} \edefn

This definition provides a concept of a
\textit{multisorted} equational theory.  Universal algebra
originated in the study of single-sorted equational
theories.

Our concern is with \textbf{multisorted equational logic}
(MSEL): a system of valid deduction for formulas in an
equational theory.

\section{Rules of inference of equational deduction}\label{rulinfms}
Goguen and Meseguer~\shortcite{gogmes} prove that the
following rules for equational deduction in multisorted
equational deduction are sound and complete.

\begin{description}
\item[reflexivity] $\Frac{}{e=_Ve}$. \item[symmetry]
$\Frac{e=_Ve'}{e'=_Ve}$. \item[transitivity]
$\Frac{e=_Ve'\quad e'=_Ve''}{e=_Ve''}$.

\item[concretion] Given a set $V$ of typed variables,
$x\in V$ and an equation $e=_Ve'$ such that $x\in
V\setminus(\VarSet[e]\union\VarSet[e'])$, and given
that $\Type[x]$ is inhabited,
$$\Frac{e=_V e'}{e=_{V\backslash \{x\}}e'}$$

\item[abstraction] Given a set $V$ of typed variables
and a variable $x$,
$$\Frac{e=_V
e'}{e=_{V\union\{x\}}e'}$$

\item[substitutivity] Given a set $V$ of typed
variables, $x\in V$, and expressions $u$ and $u'$ for
which $\Type[x]=\Type[u]=\Type[u']$,
$$\Frac{e=_Ve'\quad u=_Wu'}{e[x\leftarrow
u]=_{V\setminus\{x\}\union W}e'[x\leftarrow u']}$$

\end{description}

\section{Deductions in MSEL}\label{deducs}

We now define a deduction in MSEL of an equation $E$
from a list $(E_1,\ldots,E_n)$ of equations (called
\textbf{premises} in this context) as a rooted tree.
This definition is not as succinct as it could be, but
the form we give makes it easy to prove that every
deduction corresponds to an actual factorization
(Section~\ref{dedasfac2}).

\defn\label{deddef} Let $E$ be an equation and $P\ceq
(E_1,\ldots,E_n)$ a list of equations.  A
\textbf{deduction of} $E$ \textbf{from} $P$ has one of
the following four forms. \blist{D} \item $(E)$, where
$P=(E)$. \item $(E)$, where $P$ is the empty list and
$E$ is of the form $e=_Ve$.  (Note that reflexivity is
the only rule with empty premises.) \item $(E,D)$,
where $D$ is a deduction of an equation $E'$ from $P$
(the same list of premises) and $$\Frac{E'}{E}$$ is an
instance of a rule of inference of MSEL. \item
$(E,D_1,D_2)$, where for $i=1,2$, $D_i$ is a deduction
of an equation $E_i$ from a list of premises $P_i$,
$P=P_1P_2$ (the concatenate), and
$$\Frac{E_1\quad E_2}{E}$$ is an instance of a rule of inference
of MSEL. \elist\edefn

\chapter{Signatures to Sketches}\label{EqThChap}

We now show how to construct a finite-product sketch $\S$
corresponding to a given signature in such a way that the
categories of models of the signature and of the sketch are
naturally equivalent.\mpark{This presentation is not wrong
as it stands, but if we leave it as it is we need another
section defining the finite-product sketch of an equational
theory. Alternatively we could reword this whole chapter to
talk about equational theories instead of signatures, using
$\Tsc$ instead of $\Ssc$ and $\T$ instead of $\S$.}

\section{The sketch associated to a signature}\label{sksigg}
Given a signature $\Ssc=(\Sigma,\Omega)$, we now construct
a finite-product sketch $\Sk[\Ssc]$. This sketch, like any
finite-product sketch, determines and is determined (up to
isomorphism) by a finite-product form $\F$: Precisely (see
Chapter~\ref{skgeneral}), there is a diagram
$\delta:I\to\FinProd$ and a global element $\Fa:1\to
\Lim[\delta]$ in $\SynCat[\FinProd,\F]$ with the property
that the value of $\Fa$ in the initial model of
$\SynCat[\FinProd,\F]$ in $\Set$ consists (up to
isomorphism) of the graph, diagrams and (discrete) cones
that make up the sketch $\S$. Moreover, the finite-product
theory $\FPTh\bigl[\Sk[\Ssc]\bigr]$ (defined in
\cite{ctcs}, Section~7.5) is equivalent as a category to
the finite-product category
$\CatTh\bigl[\FinProd,\F]\bigr]$.

\section{The graphs and cones of $\Sk[\Ssc]$}\label{sksc}
In what follows, we recursively define arrows and
commutative diagrams in $\Sk[\Ssc]$ associated to
terms and equations of $\Ssc$ respectively.

\defn The set of nodes
of $\Sk[\Ssc]$ consist by definition of the following:
\blist{OS} \item Each type of $\Ssc$ is a node.
\item\label{listnode} Each list
$v=\lstsb{\typevar}{1}{n}$ that is the input type list
(see Remark~\ref{inptlist}) of at least one operation
in $\Omega$ is a node. \elist \edefn

\defn The arrows of
$\Sk[\Ssc]$ consist by definition of the
following:\blist{AS} \item Each operation $f$ in
$\Omega$ is an arrow $f: \Inp[f]\to\Outp[f]$. \item
For each list $v=\lstsb{\typevar}{1}{n}$ that is the
input type list of some operation in $\Omega$, there
is an arrow $\Proj[i]:v\to\typevar_i$ for each $i\in
1\twodots n$. (We will write $\Proj[v,i]$ for
$\Proj[i]$ if necessary to to avoid confusion, and on
the other hand we will write $p_i$ for $\Proj[i]$ in
some diagrams to save space.)\elist \edefn

\defn\label{conedef}
The cones of $\Sk[\Ssc]$ consist by definition of the
following: For each list
$v=(\typevar_1,\ldots,\typevar_n)$ that is the input
type list of some operation in $\Omega$, there is a
cone of $\Sk[\Ssc]$ with vertex $v$ and an arrow
$\Proj[i]:v\to\typevar_i$ for each $i\in 1\twodots n$.

It follows that in a model $M$ of the sketch
$\Sk[\Ssc]$, $M(v)=\prd_{i\in 1\twodots
n}M(\typevar_i)$.
\edefn

\subsection{Constants}
If the signature contains constants, then one of the
lists mentioned in OS.\ref{listnode} is the empty
list. As a consequence, the sketch will contain an
empty cone by Definition~\ref{conedef}, and the vertex
will become a terminator in a model.

\thl\label{constantlemma} If $\sigma$ is an inhabited
type of $\Ssc$, then there is a constant of type
$\sigma$ (global element of $\sigma$) in
$\FPTh\bigl[\Sk[\Ssc]\,\bigr]$.\ethl \pf The proof is
an easy structural induction.\epf

\section{Terms as arrows}
\label{termsandarr} We now describe how to associate
each term of a signature $\Ssc$ to an arrow in
$\FPTh\bigl[\Sk[\Ssc]\,\bigr]$ and each equation to a
commutative diagram or a pair of equal arrows in
$\FPTh\bigl[\Sk[\Ssc]\,\bigr]$.  The constructions
given here are an elaboration of those in \cite{ctcs},
pages~185--186.

\subsection{The arrow in
$\FPTh\bigl[\Sk[\Ssc]\,\bigr]$ corresponding to a
term}\label{conarr} We first define recursively two
arrows $\Sep[e]$ and $\Par[e]$ of
$\FPTh\bigl[\Sk[\Ssc]\,\bigr]$ for each expression
$e$, and an arrow $\Dia[t]$ of
$\FPTh\bigl[\Sk[\Ssc]\,\bigr]$ for each term $t$. The
arrow
$$\Arr[t]\ceq \Sep\bigl[\Expr[t]\,\bigr]\oc
\Par\bigl[\Expr[t]\,\bigr] \oc
\Dia[t]:\prd\InputTypes[t]\to\Outp[\Expr[t]]$$ will
then be the arrow of $\FPTh\bigl[\Sk[\Ssc]\,\bigr]$
associated with $t$; the meaning of the term $t$ in a
model of the signature is up to equivalence the same
function as the value of $\Arr[t]$ in the
corresponding model of $\Sk[\Ssc]$.

In these definitions, we suppress mention of the
universal model of $\Sk[\Ssc]$.  For example, if the
universal model is
$\UnivMod[\Ssc]:\Sk[\Ssc]\to\FPTh\bigl[\Sk[\Ssc]\,\bigr]$
and $\Theta$ is a node of $\Sk[\Ssc]$, then we write
$\Theta$ instead of $\UnivMod[\Theta]$.  We treat
arrows of $\Sk[\Ssc]$ similarly.

\defn\label{Qcon}
For an expression $e$, $\Sep[e]$ is defined
recursively by these requirements:\blist{Sep} \item If
$e$ is a variable $x$ of type $\tau$, then
$\Sep[e]\ceq\Id[\tau]$. Using the notation introduced
in Section~\ref{assu}, A.\ref{doubleind}, if
$e=x_j^i$, then $\Sep[e]\ceq\Id[\sigma^i]$. \item
Suppose $e=f\lstsb{e}{1}{n}$, where $f$ is an
operation with $\Inp[f]=\lstsb{\typevar}{1}{n}$ and
$\Outp[f]=\tau$. By definition, for $i\in(1\twodots
n)$, $\Type[e_i]=\typevar_i$.


Then $\Sep[e]$ is defined to be the arrow
{\wider\begin{equation}\label{diagtop40} \xymatrix{
\prd_{i=1}^n\Dom[\Sep\bigl[[e_i]\bigr]
\ar[rr]^{\prd_{i=1}^n \Sep[e_i]} &&
\crs{\typevar}{n}\ar[r]^{\,\,\,f} & \tau}
\end{equation}} \elist \edefn

\rems (a) Sep.2 recursively constructs the correct
parenthesization of the domain of $\Sep[e]$, as
illustrated in Examples~\ref{compex} and
Section~\ref{twoex}. (b) If $n=0$ in Sep.2, in other
words $\Inp[f]$ is empty, the composite in
Diagram~\eqref{diagtop40} becomes $$\xymatrix{ 1\ar[r]
& 1\ar[r]^f & \tau} $$ (c) ``$\Sep$'' is short for
``separated'', so named because $\Sep[e]$ represents
the expressions $e$ with variables renamed so that no
duplicates occur. It might have been better to refer
to this as the ``linearization'' of $e$, but we were
afraid this would cause confusion with linear
sketches.\erem

\defn
For an expression $e$, $\Par[e]$ is defined
recursively by these requirements:\blist{Par} \item If
$e$ is a constant of type $\tau$, then
$\Par[e]\ceq\Id[1]$, the formal identity of the formal
terminal object. \item If $e$ is a variable $x$ of
type $\tau$, then $\Par[e]\ceq\Id[\tau]$. \item If
$e=f\lstsb e1n$ as in Sep.2, then
{\begin{equation}\label{diagtop40} \xymatrixnocompile{
\prd\TypeList[e] \ar[r]^{\Ass[e]} &
\prd_{i=1}^n\left(\prd\TypeList[e_i]\right)
\ar[rr]^{\,\,\,\prd_{i=1}^n\Par[e_i]} &&
\prd_{i=1}^n\cod[\Par[e_i]]} \end{equation}}%
where $\Ass[e]$ is the canonical associativity arrow.
\elist \edefn

\rems (a) $\Par[e]$ is so called because it
parenthesizes $\prd_{i=1}^n\TypeList[e]$. (b) Note
that the definition ensures that $\Par[e]$ is an
identity arrow or a product of canonical associativity
arrows.\erems

\defn\label{dtdef} Let $t$ be an arbitrary term.
Then $$\Dia[t]:\prd\InputTypes[t]\to
\prd\TypeList\bigl[\Expr[t]\,\bigr]$$ is defined to be
the unique arrow induced by requiring that the
following diagrams commute for each pair $$(i,k)\in
\bigl(1\twodots\ab{\InputTypes[t]}\bigr)\x
\bigl(1\twodots\ab{\TypeList\left[\Expr[t]\,\right]}\bigr)$$
with the property that the $i$th variable from the
left in $\Expr[t]$ is $(\Var[t])_k$.

{\narrower\narrower\narrower\narrower
\begin{equation}\label{diagbottom40}
\xymatrix{ \prd\InputTypes[t] \ar[rr]^{\Dia[t]}
\ar[dr]_{\Proj[k]} & &
\prd\TypeList\bigl[\Expr[t]\,\bigr]
\ar[dl]^{\Proj[i]}\\
&(\InputTypes[t])_{k} &\\
}
\end{equation}
} \edefn

Alternatively, suppose $\Var[t]$ has cardinality $L$
and $\VarList\bigl[\Expr[t]\,\bigr]$ has length $M$.
Let $\phi:1\twodots M\to 1\twodots L$ be defined by
$\phi(m)=l$ if
$(\Var[t])_l=(\VarList\bigl[\Expr[t]\,\bigr])_m$
(there is a unique $l$ that makes this true). Then we
may also define
$$\Dia[t]:\prd\InputTypes[t]\to
\prd\TypeList\bigl[\Expr[t]\,\bigr]$$ to be the arrow
$(\Proj[\phi(m)]\mid m\in 1\twodots M)$.

This works because the $\left(\phi(m)\right)$th type
in $\prd\TypeList\bigl[\Expr[t]\,\bigr]$ is indeed the
type of the $\left(\phi(m)\right)$th variable in
$\VarSet\bigl[\Expr[t]\,\bigr]_m$ (see
Definition~\ref{def234}).

This is equivalent to requiring the diagrams
\eqref{diagbottom40} to commute.  The two definitions
are useful for different sorts of calculations and are
therefore included.

\rem $\Dia[t]$ is in some sense a generalized diagonal
map; hence the name.\erem

\exam Consider $e\ceq g(x^1_1,c)$, where $c$ is a
constant of type $\sigma^2$ and $g$ has type
$\sigma^3$.  Suppose
$$t=(g(x^1_1,c),\{x^1_1,x^4_1\},\sigma^3)$$
Then $t$ corresponds to the arrow {\taller\taller $$
\xymatrix{ & \sigma^1\x\sigma^4
\ar[d]_{\Proj[1]}^{\brspace\bigg\} \Dia[t]}\\ &
\sigma^1 \ar[d]_{<\Id[\sigma^1],!>}^{\brspace\bigg\}
\Par[e]}\\ & \sigma^1\x 1 \ar[d]_{\Id[\sigma^1]\xord
c} \ar@{}[dd]^{\brspace\biggg\} \Sep[e]} \\ &
\sigma^1\x\sigma^2 \ar[d]_g\\ & \sigma^3\\ } $$ Note
that one does not have to consider constants in
constructing $\Dia[t]$.} \eexam

\exam\label{compex} Let $e\ceq
f\left(x_1^1,g(x_2^1,x_1^1),x_1^2\right)$ with
$\Inp[g]=(\sigma^1,\sigma^1)$, $\Outp[g]=\sigma^2$.
$\Inp[f]=(\sigma^1,\sigma^2,\sigma^2)$, and
$\Outp[f]=\sigma^5$. Let
$$t\ceq(e,\{x_1^1,x_2^1,x_1^2,x_3^4\},\sigma^5)$$
Then
$$\VarList\bigl[\Expr[t]\,\bigr]=(x_1^1,x_2^1,x_1^1,x_1^2)$$
$$\InputTypes[t]=(\sigma^1,\sigma^1,\sigma^2,\sigma^4)$$
$$\Var[t]=\{x_1^1,x_2^1,x_1^2,x_3^4\}$$
$$\TypeList\bigl[\Expr[t]\,\bigr]=(\sigma^1,\sigma^1,\sigma^1,\sigma^2)$$
and
$$\Type[\Var[t]]=\prd\InputTypes[t]=\sigma^1\xord\sigma^1\xord\sigma^2\xord\sigma^4$$

If we use the first definition of $\Dia[t]$ in
Definition~\ref{dtdef}, then the following four
triangles must commute:

{ \spreaddiagramcolumns{-1.7em}
\mathchardef\times="0202

$$ \begin{array}{@{\hspace{-1.5em}}cc}
\xymatrix{
\sigma^1\xord\sigma^1\xord\sigma^2\xord\sigma^4
\ar[rr]^{\Dia[t]} \ar[dr]_{\Proj[1]} & &
\sigma^1\xord\sigma^1\xord\sigma^1\xord\sigma^2
\ar[dl]^{\Proj[1]}\\
&\sigma^1&\\
}
 &
\xymatrix{
\sigma^1\xord\sigma^1\xord\sigma^2\xord\sigma^4
\ar[rr]^{\Dia[t]} \ar[dr]_{\Proj[2]} & &
\sigma^1\xord\sigma^1\xord\sigma^1\xord\sigma^2
\ar[dl]^{\Proj[2]}\\
&\sigma^1&\\
}
\end{array}$$
\begin{equation}\label{page41}
 \begin{array}{@{\hspace{-1.5em}}cc}
\xymatrix{
\sigma^1\xord\sigma^1\xord\sigma^2\xord\sigma^4
\ar[rr]^{\Dia[t]} \ar[dr]_{\Proj[1]} & &
\sigma^1\xord\sigma^1\xord\sigma^1\xord\sigma^2
\ar[dl]^{\Proj[3]}\\
&\sigma^1&\\
} & \xymatrix{
\sigma^1\xord\sigma^1\xord\sigma^2\xord\sigma^4
\ar[rr]^{\Dia[t]} \ar[dr]_{\Proj[3]} & &
\sigma^1\xord\sigma^1\xord\sigma^1\xord\sigma^2
\ar[dl]^{\Proj[4]}\\
&\sigma^2&\\
}
 \end{array}
\end{equation}
} It follows that $\Dia[t]$ is given by the following
diagram, where to save space we write $p_k$ for
$\Proj[k]$.
\begin{equation}\label{pppp41} \xymatrix{
\sigma^1\xord\sigma^1\xord\sigma^2\xord\sigma^4
\ar[rrr]^{<p_1,p_2,p_1,p_3>} &&&
\sigma^1\xord\sigma^1\xord\sigma^1\xord\sigma^2\\
}
\end{equation}
and that $\Arr[t]$ is   the composite
$$\xymatrix{
\sigma^1\xord\sigma^1\xord\sigma^2\xord\sigma^4
\ar[d]_{<p_1,p_2,p_1,p_3>}^{\brspace\bigg\} \Dia[t]} \\
\sigma^1\xord\sigma^1\xord\sigma^1\xord\sigma^2
\ar[d]_{<p_1,<p_2,p_3>,p_4>}^{\brspace\bigg\} \Par[e]}
\\ \sigma^1\xord(\sigma^1\xord\sigma^1)\xord\sigma^2
\ar[d]_{\Id[\sigma^1]\xord g\xord\Id[\sigma^2]}
\ar@{}[dd]^{\brspace\biggg\} \Sep[e]}
\\
\sigma^1\xord\sigma^2\xord\sigma^2 \ar[d]_f
\\
\sigma^5 }$$\eexam \rem The first stage of recursion
in constructing the third node in the preceding
diagram gives
$$\Dom\bigl[\Sep[\Id[\sigma^1]\bigr]\x
\Dom\bigl[\Sep[g(x^1_2,x^1_1)]\bigr]\x\Dom\bigl[\Sep[\Id[\sigma^1]\bigr]$$\erem

\rem The definition of
$$\Arr[t]:\prd\InputTypes[t]\to\Outp[\Expr[t]]$$ determines a
unique arrow, once a choice of the product
$\prd\InputTypes[t]$ is made.  Changing the choice of
the products occurring in the intermediate stages of
the definition, namely
$\prd_{i=1}^n\Dom[\Sep\bigl[[e_i]\bigr]$ (in the case
where recursion is required) and
$\prd_{i=1}^n\TypeList[e]$ do not change $\Arr[t]$
because $\Par[e]$ is the unique associativity
isomorphism, which commutes with the isomorphisms
between different choices of products. \erem


\section{The finite-product sketch associated with an equational theory}\label{diageqass}
\mpar{changed the name of this section}
Let the equation $E\ceq e=_Ve'$ be given. Define the
terms $t_1$ and $t_2$ by $t_1=(e,\Var[E],\Type[E])$
and $t_2=(e',\Var[E],\Type[E])$ (using the notation
of~\ref{termdeff}). Recall that
$\Type[E]=\Type[e]=\Type[e']$. The notation
$\InputTypes[E]$ will denote the list
$\InputTypes[t_1]$, which is the same as
$\InputTypes[t_2]$. As in~\ref{conarr}, we have arrows
$\Arr[t_1]$ and $\Arr[t_2]$ with the same domain and
codomain. We will associate the diagram
\begin{equation}\label{p42first}
\xymatrix{
\InputTypes[E]\ar[rr]<1ex>^{\quad\Arr[t_1]}
\ar[rr]<-1ex>_{\quad\Arr[t_2]} &&
\Type[E]\\
}
\end{equation}
to the equation $E$.  By~\ref{conarr}, this is the
same as {\taller
\begin{equation}\label{p42second}
\xymatrix{
\InputTypes[E]\ar[rr]^{\quad \Dia[t_1]} \ar[d]|{\Dia[t_2]} &&
\TypeList[e]\ar[d]|{\Sep[t_1]\oc
\Par[t_1]}\\
\TypeList[e'] \ar[rr]_{\quad \Sep[t_2]\oc \Par[t_2]} && \Type[E] \\
}
\end{equation}
}

This completes the translation.

\rem The\mpar{Removed reference to 12.3 which I don't understand.} commutative diagram as exhibited above can
also be viewed as a pair of formally equal arrows as
in Diagram~\eqref{p42first}, and in what follows we
will use this description frequently.\erem


\defn
The finite-product sketch associated with an equational theory consists by definition of the sketch associated with the signature of the theory (defined in Section~\ref{sksigg}), to which is adjoined the diagram just defined that is associated to each equation of the theory.
\edefn
\todo{At this point we need to show how to construct the equational theory associated with each finite-product sketch.  See the note at the end of Section~\ref{formssec}.}

\chapter{Substitution}\label{subssec}

As terms are defined recursively, substitution may be
defined either by structural recursion or using
composition. These two ways of defining substitution are
convenient for different purposes.  Here we establish the
equivalence of the two procedures.  It may be useful to the
reader to compare the following constructions to the
examples in Section~\ref{twoex}.

\section{Recursive definition} We define substitution in
expressions, then in terms.

Given an expression $e$, the result of substituting an
expression $u$ for the variable $x$ in $e$ is denoted
by $e[x\leftarrow u]$ and is defined in this way:
\blist{S} \item If $c$ is a constant, $c[x\leftarrow
u] \ceq c$.

\item If $x$ is a variable, $x[x\leftarrow u]\ceq u$.

\item If $x$ and $y$ are different variables,
$x[y\leftarrow u]\ceq x$.

\item If $e=f(e_1,\ldots,e_n)$, then $$e[x\leftarrow
u] \ceq f\left(e_1[x\leftarrow
u],\ldots,e_n[x\leftarrow u]\right)$$

\elist

Now let $t\ceq \left(e,V,\sigma \right)$ be a term.
The expression $t\left[x\leftarrow \left(u,W,\tau
\right)\right]$ denotes the result of substituting the
term $ \left(u,W,\tau \right) $ for $x$ in $t$.  This
expression is defined only if $\Type[x]=\tau$, and it
is defined in this way:

$$t\left[x\leftarrow
\left(u,W,\tau \right)\right]\ceq \bigl(
e\left[x\leftarrow\left(u,W,\tau \right) \right],
\left(V\setminus\{x\} \right)\union W, \tau \bigr)$$

In particular,
$$\Arr\bigl[t\left[x\leftarrow
\left(u,W,\tau \right)\right]\bigr]\ceq
\Arr\bigl[\bigl( e\left[x\leftarrow\left(u,W,\tau
\right) \right], \left(V\setminus\{x\} \right)\union
W, \tau \bigr)\bigr]$$

\section{Direct definition}\label{dirdef} The alternative
way suggested by Examples~\ref{fexam} and~\ref{sexam} is to
define $\Arr\bigl[t\left[x\leftarrow \left(u,W,\tau
\right)\right]\bigr]$ directly, given
\begin{gather} \Arr \left[e,V,\sigma
\right]=\Sep[e]\Par[e]\Dia[(e,V,\sigma)]\\ \Arr
\left[u,W,\tau \right]=\Sep[u]\Par[u]\Dia[(u,W,\tau)]
\end{gather}

Note that $\Arr\left[e,V\union W,\sigma\right]=
\Sep[e]\Par[e]\Dia[(e,V\union W,\sigma)]$, so that
$\Arr\left[e,V\union W,\sigma\right]$ differs from
$\Arr\left[e,V,\sigma\right]$ only in the
$D$-composand of the arrow.

We have
$$\Dom\bigl[D \left[(e,V\union W,\sigma)  \right]
\bigr] = \Type \left[V\union W \right]
=\prd\InputTypes\left[e,V\union W,\sigma\right]$$

Now we define an arrow
$$\Insert\left[u,x,V,W\right]:
\Type\bigl[\left(V\setminus\{x\} \right)\union
W\bigr]\to \Type \left[V\union W \right]$$ It is
defined whenever $V$ and $W$ are sets of variables,
$x$ is a variable, and $u$ is an expression with
$\Type[u]=\Type[x]$. If $x\notin V$ we take
$\Insert\left[u,x,V,W\right]$ to be the identity
arrow. Otherwise,  choose $I\in
1\twodots\Length[V\union W]$ such that $(V\union
W)_I=x$. Then, for all $i\in \bigl(1\twodots\Length
\left(V\union W \right)\bigr)\setminus\{I\}$,
$$\Insert\left[u,x,V,W\right]_i\ceq \begin{cases}\Arr \left[u,
\left(V\setminus\{x\} \right) \union W, \tau \right] &
\text{if $i=I$}\\ \Proj[i-1] & \text {if $i\gt I$}\\
\Proj[i] & \text{otherwise}

\end{cases}$$

Note that $\Type \bigl[ \left(V\setminus\{x\}
\right)\union W\bigr]$ and $\Type \left[V\union
W\right]$ can differ in at most one factor depending
on whether $x\in W$ or not.

Finally, we define $$\Arr\bigl[t\left[x\leftarrow
\left(u,W,\tau \right)\right]\bigr]\ceq \Arr \left[e,
V\union W, \sigma\right]\oc
\Insert\left[u,x,V,W\right]$$ We have given two
methods of obtaining the arrow corresponding to the
term for which substitution has been made.  It remains
to be seen that these two methods give the same arrow.

\section{Proof of the equivalence of the two constructions}

The proof will be by structural induction.

\paragraph*{S.1} If $f$ is an arrow that factors
through a terminal object, then $f\oc g$ also factors
through the terminal object for any $g$, in particular
for $\Insert$.

\paragraph*{S.2} $t= \left(x,V,\sigma \right)$ and
$\sigma=\tau$. We note that
\begin{align*}
\Arr[x,V\union W,\sigma] &=
\Sep[x]\oc \Par[x]\oc \Dia[x,V\union W,\sigma]\\
&= \Id[\sigma]\oc\Id[\sigma]\oc
\proj[I]\quad\text{(where $ \left(V\union W
\right)_I=x$)}\\ &=\proj[I] \end{align*}

From the direct definition we have
\begin{multline*}
\Arr\bigl[t[x\leftarrow (u,W,\tau)]\bigr] = \Arr
\left[x, V\union W,\tau  \right]\oc
\Insert\left[u,x,V,W\right]\\
=\proj[I]\oc
\Insert\left[u,x,V,W\right]\\
=\Arr \bigl[u, \left(V\setminus\{x\} \right)\union
W,\tau\bigr]
\end{multline*}

by the definition of $\Insert\left[u,x,V,W\right]$,
which agrees with the recursive definition.

\paragraph*{S.3}
If $x\notin V$, $\Insert[u,x,V,W]$ is the identity
arrow by definition. If $x\in V$, then as in the proof
for S.2 above, $\Arr[x,V\union W,\sigma]=\Proj[I]$,
and since $y\neq x$, $y$ has a different index in
$V\union W$, so that
$\Proj[I]\oc\Insert[u,y,V,W]=\Proj[I]$.

\paragraph*{S.4}
In this case, we assume $t= \left(f(e_1,\ldots,e_n),
V,\sigma \right)$, where $\Outp[f]=\sigma$ and for all
$i\in 1\twodots n$, $\Outp[e_i]=\gamma_i$.

To begin with, suppose we have a term $t'$ that is
just like $t$ except for the set of variables, so that
$$t'\ceq \left(f(e_1,\ldots,e_n) ,V',\sigma\right)$$
where $\VarSet \left[f(e_1,\ldots,e_n)
\right]\includedin V'$. For each $i=1,\ldots,n$, we
require that the following diagram commute:
\begin{equation}
\xymatrix{ \prd\InputTypes[t'] \ar[d]_{\Dia[t']}
\ar[ddr]^{\phi}
\\
\prd\TypeList[t'] \ar[d]_{\Ass[e]}
\\
\prd_{i=1}^n\TypeList[e_i] \ar[dd]_\alpha
\ar[r]_-{\Proj[i]} & \TypeList[e_i]
\ar[d]^{\Par[e_i]}\\
& \Dom[\Sep[e_i] \ar[d]^{\Sep[e_i]}\\
\prd_{i=1}^n\gamma_i \ar[r]_{\Proj[i]} \ar[d]_f &
\gamma_i
\\
\tau }\end{equation}

Some observations about these diagrams: \alist \item
$\prd\TypeList[t']=\prd\TypeList[t]=\prd\TypeList[f(e_1,\ldots,e_n)]$.
\item $\alpha=\prd_{i=1}^n\Sep[e_i]\oc\Par[e_i]=
\bigl<\Sep[e_1]\oc\Par[e_1]\oc\Proj[e_1],\ldots,\Sep[e_n]\oc\Par[e_n]\oc\Proj[e_n]\bigr>$.
\item $\phi=\Dia[e_i,V',\sigma]$. \ealist

It follows that
\begin{align}\label{cruxeq}
\Arr[t']\nonumber &=f\oc
\bigl<\Sep[e_1]\oc\Par[e_1]\oc\Proj[1],\ldots,\Sep[e_n]\oc\Par[e_n]\oc\Proj[e_n]\bigr>
\oc\Ass[e]\oc\Dia[t']\\ \nonumber &= f\oc
\bigl<\Sep[e_1]\oc\Par[e_1]\oc\Dia[e_1,V',\gamma_1],
\ldots,\Sep[e_n]\oc\Par[e_n]\oc\Dia[e_n,V',\gamma_n]\bigr>
\\
&=
f\oc\bigl<\Arr[e_1,V',\gamma_1],\ldots,\Arr[e_n,V',\gamma_n]\bigr>
\end{align}

Now we return to our assumption that $t=
\left(f(e_1,\ldots,e_n), V,\sigma \right)$. By
induction hypothesis, we have, for all $i\in 1\twodots
n$,
\begin{equation*}\Arr \bigl[
e_i[x\leftarrow u], \left(V\setminus\{x\}
\right)\union W,\gamma_i \bigr] =\Arr \bigl[
e_i,V\union W,\gamma_i\bigr]\oc
\Insert\left[u,x,V,W\right]
\end{equation*}

By the direct definition,
$\Arr\bigl[t\left[x\leftarrow
\left(u,W,\tau\right)\right]\bigr]$ is
\begin{multline*} \Arr \left[f(e_1\ldots,e_n), V\union
W,\sigma  \right]\oc \Insert\left[u,x,V,W\right]\\ =
f\oc \bigl< \Arr \left[e_1,V\union W,\gamma_1
\right],\ldots, \Arr \left[e_n,V\union W,\gamma_n
\right]
\bigr>\oc \Insert\left[u,x,V,W\right]\\
\shoveright{\text{(by (\ref{cruxeq}))}}\\
=f\oc \bigl< \Arr \left[e_1,V\union W,\gamma_1
\right]\oc
\Insert\left[u,x,V,W\right],\\
\quad\quad\ldots, \Arr \left[e_n,V\union W,\gamma_n
\right]\oc
\Insert\left[u,x,V,W\right]\bigr>\\
= f\oc \bigl< \Arr \left[e_1[e\leftarrow u],
\left(V\setminus\{x\}\right)\union W, \gamma_1 \right]
, \ldots, \Arr \left[e_n[e\leftarrow u],
\left(V\setminus\{x\}\right)\union W, \gamma_n \right]
\bigr>\\
\shoveright{\text{(by induction hypothesis)}}\\
=\Arr\bigl[f \left(e_1[x_1\leftarrow u] ,\ldots,
e_n[x_n\leftarrow u] \right),
\left(V\setminus\{x\} \right)\union W,\sigma\bigr]\\
\shoveright{\text{(by (\ref{cruxeq}))}}\\
\end{multline*}
which is $\Arr\bigl[t\left[x\leftarrow
\left(u,W,\tau\right)\right]\bigr]$ by the recursive
definition. This completes the proof of the
equivalence of the two definitions.

Later, we shall use the equivalence of these two
methods of obtaining the arrow corresponding to the
term in which substitution has been made. To
facilitate reference we record this in the form of a
lemma. \thl\label{facil} Let $t\ceq \left(e,V,\sigma
\right)$ and $t'\ceq \left(u,W,\tau\right)$ be terms,
suppose $x\in V$ and suppose $\Outp[u]=\Type[x]$ so
that $u$ may be substituted for $x$. Then $$\Arr\bigl[
t \left[x\leftarrow
u\right],\left(V\setminus\{x\}\right)\union W,\tau
\bigr] = \Arr\bigl[ \left(e,V\union W,\sigma
\right)\bigr]\o \Insert[u,x,V,W]$$ \ethl

\section{Two extended examples of the
constructions}\label{twoex}

\exam\label{fexam} $$ \begin{array}{l} e \ceq f\left( \ds
x11, \ds x43, \ds x32, \ds x11, g( \ds x11, \ds x32), \ds
x21\right)\\ \\
\begin{cases}\Inp[f] = \sigma^1\x \sigma^4\x \sigma^3\x \sigma^1\x \sigma^5\x
\sigma^2 & \\ \Outp[f] = \sigma^5\end{cases}\\ \\
\begin{cases}\Inp[g] =
\sigma^1\x \sigma^3&\\ \Outp[g]= \sigma^5&\end{cases}\\
\\ V=\{\ds x11,\ds
x12,\dsu x13,\ds x21,\dsu x22,\dsu x31,\ds x32,\ds x43\}
\end{array}$$ The underlined variables are redundant; that
is, they do not appear in the expression $e$.
$$ \begin{array}{l}
\begin{cases}u \ceq h( \ds x21, \ds
x32)&\\ \Inp[h] = \sigma^2\x \sigma^3&\\ \Outp[h] = \sigma^3&\end{cases}\\ \\
W \ceq\{ \dsu x11, \ds x21, \dsu x22, \dsu x23, \dsu x31,
\ds x32, \dsu x33\} \end{array}$$ $ \ds x32$ is a variable
for which we are making a substitution. We wish to
calculate
\begin{align*}(e,V, \sigma^5)\left[ \ds x32 \leftarrow
(u,W, \sigma^3)\right] &= \bigl(e[ \ds x32\leftarrow
u],(V\setminus\{ \ds x32\})\union W,
\sigma^5\bigr)\end{align*} By direct calculation,
$$(V\setminus\{ \ds x32\})\union W=
\{ \ds x11, \ds x12, \dsu x13, \ds x21, \dsu  x22, \dsu
x23, \dsu x31, \ds x32, \dsu x33, \ds x43\}$$ and
$$e( \ds x32\leftarrow u)=
f\left( \ds x11, \ds x43, e( \ds x21, \ds  x32), \ds x11,
g( \ds x12, e( \ds x21, \ds x32)), \ds x21\right)$$ We now
exhibit the arrows for $e$ and $u$ over $V$: $e = f\left(
\ds x11, \ds x43, \ds x32, \ds x11, g( \ds x11, \ds x32),
\ds x21\right)$: {\def\brspace{\hspace{5em}}
$$\xymatrix{
\sigma^1\xord\sigma^1\xord\mk{\sigma}^1\xord\sigma^2
\xord\mk{\sigma}^2\xord\mk{\sigma}^3\xord\sigma^3
\xord\sigma^4
\ar[d]_{<p_1,p_8,p_7,p_1,p_2,p_7,p_4>}^{\brspace\bigg\}
\Dia[e]}
\\
\sigma^1\xord\sigma^4\xord\sigma^3\xord\sigma^1
\xord\sigma^1\xord\sigma^3 \xord\sigma^2
\ar[d]_{<p_1,p_2,p_3,p_4,<p_5,p_6>,p_7>}^{\brspace\bigg\}
\Par[e]}
\\
\sigma^1\xord\sigma^4\xord\sigma^3\xord\sigma^1\xord(\sigma^1\xord\sigma^3)
\xord\sigma^2 \ar[d]_{\Id[\sigma^1]\xord
\Id[\sigma^4]\xord\Id[\sigma^3]\xord \Id[\sigma^1]\xord
g\xord\Id[\sigma^2]} \ar@{}[dd]^{\brspace\biggg\} \Sep[e]}
\\
\sigma^1\xord\sigma^4\xord\sigma^3\xord\sigma^1\xord\sigma^5
\xord\sigma^2 \ar[d]_f
\\
\sigma^5 }$$ $u=h(\ds x21, \ds x32)$ (over $W$):}
$$\xymatrix{
\mk{\sigma}^1\xord\sigma^2\xord\mk{\sigma}^2
\xord\mk{\sigma}^2\xord\mk{\sigma}^3\xord\sigma^3
\xord\mk{\sigma}^3
\ar[d]_{<p_2,p_6>}^{\brspace\bigg\} \Dia[u]}\\
\sigma^2\xord\sigma^3
\ar[d]_{<p_1,p_2>=\Id[\sigma^2\xord\sigma^3]
=\Id[\sigma^2]\xord\Id[\sigma^3]}^{\brspace\bigg\} \Par[u]}\\
\sigma^2\xord\sigma^3
\ar[d]_h^{\brspace\bigg\} \Sep[u]}\\
\sigma^3}$$

Therefore
$\Arr[u,W,\sigma^3]=h<p_1,p_2><p_2,p_6>=h<p_2,p_6>$ and
$u=h( \ds x21, \ds x32)$ (over $\left(V\setminus\{ \ds
x32\}\right)\union W$) is the arrow
$$\xymatrix{
\mk{\sigma}^1\xord\mk{\sigma}^1\xord\mk{\sigma}^1
\xord\sigma^2\xord\mk{\sigma}^2\xord\mk{\sigma}^2
\xord\mk{\sigma}^3\xord\sigma^3\xord\mk{\sigma}^3
\xord\mk{\sigma}^4
\ar[d]_{<p_4,p_8>}^{\brspace\bigg\} \Dia[u]}\\
\sigma^2\xord\sigma^3
\ar[d]_{<p_1,p_2>=\Id[\sigma^2\xord\sigma^3]
=\Id[\sigma^2]\xord\Id[\sigma^3]}^{\brspace\bigg\}\Par[u]}\\
\sigma^2\xord\sigma^3
\ar[d]_h^{\brspace\bigg\}\Sep[u]}\\
\sigma^3}$$ so that $\Arr[u,\left(V\setminus\{ \ds
x32\}\right)\union W,\sigma^3]
=h<p_1,p_2><p_4,p_8>=h<p_4,p_8>$.

Now suppose that $\sigma^1$, $\sigma^2$, $\sigma^3$ and
$\sigma^4$ are inhabited.  Then by
Lemma~\ref{constantlemma}, we may choose constants
$c_i:1\to\sigma^i$ for $i\in 1\twodots4$. Then we have the
maps
\begin{gather*}
\alpha:\Type[W]\to\Type\bigl[\left(V\setminus\{ \ds
x32\}\right)\union W\bigr]\\
\beta:\Type\bigl[\left(V\setminus\{ \ds x32\}\right)\union
W\bigr]\to\Type[W] \end{gather*} This is $\alpha$:
{\taller\taller$$\xymatrix{
\Type[W]=\mk{\sigma}^1\xord\sigma^2\xord\mk{\sigma}^2
\xord\mk{\sigma}^2\xord\mk{\sigma}^3\xord\sigma^3
\xord\mk{\sigma}^3
 \ar[d]|{\alpha=< c_1!, c_1!, c_1!, p_2, c_2!, c_3!, c_3!,
p_6, c_3!, c_4!>}\\ \mk{\sigma}^1\xord\mk{\sigma}^1
\xord\mk{\sigma}^1\xord\sigma^2\xord\mk{\sigma}^2\xord\mk{\sigma}^2
\xord\mk{\sigma}^3\xord\sigma^3\xord\mk{\sigma}^3\xord\mk{\sigma}^4
}$$ where the codomain is
$$\Type\bigl[\left(V\setminus\{ \ds x32\}\right)\union W\bigr]$$
and $!$ is the unique map with domain $\Type[W]$ and
codomain $1$. The map $\beta $ is similarly defined. Both
are the maps given by Lemma~\ref{lemmb}.

It follows that $$<p_4,p_8>\alpha=<p_2,p_6>$$ and
$$<p_2,p_6>\beta=<p_4,p_8>$$ and that
$$\Arr[u,W,\sigma^3]=\Arr\bigl[u,\left(V\setminus\{ \ds
x32\}\right)\union W, \sigma^3\bigr]\oc\alpha$$ and
$$\Arr\bigl[u,\left(V\setminus\{ \ds x32\}\right)\union W,
\sigma^3\bigr]=\Arr[u,W,\sigma^3]\oc\beta$$

We now proceed with our example.  After substitution,
$$e[\ds x32 \leftarrow u]=f\left(\ds x11, \ds x43, h(
\ds x21, \ds x32), \ds x11, g( \ds x12, h( \ds x21, \ds
x32)), \ds x21\right)$$ This corresponds to the arrow
{\def\brspace{\hspace{7em}}
$$ \xymatrix{ \sigma^1 \xord \sigma^1 \xord \mk{\sigma}^1 \xord
\sigma^2 \xord \mk{\sigma}^2 \xord \mk{\sigma}^2 \xord
\mk{\sigma}^3 \xord \sigma^3 \xord \mk{\sigma}^3 \xord
\sigma^4 \ar[d]_{<p_1, p_{10}, p_4, p_8, p_1, p_2, p_4,
p_8, p_4>}^{\brspace\bigg\}\Dia[e]}\\ \sigma^1 \xord
\sigma^4 \xord \sigma^2 \xord \sigma^3 \xord \sigma^1 \xord
\sigma^1\xord\sigma^2 \xord \sigma^3 \xord \sigma^2
\ar[d]_{\left<p_1, p_2, <p_3, p_4>, p_5,\left<p_6, <p_7,
p_8>\right>, p_9\right>}^{\brspace\bigg\}\Par[e]}\\
\sigma^1 \xord \sigma^4 \xord (\sigma^2 \xord \sigma^3 )
\xord \sigma^1 \xord \left(\sigma^1 \xord (\sigma^2 \xord
\sigma^3)\right) \xord \sigma^2 \ar[d]|{\Id[\sigma^1] \xord
\Id[\sigma^4] \xord \Id[\sigma^2 \xord \sigma^3] \xord
\Id[\sigma^1] \xord \left(\Id[\sigma^1] \xord h\right)
\xord \Id[\sigma^2]\phantom{\left(\Id[\sigma^1] \xord
h\right) \xord \Id[\sigma^2]}}
\ar@{}[ddd]^{\brspace\hspace{1em}\Biggggg\} \Sep[e]}\\
\sigma^1 \xord \sigma^4 \xord (\sigma^2 \xord \sigma^3 )
\xord \sigma^1 \xord (\sigma^1 \xord \sigma^3) \xord
\sigma^2 \ar[d]_{\Id[\sigma^1] \xord \Id[\sigma^4] \xord h
\xord
\Id[\sigma^1] \xord g \xord \Id[\sigma^2]}\\
\sigma^1 \xord \sigma^4 \xord \sigma^3 \xord \sigma^1 \xord
\sigma^5 \xord \sigma^2
\ar[d]_f\\
\sigma^5 }$$ }

We now calculate
\begin{equation*}\begin{split}
e\left[x^3_2\leftarrow u\right] &= f\left(\Id[\sigma^1]
\xord \Id[\sigma^4] \xord h \xord
\Id[\sigma^1] \xord g \xord \Id[\sigma^3]\right)\\
& \quad\quad \oc \left( \Id[\sigma^1] \xord \Id[\sigma^4]
\xord \Id[\sigma^2 \xord \sigma^3] \xord \Id[\sigma^1]
\xord ( \Id[\sigma^1] \xord h)
\xord \Id[\sigma^2]\right) \\
& \quad\quad \oc < p_1, p_2, < p_3, p_4 >, p_5, < p_6, < p_7, p_8>>, p_9>\\
& \quad\quad \oc
< p_1, p_{10}, p_4, p_8, p_1, p_2, p_4, p_8, p_4>\\
&= f\bigl( \Id[\sigma^1] \xord \Id[\sigma^4] \xord
\Id[\sigma^3] \xord \Id[\sigma^1] \xord g \xord
\Id[\sigma^2]
\bigr)\\
& \quad\quad \oc
< p_1, p_2, h<p_3, p_4>, p_8, <p_6, h<p_7, p_8>>, p_9>\\
& \quad\quad \oc
<p_1, p_{10}, p_4, p_8, p_1, p_2, p_4, p_8, p_4>\\
&= f\bigl( \Id[\sigma^1] \xord \Id[\sigma^4] \xord
\Id[\sigma^3] \xord \Id[\sigma^1] \xord g \xord
\Id[\sigma^2]
\bigr)\\
& \quad\quad \oc
< p_1 , p_2 , p_3 , p_4 , < p_5 , p_6> , p_7>\\
& \quad\quad \oc
\left<p_1 , p_2, h<p_3, p_4>, p_5, p_6, h<p_7, p_8>, p_9\right>\\
& \quad\quad \oc
<p_1, p_{10}, p_4, p_8, p_1, p_2, p_4, p_8, p_4>\\
&= f\bigl( \Id[\sigma^1] \xord \Id[\sigma^4] \xord
\Id[\sigma^3] \xord \Id[\sigma^1] \xord g \xord
\Id[\sigma^2]
\bigr)\\
& \quad\quad \oc
\left<p_1, p_2, p_3, p_4, <p_5, p_6>, p_7\right>\\
& \quad\quad \oc
\left<p_1, p_{10}, h<p_4, p_8>, p_1, p_2, h<p_4, p_8>, p_4\right>\\
&= f\bigl( \Id[\sigma^1] \xord \Id[\sigma^4] \xord
\Id[\sigma^3] \xord \Id[\sigma^1] \xord g \xord
\Id[\sigma^2]
\bigr)\\
& \quad\quad \oc \left< p_1, p_2, p_3, p_4, <p_5, p_6>, p_7\right>\\
& \quad\quad \oc
<p_1, p_{10}, p_8, p_1, p_2, p_8, p_4>\\
& \quad\quad \oc <p_1, p_2, p_3, p_4, p_5, p_6, p_7, h<p_4,p_8>, p_9, p_{10}>\\
&=\Arr\left[e,V\union W, \Type[e]\right]\\
& \quad\quad \oc \left<p_1, p_2, p_3, p_4, p_5, p_6, p_7,
\Arr\left[u,\left(V\setminus\{x^3_2\}\right)\union W,
\Type[u]\right],p_9,p_{10}\right>
\end{split}\end{equation*}

\exam\label{sexam}

$$u=m( x^2_1, x^2_1, x^4_4)$$
$$\begin{cases} \Inp[u]= {\sigma}^2 \xord \sigma^2 \xord \sigma^4 &\\
\Outp[u] = \sigma^3 &\\
\end{cases}$$
$$W=\{ \mk{x}^1_1, x^2_1, \mk{x}^2_2, \mk{x}^2_3, \mk{x}^3_1, \mk{x}^3_3,
x^4_4\}$$ This is different from Example~\ref{fexam}
because the variable $x^3_2$ (in $e$) for which we are
making the substitution does not reappear in $u$.

This is the arrow for $u=m(x^2_1, x^2_1, x^4_4)$ over $W$:
$$\xymatrix{
\sigma^1 \xord \sigma^2 \xord \mk{\sigma}^2 \xord
\mk{\sigma}^2 \xord \mk{\sigma}^3 \xord \mk{\sigma}^3 \xord
\sigma^4
\ar[d]_{<p_1, p_1, p_7>}^{\brspace\bigg\} \Dia[u]}\\
\sigma^2 \xord \sigma^2 \xord \sigma^4 \ar[d]_{<p_1, p_2,
p_3> =\Id[\sigma^2 \xord \sigma^2 \xord \sigma^4]
= \Id[\sigma^1] \xord \Id[\sigma^2] \xord \Id[\sigma^4]}^{\brspace\bigg\} \Par[u]}\\
\sigma^2 \xord \sigma^2 \xord \sigma^4
\ar[d]_m^{\brspace\bigg\} \Sep[u]}\\
\sigma^5 }$$ We calculate
$$\left(V\setminus\{x^3_2\}\right)\union W=
\{x^1_1, x^1_2, \mk{x}^1_3, x^2_1, \mk{x}^2_2, \mk{x}^2_3,
\mk{x}^3_1, \mk{x}^3_3, x^4_3, x^4_4\}$$ Then $u=m(x^2_1,
x^2_1, x^4_4)$ (over $\bigl(V\setminus\{x^3_2\}\bigr)\union
W$) gives the arrow
$$\xymatrix{ \sigma^1 \xord \sigma^1 \xord \mk{\sigma}^1 \xord
\sigma^2 \xord \mk{\sigma}^2 \xord \mk{\sigma}^2 \xord
\mk{\sigma}^3 \xord \mk{\sigma}^3 \xord \sigma^4 \xord
\sigma^4 \ar[d]_{<p_4, p_4, p_{10}>}^{\brspace\bigg\}
\Dia[u]}\\ \sigma^2 \xord \sigma^2 \xord \sigma^4
\ar[d]_{<p_1, p_2, p_3> =\Id[\sigma^2 \xord \sigma^2 \xord
\sigma^4] = \Id[\sigma^1] \xord \Id[\sigma^2] \xord
\Id[\sigma^4]}^{\brspace\bigg\} \Par[u]}\\ \sigma^2 \xord
\sigma^2 \xord \sigma^4 \ar[d]_m^{\brspace\bigg\}
\Sep[u]}\\ \sigma^5 }$$

After substitution,
$$e[x^3_2\leftarrow u]=
f\bigl( x^1_1, x^4_3, m\left( x^2_1, x^2_1, x^4_4\right),
x^1_1, g\left(x^1_2, m\left(x^2_1, x^2_1,
x^4_4\right)\right), x^2_1\bigr)$$ This corresponds to the
arrow shown below.
{\taller\taller\def\brspace{\hspace{10em}}$$ \xymatrix{
\sigma^1 \xord \sigma^1 \xord \mk{\sigma}^1 \xord \sigma^2
\xord \mk{\sigma}^2 \xord \mk{\sigma}^2 \xord \mk{\sigma}^3
\xord \mk{\sigma}^3 \xord \sigma^4 \xord \sigma^4
\ar[d]|{\left<p_1, p_9, p_{10}, p_4, p_4, p_{10}, p_1, p_2,
p_4, p_4, p_{10}, p_2\right>}^{\brspace\bigg\} \Dia[e]}\\
\sigma^1 \xord \sigma^4 \xord \sigma^2 \xord \sigma^2 \xord
\sigma^4 \xord \sigma^1 \xord \sigma^1 \xord \sigma^2 \xord
\sigma^2 \xord \sigma^4 \xord \sigma^2 \ar[d]|{\bigl<p_1,
p_2, \left<p_3, p_4, p_5\right>, p_6, \left<p_7, <p_8, p_9,
p_{10}>\right>, p_{11}\bigr>}^{\brspace\bigg\} \Par[e]}\\
\sigma^1 \xord \sigma^4 \xord \left( \sigma^2 \xord
\sigma^2 \xord \sigma^2 \right) \xord \sigma^1 \xord \bigl(
\sigma^1 \left( \sigma^2 \xord \sigma^2 \xord
\sigma^4\right)\bigr) \xord \sigma^2 \ar[d]|{\Id[\sigma^1]
\xord \Id[\sigma^4] \xord \Id\left[\sigma^2 \xord \sigma^2
\xord \sigma^4\right] \xord \Id[\sigma^2] \left(
\Id[\sigma^1] \xord h \right) \xord \Id[\sigma^2]}\\
\sigma^1 \xord \sigma^4 \xord \left( \sigma^2 \xord
\sigma^2 \xord \sigma^2 \right) \xord \sigma^1 \xord \bigl(
\sigma^1 \xord \mk{\sigma}^3 \bigr) \xord \sigma^2
\ar[d]|{\Id[\sigma^1] \xord \Id[\sigma^4] \xord h \xord
\Id[\sigma^1] \xord g \xord
\Id[\sigma^2]}^{\brspace\biggggg\} \Sep[e]}\\ \sigma^1
\xord \sigma^4 \xord \mk{\sigma}^3 \xord \sigma^1 \xord
\mk{\sigma}^5 \xord \sigma^2 \ar[d]^f\\ \tau }$$}

We next re-express this in a convenient form:

\begin{equation*}\begin{split}
e[x^3_2\leftarrow u]&= f \oc \bigl( \Id[\sigma^1] \xord
\Id[\sigma^4] \xord h \xord
\Id[\sigma^2] \xord g \xord \Id[\sigma^2] \bigr)\\
 & \quad\quad \oc \bigl(
\sigma^1 \xord \sigma^4 \xord \Id\left[ \sigma^2 \xord
\sigma^2 \xord \sigma^4\right] \xord \Id[\sigma^1] \xord
\left(  \Id[\sigma^1] \xord h \right) \xord \Id[\sigma^2] \bigr)\\
& \quad\quad \oc \bigl< p_1, p_2, \left< p_3, p_4, p_5
\right>, p_6, \left< p_7, \left< p_8, p_9, p_{10} \right>
\right>,
p_{11}\bigr>\\
 & \quad\quad \oc \bigl< p_1, p_9, p_4, p_4, p_{10},
p_1, p_2, p_4, p_4, p_{10}, p_2 \bigr>\\
 &=f \oc \bigl(
\Id[\sigma^1] \xord \Id[\sigma^4] \xord \Id[\sigma^3] \xord
\Id[\sigma^1] \xord g \xord \Id[\sigma^2] \bigr)\\
 & \quad\quad \oc
\bigl< p_1, p_2, h<p_3, p_4, p_5>, p_6, \left< p_7, h<p_8,
p_9,
p_{10}> \right>, p_{11} \bigr>\\
 & \quad\quad \oc \bigl< p_1, p_9,
p_4, p_4, p_{10}, p_1, p_2, p_4, p_4, p_{10}, p_2 \bigr>\\
 &= f
\oc \bigl( \Id[\sigma^1] \xord \Id[\sigma^4] \xord
\Id[\sigma^3]
\xord \Id[\sigma^1] \xord g \xord \Id[\sigma^2] \bigr)\\
& \quad\quad \oc \bigl< p_1, p_9, h<p_4, p_4, p_{10}>, p_1,
\left<
p_2, h<p_4, p_4, p_{10}> \right>, p_2 \bigr> \\
 &= f \oc \bigl(
\Id[\sigma^1] \xord \Id[\sigma^4] \xord \Id[\sigma^3] \xord
\Id[\sigma^1] \xord g \xord \Id[\sigma^2] \bigr)\\
 & \quad\quad \oc
\bigl< p_1, p_2, p_3, p_4, \left<p_5, p_6 \right>, p_7 \bigr>\\
& \quad\quad \oc \bigl< p_1, p_{10}, p_8, p_1, p_2, p_8, p_4
\bigr>\\
 & \quad\quad\oc \bigl< p_1, p_2, p_3, p_4, p_5, p_6, p_7,
h<p_4, p_4, p_{10}>, p_9, p_1, p_{11} \bigr>\\
&= \Arr \left( e,
V \union W, \sigma^5 \right)\\
 & \quad\quad \oc < p_1, p_2, p_3,
p_4, p_5, p_6, p_7, \Arr \bigl(u, \left(
V\setminus\{x^3_2\} \right) \union W \bigr), p_9, p_{10},
p_{11} \bigr>
\end{split}\end{equation*}

\rem In Examples~\ref{fexam} and~\ref{sexam} we may define
a map
$$A(e,u):\prd\TypeList\bigl[\left(V\setminus\{x^3_2\}\right)\union W\bigr]\to
\prd\TypeList\left[V\union W\right]$$ as follows. Choose
$I\in 1\twodots\Length \left[V\union W\right]$ such that
$$\bigl(\TypeList[V\union W]\bigr)_I=x^3_2$$
We next define for all $ i\in1\twodots\Length\left[V\union
W\right]$
$$\left(A(e,u) \right)_i=\Proj[i]$$ and $$\left(A(e,u)
\right)_I =\Arr\bigl(u, \left(V\setminus\{x^3_2\} \right)
\union W,\sigma^5\bigr)$$ With this definition in the
previous two examples we have $$\Arr\bigl[
\left(e,V,\sigma^5 \right)\left[x^3_2\leftarrow \left(u, W,
\sigma^3 \right)\right]\bigr] =\Arr\left[e,V\union
W,\sigma^5\right]\oc A(e,u)$$ This is an example of the
construction in~\ref{dirdef}.\erem

\chapter{Logic of equational theories}
\section{A lemma}
The following lemma will be used later in our discussion of
the rules ``concretion'' and ``abstraction'' that have to
do with including extraneous variables in and excluding
them from the list of variables of some term. The proof may
best be understood by considering the examples in
Section~\ref{subssec}.


\thl\label{lemmb} Let $e$ be an expression of type $\tau$,
and let $V_1$ and $V_2$ be lists of variables such that
$\VarSet[e]\includedin V_1$ and $\VarSet[e]\includedin
V_2$. Let
$$t_1\ceq\left[e,V_1,\tau\right]$$ and
$$t_2\ceq\left[e,V_2,\tau\right]$$
Then there are arrows
$$\alpha_{12}:\prd\TypeList[t_1]\to\prd\TypeList[t_2]$$
and
$$\alpha_{21}:\prd\TypeList[t_2]\to\prd\TypeList[t_1]$$
for which $\Arr[t_1]=\Arr[t_2]\oc\alpha_{12}$ and
$\Arr[t_2]=\Arr[t_1]\oc\alpha_{21}$. \ethl

\pf Assume $\Type[x]=\sigma$. The arrow $\alpha_{12}$ is
defined by the following requirements: \blist{$\alpha$}
\item If $x\in V_2\setminus\VarSet[e]$ and $x$ is in the
$j$th place in $\Var[t_2]$ (so that
$\TypeList[t_2]_j=\sigma$), then the diagram {
\begin{equation} \xymatrix{ \prd\TypeList[t_1]
\ar[d]_{\alpha_{12}} \ar[r]^-{!} & 1 \ar[d]^{c_j}\\
\prd\TypeList[t_2] \ar[r]_-{\Proj[j]} & \sigma
}%
\end{equation}
}%
must commute, where $c_j$ is some constant of type $\sigma$
(there must be one by Lemma~\ref{constantlemma}).  Note
that it does not matter which constant of type $\sigma$ is
chosen.

\item If $x\in \VarSet[e]$ and $x$ is in the $i$th place in
$\Var[t_1]$ and in the $j$th place in $\Var[t_2]$, then the
diagram { \spreaddiagramcolumns{-2em}
\begin{equation}
\xymatrix{ \prd\TypeList[t_1] \ar[rr]^{\alpha_{12}}
\ar[dr]_{\Proj[i]} & & \prd\TypeList[t_2]
\ar[dl]^{\Proj[j]}\\
&\sigma\\
}%
\end{equation}
}%
must commute. \elist\epf

\section{Rules of inference of MSEL as
factorizations}\label{rinfasfac} In this section we show
how the rules of inference of multisorted equational logic
can be codified into our present system.  This is a
two-step process. First, we show that for each rule of
inference the pair of equal arrows corresponding to the
conclusion of the rule of inference can be constructed
using the rules of construction of graph-based logic from
the single arrow or the product of the equal pairs of
arrows that form the hypothesis of that rule of inference.
Next, we exhibit the construction as an actual
factorization, where the nodes and arrows appearing in the
positions corresponding to the various labels on the
diagram~\eqref{actfact} are the appropriate instances of
the hypothesis, claim, workspace and so on for the rule in
question.
\begin{equation}\label{actfact}\xymatrix{
& \hyp \ar[d]^{\claimcon}\\
\claim \ar[ur]^{\verif} \ar[r]_{\hypcon} &
\wksp}\end{equation} While some of these are done in detail
some others are not.  For our purposes, it is enough to
prove that a codification as an actual factorization in
$\SynCat\bigl[\FinProd,\Fa\bigr]$ (as defined in
Section~\ref{sksigg}) is possible. In general, this may be
done in more than one way. Symmetry and reflexivity are
treated separately. Transitivity, concretion, abstraction
and substitutivity are all treated in
Section~\ref{alltreat}, as they all are special instances
of Example~\ref{examthm}.

\subsection{Reflexivity}

The equational rule of inference is
$$\Frac{}{e=_Ve}$$
\astep Translated into the present context, as an instance
of the rule of construction REF, this is represented as
$$\Rule{
\xymatrix{ A \ar[r]^f &   B\\
}}{
\xymatrix{ A \rtwo ff & B\\
}}{REF}{}
$$
where \begin{align*}f&\ceq
\Arr\left[e,V,\Type[e]\right]\\
A&\ceq \Type[V]\\
B&\ceq \Type[e]
\end{align*}
This concludes the first step.

\bstep The corresponding actual factorization: {\wider
\begin{equation}\label{refleq}\xymatrix{
&(\arr\x\arr)^{<f,f>}\ar[d]^{\Proj[1]}\\
\arr^f\ar[r]_{\Id[\arr]} \ar[ur]^{\Delta}& \arr^f
}\end{equation}
}%
where $\arr$ is the object of arrows in the sketch for
categories.

Note that one can also use $\Proj[2]$ as $\claimcon$.
This factorization actually occurs in
$\SynCat[\Cat,\F]$ and is inherited by
$\SynCat[\FinProd,\F]$. A similar remark is true of
the constructions for symmetry and transitivity.

\subsection{Symmetry}

Although Chapter~\ref{gbl}, we not use a rule of
construction corresponding to symmetry, we shall record an
actual factorization for this to facilitate later
discussion (in this section) on proofs as actual
factorizations. The rule in equational deduction is
$$\Frac{e=_Ve'}{e'=_Ve}$$
We define \begin{align*}f&\ceq \Arr\left[e,V,\Type[e]\right]\\
f'&\ceq \Arr\left[e',V,\Type[e]\right]\end{align*}
then the actual factorization is as exhibited
below:\mpark{I just realized that I don't understand
why this says that if f=g in a model then g=f in a
model.  I am running out of time so I will think about
it later.} {\wider\wider\wider
$$\xymatrix{
&(\arr\x\arr)^{<g,f>}\ar[d]^{<\Proj[2],\Proj[1]>}\\
(\arr\x\arr)^{<f,g>} \ar[r]_{\Id[\arr]
\xord\Id[\arr]}\ar[ur]|{<\Proj[2],\Proj[1]>}&(\arr\x\arr)^{<f,g>}
}$$ }

\subsection{Transitivity}

The equational rule of inference is
$$\Frac{e=_Ve'\quad e'=_Ve''}{e=_Ve''}$$
For the first step we define
\begin{align*}f:D\rightarrow C&\ceq
\Arr\left[e,V,\Type[e]\right]\\ g:D\rightarrow C&\ceq
\Arr\left[e',V,\Type[e]\right]\\ h:D\rightarrow C&\ceq
\Arr\left[e'',V,\Type[e]\right]\end{align*} Note that
$f$, $g$, and $h$ have the same domain and the same
codomain as $e$, $e'$, and $e''$ have the same type
and as $V$ is the same in each of the terms exhibited
below:

$$\Rule{\xymatrix{ D \ar[r]<1ex>^f \ar[r]<-1ex>_g & C
& D \ar[r]<1ex>^g \ar[r]<-1ex>_h & C\\
}}{\xymatrix{ D \ar[r]<1ex>^f \ar[r]<-1ex>_h & C\\
}}{TRANS}{}$$

The corresponding actual factorization is provided in
Section~\ref{alltreat}.


\subsection{Concretion}

In this case the equational inference rule reads

Given a set $V$ of typed variables, $x\in V$ and an
equation $e=_Ve'$ such that $x\in
V\setminus(\VarSet[e]\union\VarSet[e'])$, and given
that $\Type[x]$ is inhabited,
$$\Frac{e=_V e'}{e=_{V\backslash \{x\}}e'}$$

We define $\tau\ceq \Type[e]=\Type[e']$ and
$\sigma\ceq \Type[x]$, and
\begin{align*}f:P\rightarrow \tau&\ceq \Arr\left[e,V,\tau\right]\\
f':P\rightarrow \tau&\ceq \Arr\left[e',V,\tau\right]\\
g:Q\rightarrow \tau&\ceq \Arr\left[e,V\backslash\{x\},\tau\right]\\
g':Q\rightarrow \tau&\ceq
\Arr\left[e',V\backslash\{x\},\tau\right]\end{align*}

Using Lemma~\ref{lemmb}, we may choose a map
$$h:\prd\left(V\setminus\{x\}\right)\to\prd V$$ such that
\begin{align*}
g=f\oc h\\
g'=f'\oc h
\end{align*}
Thus coded as arrows, the rule reads
$$\frac{f=f'}{f\oc h=f'\oc h}$$

\subsection{Abstraction}

In this case the equational rule of inference reads

Given a set of typed variables and
$x\in\Vbl[S]\setminus V$, $$\Frac{e=_V
e'}{e=_{V\union\{x\}}e'}$$

We define $\tau\ceq \Type[e]=\Type[e']$ and
$\sigma\ceq \Type[x]$, and
\begin{align*}f:P\rightarrow \tau&\ceq
\Arr\left[e,V,\tau\right]\\ f':P\rightarrow \tau&\ceq
\Arr\left[e',V,\tau\right]\\ g:Q\rightarrow \tau&\ceq
\Arr\left[e,V\cup\{x\},\tau\right]\\ g':Q\rightarrow
\tau&\ceq
\Arr\left[e',V\cup\{x\},\tau\right]\end{align*}

Using Lemma~\ref{lemmb}, we may choose a map
$$h:\prd\left(V\union \{x\}\right)\to V$$ such that
$g=f\oc h$ and $g'=f'\oc h$.  Thus coded as arrows the
rule reads $$\frac{f=f'}{f\oc h=f'\oc h}$$

\subsection{Substitutivity}

Given a set $V$ of typed variables, $x\in V$, and
expressions $u$ and $u'$ for which
$\Type[x]=\Type[u]=\Type[u']$ and
$\Type[e]=\Type[e']=\tau$,
$$\Frac{e=_Ve'\quad u=_Wu'}{e[x\leftarrow
u]=_{V\setminus\{x\}\union W}e'[x\leftarrow u']}$$

We already have the representations

\begin{align*}f&\ceq\Arr\left[e,V,\tau\right]=\Sep[e]\Par[e]\Dia\left[e,V,\tau\right] \\
f'&\ceq\Arr\left[e',V,\tau\right]=\Sep[e']\Par[e']\Dia[e',V,\tau]\\
g&\ceq\Arr\left[u,W,\Type[u]\right]=\Sep[u]\Par[u]\Dia\left[u,W,\Type[u]\right]\\
g'&\ceq\Arr\left[u',W,\Type[u]\right]=\Sep[u']\Par[u']\Dia\left[u',W,\Type[u]\right]\\
h&\ceq\Arr\bigl[e[x\leftarrow
u],\left(V\setminus\{x\}\right)\union W,\tau\bigr]\\
h'&\ceq\Arr\bigl[e'[x\leftarrow
u'],\left(V\setminus\{x\}\right)\union W,\tau\bigr]
\end{align*}

In view of Lemma~\ref{facil}, we may choose arrows
$A\ceq\Insert[u,x,V,W]$ and $A'\ceq\Insert[u',x,V,W]$
for which $h=f\oc A$ and $h'=f'\oc A'$. Note that the
assumptions $e=_Ve'$ and $u=_Wu'$ are equivalent to
assuming that $f=f'$ and $A=A'$. It follows that when
coded in terms of arrows, the rule reads

$$\frac{f=f' \quad\quad A=A'}{f\oc A=f'\oc A'}$$

\subsection{Transitivity, concretion, abstraction and
substitutivity as  actual factorizations}\label{alltreat}

Transitivity may be viewed as a special case of
Theorem~\ref{examthm} once equations are interpreted as
commutative diagrams as shown in Diagram~\eqref{twa}:
\begin{equation}\label{twa}\xymatrix{ D
\ar[d]_{\Id[D]} \ar[r]^h & C \ar[d]^{\Id[C]}\\ D
\ar[ur]|g \ar[r]_f &  C}\end{equation} The fact that
the two triangles commute means that $h=g$ and $g=f$.
That the outside square commutes means that $h=f$.

In view of Lemma~\ref{lemmb}, concretion and
abstraction can be seen to be special cases of the
following: For every pair of formally equal arrows
$f,f':D\to C$ and for every $h:E\to D$, $f\oc h$ and
$f'\oc h$ are formally equal.  This can also be
realized as a special case of the commutativity of
Diagram~(\ref{examthm})\mpar{fixed reference here and in next paragraph}, with choices as shown:
\begin{equation}\label{twb}\xymatrix{ E \ar[d]_h
\ar[r]^h & D \ar[d]^f\\ D \ar[ur]|{\Id[D]} \ar[r]_g &
C}\end{equation} Particular choices for $h$ yield
concretion and abstraction.

In substitutivity, in view of Lemma~\ref{lemmb}, we
have the following in terms of arrows:  For every pair
of formally equal arrows $f,f':D\to C$, and for every
pair of formally equal arrows $A,A':E\to D$, $f\oc A$
and $f'\oc A'$ are formally equal.  This is also a
special case of Diagram~(\ref{examthm}) as shown below:
\begin{equation}\label{twc}\xymatrix{
E \ar[d]_{A'} \ar[r]^A
& D \ar[d]^f\\
D \ar[ur]|{\Id[D]} \ar[r]_{f'} & C}\end{equation} On
the basis of the preceding discussion we conclude that
we may make choices for all nodes and arrows in the
diagram
\begin{equation*}\label{actfact2}\xymatrix{ & \hyp
\ar[d]^{\claimcon}\\ \claim \ar[ur]^{\verif}
\ar[r]_{\hypcon} & \wksp}\end{equation*} so that the
actual factorization in $\SynCat[\FinProd,\F]$ codes
transitivity, concretion, abstraction and
substitutivity respectively.

\section{Deductions as factorizations}\label{dedasfac2}
We now show that deductions in MSEL (Section~\ref{deducs})
correspond to actual factorizations in
$\SynCat[\FinProd,\Fa]$.

We first need two lemmas.

\thl\label{evlem} Given two actual factorizations in
any syntactic category
\begin{equation}\label{tochain}\begin{array}{cc}
\xymatrix{ & B \ar[d]^{c_1}\\ H \ar[ur]^{u_1}
\ar[r]_{h_1} & W_1} & \xymatrix{ & C \ar[d]^{c_2}\\ B
\ar[ur]^{u_2} \ar[r]_{h_2} & W_2}
\end{array}\end{equation} there is a node $W$ and
arrows $c$ and $h$ for which
\begin{equation}\label{chained}\xymatrix{ & C \ar[d]^{c}\\ H
\ar[ur]^{u_2\oc u_1} \ar[r]_{h} & W}\end{equation} is
an actual
factorization.\ethl 

\pf As every node in a syntactic category (or
$\SynCat[\E,\F]$) is the vertex of a limit cone over
some diagram in $\E$, we may choose \begin{align*}
\Delta_1 &=\BsDiag[W_1]\\ \Delta_2 &=\BsDiag[W_2]\\
\Delta_B &=\BsDiag[B]
\end{align*}
to get the following in the category of diagrams  of
$\E$:

$$ \xymatrix{
& \Delta_1 \ar[d]^{\alpha_1}\\
\Delta_2 \ar[r]_{\alpha_2} & \Delta_B }$$ where
$\alpha_1$ and $\alpha_2$ are the morphisms of
diagrams that give rise to $c_1$ and $h_2$. As the
category of diagrams in a category is small complete,
we may form the pullback as shown:
\begin{equation}\label{pbincd} \xymatrix{ \Delta
\ar[r]^{\beta_1} \ar[d]_{\beta_2} & \Delta_1
\ar[d]^{\alpha_1}\\ \Delta_2 \ar[r]_{\alpha_2} & \Delta_B
}\end{equation} Taking the limit over the diagrams
corresponding to the vertices in~\eqref{pbincd} and using
the lemmas in Chapter~\ref{limdiagchap}, we get the
following diagram in $\SynCat[\E,\F]$:
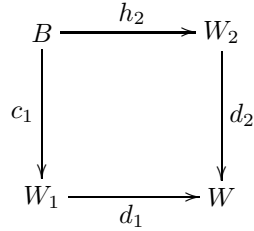
$$ \xymatrix{
B \ar[r]^{h_2} \ar[d]_{c_1}
& W_2 \ar[d]^{d_2}\\
W_1 \ar[r]_{d_1} & W }$$ This gives the following
diagram in $\SynCat[\E,\F]$:
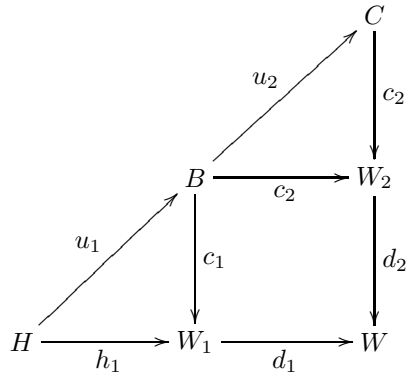
$$ \xymatrix{
&& C \ar[d]^{c_2}\\
& B \ar[ur]^{u_2} \ar[r]_{c_2} \ar[d]^{c_1}
& W_2 \ar[d]^{d_2}\\
H \ar[ur]^{u_1} \ar[r]_{h_1} & W_1 \ar[r]_{d_1} & W}$$
The lemma follows by setting $c \ceq d_2\oc c_2$ and $h
\ceq d_1\oc h_1$.\epf

\defn The Diagram~(\ref{chained}) is said
to be obtained by \textbf{chaining} the diagrams
in~(\ref{tochain}).\edefn

\rem Chaining is an analogue of getting the deduction
$\frac{E_1}{E}$ given the deduction $\frac{E_1}{E_2}$
and $\frac{E_2}{E}$.\erem

\thl\label{givlem} Given two actual factorizations in
any syntactic category
$$\begin{array}{cc}
\xymatrix{ & C_1
\ar[d]^{c_1}\\
H_1 \ar[ur]^{u_1} \ar[r]_{h_1} & W_1} & \xymatrix{
& C_2 \ar[d]^{c_2}\\
H_2 \ar[ur]^{u_2} \ar[r]_{h_2} & W_2}
\end{array}$$
we have the factorization
$$ \xymatrix{
& C_1\x C_2 \ar[d]^{c_1\x c_2}\\
H_1\x H_2 \ar[ur]^{u_1\x u_2} \ar[r]_{h_1\x h_2} &
W_1\x W_2 }$$\ethl \pf Omitted.\epf

\section{From deduction to factorization}

\todo{Restore the construction of a normal form for equational deductions that appeared in \cite{logstr2}.}
Now suppose we have a deduction $D$ of an equation $E$
from a list $P\ceq (E_1,\ldots,E_n)$ of equations.  We
show how to give a factorization corresponding to the
deduction for each of the four parts of
Definition~\ref{deddef} of deduction.

\paragraph*{D.1} $D=(E)$ and $P=(E)$.  Let $\n$ be the node of
$\SynCat[\FinProd,\F]$ corresponding to $E$.  Then the
factorization is
$$\xymatrix{
& \n \ar[d]^{\Id[\n]}
\\
\n \ar[r]_{\Id[\n]} \ar[ur]^{\Id[\n]} & \n
\\
}$$

\paragraph*{D.2} $D=(E)$ and $P$ is the empty list.  Then $E$ is
$e=_Ve$ and the factorization is given in
Diagram~(\ref{refleq}), where $f=\Arr[e,V,\Type[e]$.

\paragraph*{D.3}  $D=(E,D_1)$, where $D_1$ is a deduction of an
equation $E_1$ from $P$  and
$$\Frac{E_1}{E}$$ is an instance of a rule of inference $R$ of
MSEL.

In this case we have constructed the factorization for
$R$ in Section~\ref{rinfasfac}, and there is a
factorization for $D_1$ by induction hypothesis.
These may be chained (Lemma~\ref{evlem}) to get the
factorization for $D$.

\paragraph*{D.4}
$D=(E,D_1,D_2)$, where for $i=1,2$, $D_i$ is a
deduction of an equation $E_i$ from a list of premises
$P_i$, $P=P_1P_2$, and
$$\Frac{E_1\quad E_2}{E}$$ is an instance of a rule of inference
$R$ of MSEL.

By induction hypothesis, there are factorizations
$F_1$ and $F_2$ for $D_1$ and $D_2$, and we have
constructed a factorization $F$ for $R$ in
Section~\ref{rinfasfac}. $F_1$ and $F_2$ may be
combined into a single factorization by
Lemma~\ref{givlem}, and the resulting factorization
may be chained with $F$ to obtain a factorization for
$D$.\mpar{}
\section{From factorization to deduction}
\todo{Write this section.}

%
%

\chapter{Future work}\label{futwor}

\section{More explicit rules of construction} It is
noteworthy that the rules of construction for
constructor spaces given in Section~\ref{rulesapp}
correspond to arrows of  $\FinLim$, although not in a
one-to-one way (see Remark~\ref{rulerem}). The rules
are given here in a form that requires pattern
recognition (recall the discussion in~\ref{vfdisc}),
but they clearly could be given at another level of
abstraction as arrows or families of arrows of
$\FinLim$.  \mpark{Omitted: ``We expect to make this
explicit in a later work.''}

\section{Comparison with Makkai's work} M. Makkai
\shortcite{msdp}, \shortcite{maksyn} has produced an
approach to explicating the logic of sketches that is
quite different from that presented here. Both are
attempts at codifying the process of
diagram-manipulation used to prove results valid in
particular structured categories.  In both cases,
structured categories or doctrines form the semantic
universes with a pre-existing notion of validity. Both
approaches are motivated by the desire to formulate a
syntactic notion of deducibility. After that is said,
the two approaches are very different and a detailed
investigation of the relationship between them is
desirable.

The following points are however worth mentioning.
\begin{itemize}
\item Both approaches require a generalization of Ehresmann
sketches: Forms as in Definition~\ref{skedef} here
and sketch category
over a category $G$ in \cite{maksyn}.
\item In Makkai's approach, the sketch axioms are different for
different doctrines and serve both as axioms and as rules of
inference.  In contrast, here the rules of construction are the
{\it same\/} for all doctrines; what distinguishes doctrines is
their specification as CS-sketches.  This feature is a departure
from the usual practice in symbolic logic.\end{itemize}

\section{Equivalences with other logics} In
Chapter~\ref{EqThChap}, we worked out the details of
the equivalence of multisorted equational logic  and
$\FPTh\bigl[Sk[\Ssc]\bigr]$. The method used there
should\mpark{Replaced ``will'' by ``should''} work for
any logical system that can be described as a
constructor-space sketch. Thus, in general, we shall
have some logical system $L$ and a category
$\CatTh[\E_{L},\F]$  in which $\E_L$ is the kind of
category in which the models of $L$ are.  For
instance, if $L$ is the typed $\lambda$-calculus,
$E_L$ would be $\CCC$, and if $L$ is intuitionistic
type theory, then $E_L$ would be a constructor space
for toposes.

Given any sound and complete deductive system for $L$,
if we interpret terms as arrows and encode them in
$\CatTh[\E_L,\F]$ as we have done here, then we
conjecture that the method will show that all theorems
of $L$ can be realized as actual factorizations in
$\CatTh[\E_L,\F]$.  (Indeed it appears nearly obvious
that this will happen if we know that $L$ and $E_L$
have equivalent models; a detailed proof, is of course
necessary to settle the matter.) In the examples of
the preceding paragraphs, we might use the deductive
systems formulated in \cite{lambekscottbook}. The
method used here is quite general.

\section{Syntax by other doctrines} The rules of
construction given herein take place in the doctrine
of finite limits.  Most of the syntax and rules of
deduction of string-based logical theories are clearly
expressible using context-sensitive grammars, which
intuitively at least can be modeled using finite
limits. (Context-free grammars can be modeled using
only finite products \cite{wellsbarr}.) However, one
could imagine extensions of this doctrine: \romlist
\item Use coproducts to allow the specification of
conditional compilation or other syntactical
alternatives.  \item Use doctrines involving
epimorphic covering families to allow the intentional
description of ambiguous statements.  \item Use some
extension of finite limit doctrines to allow the
treatment of categories whose structure is determined
only up to isomorphism (the usual approach in
category theory) instead of being specified. There are
two approaches in the literature:  the use of
categories enriched over groupoids and universal
sketches as described in \cite{limsk}.  It is not
clear that such methods are necessary for applications
in computer science, but they may be for general
applications to mathematical reasoning.\eromlist

\nocite{powerwells}
\nocite{gensk}
\nocite{lambekscottbook}
\nocite{ctcs}
\nocite{limsk}
\nocite{bastehres}
\nocite{EFT}
\nocite{ehres68}
\nocite{fourvick}
\nocite{guitartlair80}
\nocite{gogburst}
\nocite{lambekscott84}
\nocite{lambekscottbook}
\nocite{msdp}
\nocite{maksyn}
\nocite{makkaipare}
\nocite{makkaireyes}
\nocite{GL}
\nocite{powerwells}
\nocite{gensk}
\nocite{sln236}
\nocite{mfps13}

\tolerance=400
\addcontentsline{toc}{chapter}{Bibliography}
\raggedright
\bibliography{gbls}
\bibliographystyle{xcw}

\bigskip\noindent
\begin{minipage}{220pt}
Atish Bagchi\\
Department of Mathematics\\
Community College of Philadelphia\\
1700 Spring Garden St.\\
Philadelphia, PA 19130, USA\\
{\tt atish@math.upenn.edu}
\end{minipage}

\vspace{12pt}
\begin{minipage}{220pt}
Charles Wells\\
Professor Emeritus of Mathematics\\
Case Western Reserve University\\
Current address:\\
1440 Randolph Avenue\\
Apt 323\\
Saint Paul, MN 55105, USA\\
{\tt charles@abstractmath.org}
\end{minipage}

\end{document}